\documentclass{amsart}
\usepackage[a4paper, total={6in, 8in}]{geometry} 	%Paper size
\usepackage[utf8]{inputenc}
\usepackage[english]{babel}
\usepackage{mypreamble}
\newcommand{\Crit}{\mathrm{Crit}}
\newcommand{\Emb}{\mathrm{Emb}}
\newcommand{\hcEmb}{\mathrm{hcEmb_*^+}}
\newcommand{\SOEmb}{\mathrm{SOEmb}}
\newcommand{\SPD}{\mathrm{SPD}}

\title{The Goresky-Hingston Coproduct in Morse Homology with DG Coefficients}
\author{Jonathan Clivio}

\begin{document}
\begin{abstract}
    We describe the Goresky-Hingston coproduct on the free loop space with real coefficients via the quasi-isomorphism $C_*(\Lambda M)\simeq C_*(M,C_*(\Omega M))$. This lets us describe the coproduct on the Leray-Serre spectral sequence as the diagonal on the manifold and the coproduct on the based loop space. With this description, we compute the coproduct with real coefficients for manifolds of the form $M=S^n/G$. The appendix contains results on inverting morphisms of $\mathcal{A}_\infty$-modules that are required for our computations.
\end{abstract}

\maketitle

\section{Introduction}
In this article, we give a description of the Goresky-Hingston coproduct in real coefficients
\begin{align*}
    \nu_{GH}\colon H_*(\Lambda M, M)\to H_{1-n+*}((\Lambda M, M)^2)
\end{align*}
via the framework of Morse homology with differential graded coefficients as defined in \cite{barraud2023morse}. Here, $\Lambda M$ denotes the free loop space of a smooth manifold $M$ of dimension $n$. It contains $M$ as the constant loops.

The Goresky-Hingston coproduct or string topology coproduct is a coproduct on $H_*(\Lambda M, M)$ of degree $1-n$, introduced by Goresky and Hingston in \cite{goresky2009loop} and Sullivan in \cite{sullivan2004open}.

Morse homology with DG coefficients gives a quasi-isomorphism
\begin{align*}
    C_*(M,C_*(\Omega M,x_0))\simeq C_*(\Lambda M,M),
\end{align*}
where $C_*(M,C_*(\Omega M,x_0))$ denotes Morse homology with DG coefficients in the $C_*(\Omega M)$-module $C_*(\Omega M,x_0)$. We use the extension of this formalism by Riegel in \cite{riegel2024chain}. He proves that $(\eta_k)_{k\geq 1}$, a morphism of $\mathcal{A}_\infty$-modules over $C_*(\Omega M)$ of the coefficients, induces $\widetilde{\eta}$, a morphism on Morse homology with DG coefficients. 

This allows us to prove the following:

\begin{IntroThm}[Theorem \ref{thm:model of GH}]\label{thm:A}
    Let $(M,x_0)$ be a smooth, pointed manifold of dimension $n$. There exists a morphism of $\mathcal{A}_\infty$-modules over $C_*(\Omega M)$
    \begin{align*}
        (\nu_k)_{k\geq 1}\colon  C_{*}(\Omega M, x_0;\R)\to C_{1-n+*}((\Omega M,x_0)^2;\R),
    \end{align*}
    where $\nu_1$ computes the Goresky-Hingston coproduct on the based loop space.
    
    Moreover for the induced map $\widetilde{\nu}$, it holds that
    \begin{align*}
        C_*(\Lambda M,M;\R)\simeq C_*(M,C_*(\Omega M,x_0;\R))&\overset{\widetilde{\nu}}{\to} C_*(M,C_{1-n+*}((\Omega M,x_0)^2;\R)) \\
        &\overset{\Delta_*}{\to}C_*(M\times M,C_{1-n+*}((\Omega M,x_0)^2;\R))\\
        &\simeq C_{1-n+*}((\Lambda M,M)^2;\R)
    \end{align*}
    computes the Goresky-Hingston coproduct on homology of the free loop space.
\end{IntroThm}

An immediate corollary of this theorem is a description of the Goresky-Hingston coproduct in the Leray-Serre.

\begin{IntroThm}[Corollary \ref{cor:coprod spec seq}]\label{thm:B}
    The Goresky-Hingston coproduct is computed on the second page of the Leray-Serre spectral sequence $E_{p,q}^2=H_p(M,H_q(\Omega M,x_0))$ with real coefficients as the following composition
    \begin{align*}
        H_p(M,H_q(\Omega M,x_0;\R))&\overset{\nu_1}{\to}H_p(M,H_{q+1-n}((\Omega M,x_0)^2;\R))\overset{\Delta_*}{\to} H_p(M\times M,H_{q+1-n}((\Omega M,x_0)^2;\R)).
    \end{align*}
    Dually, the dual Goresky-Hingston cohomology product is computed by first taking the cup product with local coefficients in $H_{q+1-n}((\Omega M,x_0)^2)$ and then computing the cohomology product on the based loop space.
\end{IntroThm}

The coproduct on the $E^2$-page lets us describe the full coproduct in certain cases. We study in Subsection \ref{subsec:Sn/G} the case $M=S^n/G$ for a finite group $G$ acting freely on $S^n$:

\begin{IntroThm}[Theorem \ref{thm:GH on Sn/G}]
    Let $G$ be a non-trivial, finite group acting freely and orientation-preserving on $S^n$ for $n>1$. The relative homology of the free loop space of $M=S^n/G$ with real coefficients is given by
    \begin{align*}
        H_*(\Lambda M;\R)\cong&\  \bigoplus_{[g]\in \mathrm{ccl}(G)}\bigoplus_{k\geq 0} \R x_{[g],k}\oplus \R y_{[g],k} 
    \end{align*}
    with $x_{[g],k}$ in degree $k(n-1)$ and $y_{[g],k}$ in degree $k(n-1)+n$. The homology $H_*(\Lambda M,M;\R)$ has the same description except missing the $[g]=[1]$ and $k=0$ summand.
    
    The Goresky-Hingston coproduct lifts to a map 
    \begin{align*}
        \widehat\lor\colon H_*(\Lambda M;\R)&\to H_*(\Lambda M;\R)\otimes H_*(\Lambda M;\R)
    \end{align*}
    computed by
    \begin{align*}
        \widehat\lor(x_{[g],k})=\sum_{i+j=k-1} \sum_{h\in G} x_{[gh^{-1}],i}\otimes x_{[h],j}
    \end{align*}
    and
    \begin{align*}
        \widehat\lor(y_{[g],k})=\sum_{i+j=k-1} \sum_{h\in G} (x_{[gh^{-1}],i}\otimes y_{[h],j}+y_{[gh^{-1}],i}\otimes x_{[h],j}).
    \end{align*}
\end{IntroThm}

However in general, the coproduct on the $E^2$-page does not capture all information. For example for $\Sigma_g$ a surface of genus $g\geq 2$, the coproduct on the $E^2$-page is trivial but not on $H_*(\Lambda \Sigma_g,\Sigma_g)$ (see Example \ref{exmp:g>1 surface}).

Moreover on the $E^2$-page of the Leray-Serre spectral sequence, we have
\begin{align*}
    \nu^N\circ \Lambda f=(\Lambda f\times \Lambda f)\circ \nu^M
\end{align*}
for $f$ a degree $1$ map (see Corollary \ref{cor:commute with deg1}). This is in contrast to the Goresky-Hingston coproduct with integer coefficients which does not commute in general with maps that are not simple homotopies (see \cite{naef2021string,naef2024simple,wahl2019invariance,kenigsberg2024obstructions}). It is an open question where the loss of information appears: while passing from integer coefficients to real coefficients or while passing to maps of spectral sequences (see Remark \ref{rem:commute with maps}).

\subsection{The Proof Strategy}
The main technical goal is to exhibit the Goresky-Hingston coproduct as a morphism of $\mathcal{A}_\infty$-modules over $C_*(\Omega M)$. Morally, we want to prove $\Omega M$-equivariance of the coproduct in a coherent way.\\

The Goresky-Hingston coproduct on a class $\alpha\in H_*(\Omega , x_0)$ is computed as follows:
\begin{enumerate}[1.]
    \item\label{step:times I} multiply with the unit interval $I:=[0,1]$ and obtain a class $I\times \alpha\in H_*(I\times \Omega, I\times \{x_0\}\cup \partial I\times \Omega)$;
    \item\label{step:intersect} \enquote{intersect} $I\times \alpha$ with the subspace of loops with self-intersections:
    \begin{align*}
        \mathcal{F}\Omega M:=\{(t,\gamma)\in I\times \Omega M\mid \gamma(t)=\gamma(0)\};
    \end{align*}
    \item\label{step:cut} cut using the map $\mathcal{F}\Omega M\to \Omega M \times \Omega M, (t,\gamma)\mapsto (\gamma|_{[0,t]},\gamma_{[t,1]})$.
\end{enumerate}
While steps \eqref{step:times I} and \eqref{step:cut} are straight-forward to make equivariant, step \eqref{step:intersect} needs more work.

We follow the description in \cite{hingston2017product}. In this setup, \enquote{intersecting} as in \eqref{step:intersect} is done by a cap product
\begin{align*}
    e_\Omega^*\tau_{x_0}\cap \colon C_*(e_\Omega^{-1}(U_{x_0}), I\times \{x_0\}\cup \partial I\times \Omega \cup e_\Omega^{-1}(V_{x_0}^c))\to C_{-n+*}(e^{-1}_\Omega(U_{x_0}),I\times \{x_0\}\cup \partial I\times \Omega M).
\end{align*}
Here $x_0\in V_{x_0}\subsetneq U_{x_0}\subseteq M$ denotes open balls around $x_0$, $e_\Omega\colon I\times \Omega\to M$ denotes the evaluation at time $t$ and $\tau_{x_0}\in C^n(U_{x_0},V_{x_0}^c)$ is a local orientation at $x_0$. We explain this in further detail in Subsection \ref{subsec:GH coprod}.

The guiding question for our article thus is the following:
\begin{Ques}\label{ques:Ainfty for cap}
    Does there exist a morphism of $\mathcal{A}_\infty$-modules
    \begin{align}\label{eq:guiding question}
        (q_k)_{k\geq 1}\colon  C_*(e_\Omega^{-1}(U_{x_0}), I\times \{x_0\}\cup \partial I\times \Omega M\cup e_\Omega^{-1}(V_{x_0}^c))\to C_{-n+*}(e_\Omega^{-1}(U_{x_0}), I\times \{x_0\}\cup \partial I\times \Omega  M)
    \end{align}
    over $C_*(\Omega M)$ with $q_1= e_\Omega^*\tau_\Omega\cap$?
\end{Ques}

To make sense of this, we need to specify a $\Omega M$-action on all the spaces appearing in \eqref{eq:guiding question}. For this, we construct a $\Omega M$-action on itself that is homotopic to conjugation such that, when we extend it trivially to an action on $I\times \Omega M$, it sends the subspaces $e_\Omega^{-1}(U_{x_0})$ and $e_\Omega^{-1}(V_{x_0}^c)$ to themselves.

To construct such an action we define the notion of a \textit{transitive diffeomorphism lift} in Section \ref{sec:trans diff lifts}. The idea is that we define for any path $\lambda\colon [0,a]\to M$ a diffeomorphism
\begin{align*}
    \Theta(\lambda)\colon M\cong M
\end{align*}
that 
\begin{enumerate}[(i)]
    \item sends $\lambda(0)$ to $\lambda(a)$;
    \item works well with composition of paths;
    \item sends $U_{x_0}$ and $V_{x_0}^c$ to themselves if $\lambda\in \Omega M$.
\end{enumerate}
This defines
\begin{align*}
    \Phi\colon \Omega M\times \Omega M&\to \Omega M,\\
    (\gamma,\lambda)&\mapsto\Theta(\lambda)_*(\gamma),
\end{align*}
the desired $\Omega M$-action on itself making Question \ref{ques:Ainfty for cap} well-posed.\\

Producing a morphism of $\mathcal{A}_\infty$-modules as in Question \ref{ques:Ainfty for cap} comes down to finding (higher) homotopies between 
\begin{align*}
    e_\Omega^*\tau_{x_0}\cap (\Phi(\alpha\times \lambda))
\end{align*}
and 
\begin{align*}
    \Phi((e_\Omega^*\tau_{x_0}\cap \alpha)\times \lambda).
\end{align*}
The transitive diffeomorphism lift $\Theta$ defines a right $\Omega M$-action 
\begin{align*}
    \rho\colon U_{x_0}\times \Omega M&\to U_{x_0},\\
    (x,\lambda)&\mapsto\Theta(\lambda)(x).
\end{align*}
We show in Lemma \ref{lem:rewritten q1 Phi} \eqref{item:rewritten a} that
\begin{align*}
    e_\Omega^*\tau_{x_0}\cap (\Phi(\alpha\times \lambda))=((\tau_{x_0}\circ \rho(e_\Omega\times \id))\otimes \Phi)\circ\Delta_*(\alpha\times \lambda).
\end{align*}
% If $\tau_{x_0}$ were invariant under $\rho$, we would in particular find
% \begin{align*}
%     e_\Omega^*\tau_{x_0} \cap \Phi(\alpha\times \lambda)=\Phi((e_\Omega^*\tau_{x_0}\cap \alpha)\times \lambda)
% \end{align*}
% for all $\alpha\in C_*(e_\Omega^{-1}(U_{x_0}), I\times M\cup \partial I\times \Omega M\cup e_\Omega^{-1}(V_{x_0}^c))$. However in general, it is not possible to choose a $\Theta$ and $\tau_{x_0}$ such that $\tau_{x_0}$ is invariant under $\rho$.

In Section \ref{sec:inv Thom class}, we construct a local orientation in $\tau_{x_0}\in C^n(U_{x_0},V_{x_0}^c;\R)\cong C^n(\R^n,B_1(0)^c)$ that is invariant under the $\SO(n)$-action. We can homotope $\rho$ to an action $\rho'$ that factors through $\SO(n)$. 

In Lemma \ref{lem:rewritten q1 Phi} \eqref{item:rewritten b}, we show that for such $\tau_{x_0}$ and $\rho'$ it holds
\begin{align*}
    ((\tau_{x_0}\circ \rho'(e_\Omega\times \id))\otimes \Phi)\circ \Delta_*(\alpha\times \lambda)=\Phi((e_\Omega^*\tau_{x_0}\cap \alpha)\times \lambda).
\end{align*}
Therefore answering Question \ref{ques:Ainfty for cap} comes down to finding (higher) homotopies between 
\begin{align*}
    ((\tau_{x_0}\circ \rho(e_\Omega\times \id))\otimes \Phi)\circ \Delta_*(\alpha\times \lambda)
\end{align*}
and
\begin{align*}
    ((\tau_{x_0}\circ \rho'(e_\Omega\times \id))\otimes \Phi)\circ \Delta_*(\alpha\times \lambda).
\end{align*}
We can do this by constructing explicit (higher) homotopies between $\rho$ and $\rho'$.

This in particular gives a positive answer to Question \ref{ques:Ainfty for cap} and lets us exhibit the Goresky-Hingston coproduct as a morphism of $\mathcal{A}_\infty$-morphism:
\begin{align*}
    \boldsymbol{\nu}\colon  C_{*}(\Omega M, x_0;\R)\to C_{1-n+*}((\Omega M,x_0)^2;\R).
\end{align*}
Using Riegel's framework, we get the map
\begin{align*}
    \widetilde{\nu}\colon C_*(M,C_*(\Omega M,x_0;\R))\to C_*(M,C_{1-n+*}((\Omega M, x_0)^2;\R)).
\end{align*}

To prove that this map indeed exhibits the Goresky-Hingston coproduct on the free loop space, we introduce \textit{coherent chain homotopies}. This is the algebraic part isolated from Riegel's work \cite{riegel2024chain} that lets us show that a map on the total space of a fibration is exhibited by a map in Morse homology with DG coefficients. We define and study these coherent chain homotopies in Section \ref{sec:cch}.

An analogous construction to the morphism of $\mathcal{A}_\infty$-modules $\boldsymbol{\nu}$ gives a coherent chain homotopy for the Goresky-Hingston coproduct. This lets us prove our main theorem in Section \ref{sec:mainthm}.

\addtocontents{toc}{\SkipTocEntry}
\subsection*{Organisation of the Paper.}
In Section \ref{sec:background}, we recall the necessary background on Morse homology with DG coefficients, fix the definition of the Goresky-Hingston coproduct of this article and specify our sign conventions. In Section \ref{sec:trans diff lifts}, we define and construct transitive diffeomorphism lifts. In Section \ref{sec:ainfty struc}, we exhibit the Goresky-Hingston coproduct as a morphism of $\mathcal{A}_\infty$-modules. In Section \ref{sec:cch}, we define and study coherent chain homotopies. In Section \ref{sec:ortho hc emb}, we define orthogonal half-constant embeddings and explain how they play an analogous role to $\SO(n)$ for constructing a coherent chain homotopy. In Section \ref{sec:inv Thom class}, we construct a Thom class which is invariant under certain maps using the Mathai-Quillen construction. In Section \ref{sec:mainthm}, we construct a coherent chain homotopy for the Goresky-Hingston coproduct and using this prove the main theorem. Finally, the appendix contains a section about inverting morphisms of $\mathcal{A}_\infty$-modules and a section about homotopies of pointed embeddings of $\R^n$ relative to $B_1(0)^c$.

\addtocontents{toc}{\SkipTocEntry}
\subsection*{Conventions.}
Throughout this article $M$ denotes a smooth, closed, oriented manifold of dimension $n$ with a Riemannian metric. The Riemannian metric induces the Riemannian distance, a metric on $M$. We denote it by $d_M\colon M\times M\to \R$. We fix a basepoint $x_0\in M$.

For $A,B$ smooth manifolds, we denote $\Emb(A,B)$ the space of smooth embeddings.

The free loop space on $M$ is defined as
\begin{align*}
    \Lambda:=\Lambda M:=\{\gamma\colon [0,a]\to M\mid a\geq 0,\gamma(0)=\gamma(a), \gamma \text{ is piecewise smooth}\}
\end{align*}
and the based loop space is defined as
\begin{align*}
    \Omega:=\Omega M:=\{\gamma\colon [0,a]\to M\mid a\geq 0, \gamma(0)=\gamma(a)=x_0,  \gamma \text{ is piecewise smooth}\}.
\end{align*}
We denote $\mu\colon \Omega M\times \Omega M\to \Omega M$ for the composition.

The path space on a smooth manifold $X$ is denoted by
\begin{align*}
    \mathcal{P}X:=\{\gamma\colon [0,a]\to X\mid \gamma \text{ is piecewise smooth}\}
\end{align*}
and for subspaces $A,B\subseteq X$, we denote
\begin{align*}
    \mathcal{P}_{A\to B}X:=\{(\gamma\colon [0,a]\to X)\in \mathcal{P}X\mid \gamma(0)\in A,\gamma(a)\in B\}.
\end{align*}
In this article, we work with cubical chains. A cubical chain $\sigma\colon I^k\to X$ is degenerate if there exists $1\leq i\leq k$ such that
\begin{align*}
    \sigma(v_1,\dots,v_i,\dots,v_k)=\sigma(v_1,\dots,v_i',\dots,v_k)
\end{align*}
for all $v_1,\dots, v_k,v_i'\in I$.

We denote by $C_k(X)$ the real cubical chains on a topological space. It is the real vector space generated by continuous maps $\sigma\colon I^k\to X$ relative to the subspace generated by all degenerate chains.

\addtocontents{toc}{\SkipTocEntry}
\subsection*{Acknowledgements.}
I thank Robin Riegel for extensive discussions figuring out how to adapt his framework to the present paper. I also would like to thank my advisor Nathalie Wahl for suggesting this topic and her encouragement and guidance. I further want to thank Alexandru Oancea for helpful discussions and providing new perspectives when I was stuck. I also want to thank him and the whole department at the University of Strasbourg for hosting me during the spring of 2025 where much of this work was completed. Furthermore, I want to thank Coline Emprin for being my correspondent to the $\mathcal{A}_\infty$-world. Finally, I would like to thank Isaac Moselle for helpful comments on the introduction. I was supported by the Danish National Research Foundation through the Copenhagen Centre for Geometry and Topology (DNRF151). 

\tableofcontents
\section{Background}\label{sec:background}

In this section, we recall the machinery of Morse homology with DG coefficients as introduced in \cite{barraud2023morse}. Moreover, we recall the definition of the Goresky-Hingston coproduct on the free loop space and based loop space of a manifold as first introduced in \cite{goresky2009loop} and \cite{sullivan2004open}. Finally in this section, we fix the sign conventions we use in this article for maps of chain complexes which do not preserve degree.

\subsection{Morse Homology with DG Coefficients}\label{subsec:Morse homo with DG}

In \cite{barraud2023morse}, Barraud, Damian, Humili\`{e}re and Oancea describe the framework of Morse homology with differential graded coefficients. We recall in this subsection the parts we need in this article. In particular for a fibration
\begin{center}
    \begin{tikzcd}
        F\ar[r] & E\ar[d, "\pi"]\\
        & M,
    \end{tikzcd}
\end{center}
we construct a complex $C_*(M,C_*(F))$ with a quasi-isomorphism
\begin{align*}
    \Psi\colon C_*(M,C_*(F))\to C_*(E).
\end{align*}
\subsubsection{Data and Choices}
For this, we choose the following data:
\begin{itemize}
    \item a basepoint $x_0\in M$;
    \item a Morse-Smale pair $(f,\xi)$ on $M$;
    \item an orientation on the unstable manifolds $W^u(x)$ of $(f,\xi)$;
    \item a tree $T\subseteq M$ containing $x_0$ and the critical points of $f$;
    \item a homotopy inverse $\theta$ to the projection $p\colon M\to M/T$;
    \item a representing chain system $\{s_{x,y}\}$ (see Definition \ref{def:rep chain sys});
    \item a compatible representing system $\{s_x\}$ (see Definition \ref{def:comp rep sys});
    \item a transitive lifting function $\Phi\colon E\,{}_\pi\!\times_{e_0}\mathcal{P}M\to E$ (see Definition \ref{def:trans lift fun}).
\end{itemize}

 For any pair of critical points $x,y \in \Crit(f)$, denote the space of Morse trajectories of $\xi$ between $x$ and $y$ by $\mathcal{L}(x,y)$ and define the \textit{moduli space of broken Morse trajectories}
 \begin{align*}
     \overline{\mathcal{L}}(x,y):= {\mathcal{L}}(x,y)\cup \bigcup_{z_1,\dots,z_k\in \Crit(f)} \mathcal{L}(x,z_1)\times\mathcal{L}(z_1,z_2)\times \dots \times \mathcal{L}(z_k,y).
 \end{align*}

The space $\overline{\mathcal{L}}(x,y)$ admits the structure of a manifold with corners of dimension $|x|-|y|-1$ with interior $\mathcal{L}(x,y)$ (see e.g.\ \cite[Proposition 2.15]{latour1994existence}). Moreover, the orientation on the unstable manifolds induce an orientation $[\overline{\mathcal{L}}(x,y)]\in H_{|x|-|y|-1}(\overline{\mathcal{L}}(x,y),\partial\overline{\mathcal{L}}(x,y))$.

\begin{Def}\label{def:rep chain sys}
    A \textit{representing chain system} of the moduli space of broken Morse trajectories is a collection $\{s_{x,y}\in C_{|x|-|y|-1}(\overline{\mathcal{L}}(x,y))\mid x,y\in \Crit(f)\}$ such that
    \begin{enumerate}[(i)]
        \item $s_{x,y}$ represents a cycle in $C_{|x|-|y|-1}(\overline{\mathcal{L}}(x,y),\partial \overline{\mathcal{L}}(x,y))$ and $[s_{x,y}]=[\overline{\mathcal{L}}(x,y)]$ in homology;
        \item the following equation holds
        \begin{align*}
            d s_{x,y}=\sum_{z\in \Crit(f)} (-1)^{|x|-|z|}s_{x,z}\times s_{z,y}.
        \end{align*}
    \end{enumerate}
\end{Def}
A representing chain system always exists \cite[Proposition 5.2.6]{barraud2023morse}.

Following \cite[Lemma 5.2.10]{barraud2023morse} we denote $\Gamma\colon \overline{\mathcal{L}}(x,y)\to \mathcal{P}_{x\to y}M$ the map that sends a trajectory to the path from $x$ to $y$ parametrised by $f$. We note that this is well-defined as this parametrisation produces smooth paths. We set
\begin{align*}
    m_{x,y}'&:=p_*\circ \Gamma_*(s_{x,y})\in C_{|x|-|y|-1}(\Omega (M/T)),\\
    m_{x,y}&:=\theta_*(m_{x,y}')\in C_{|x|-|y|-1}(\Omega M).
\end{align*}
By \cite[Definition 5.2.11]{barraud2023morse}, the collection $\{m_{x,y}\mid x,y\in \Crit(f)\}$ is a \textit{twisting cocycle} meaning that they satisfy the equation
\begin{align*}
    d (m_{x,y})=\sum_{z\in \Crit(f)}(-1)^{|x|-|z|}\mu(m_{x,z}\otimes m_{z,y}).
\end{align*}

\begin{Def}
    For $x\in \Crit(f)$, the \textit{Latour cell} is defined as
    \begin{align*}
        \overline{W^u}(x):=W^u(x)\cup \bigcup_{y\in \Crit(f)} \overline{\mathcal{L}}(x,y)\times W^u(y).
    \end{align*}
\end{Def}
The natural topology on the Latour cell makes it into a compact manifold with corners which is homeomorphic to the disk $D^{|x|}$ (see \cite[Proposition 2.11]{latour1994existence}).

\begin{Def}\label{def:comp rep sys}
    A \textit{compatible representing chain system} is a collection $\{s_x\in C_{|x|}(\overline{W^u}(x))\mid x\in \Crit(f)\}$ such that
    \begin{enumerate}[(i)]
        \item $s_{x}$ represents a cycle in $C_{|x|}(\overline{W^u}(x),\partial \overline{W^u}(x))$ and $[s_{x}]=[\overline{W^u}(x)]$;
        \item the following equation holds
        \begin{align*}
            d s_{x}=\sum_{x\in \Crit(f)}s_{x}\times s_{x,y}.
        \end{align*}
    \end{enumerate}
\end{Def}
A compatible representing chain system always exists \cite[Lemma 7.3.2]{barraud2023morse}. We denote $m_x':=p_*\circ \Gamma_*(s_x)\in C_{|x|}(\Omega M/T)$.

\begin{Def}\label{def:trans lift fun}
    A \textit{transitive lifting function} is a map
    \begin{align*}
        \Phi\colon E \,{}_\pi\!\times_{e_0}\mathcal{P}M\to E
    \end{align*}
    such that
    \begin{enumerate}[(i)]
        \item $\pi\circ \Phi=e_1\circ \mathrm{pr}_2$;
        \item for all $e\in E$, we have $\Phi(e,c_{\pi(e)})=e$, where $c_{\pi(e)}\in \mathcal{P}M$ is the constant path at $\pi(e)\in M$;
        \item for all $e\in E$ and composable $\gamma,\delta\in \mathcal{P} M$ with $\gamma(0)=\pi(e)$, we have
        \begin{align*}
            \Phi(\Phi(e,\gamma),\delta)=\Phi(e,\gamma\cdot \delta).
        \end{align*}
    \end{enumerate}
\end{Def}

Any fibration is fibre homotopy equivalent to a fibration that admits a transitive lifting function \cite[Proposition 5.5]{dyer1969some}.

A transitive lifting function in particular defines a right $\Omega M$-action on $F$.

\subsubsection{Construction}

We now describe the complex $C_*(M,C_*(F))$ and the map $\Psi\colon C_*(M,C_*(F))\to C_*(E)$: the complex $C_*(M,C_*(F))$ has underlying graded abelian group
\begin{align*}
    C_*(F)\otimes \Z\Crit(f)
\end{align*}
with differential given by
\begin{align*}
    d_F(\alpha\otimes x)&=d \alpha+(-1)^{|\alpha|}\sum_{y\in \Crit(f)} \Phi_*(\alpha\otimes m_{x,y})\otimes y.
\end{align*}
We denote the degree $-1$ map
\begin{align*}
    \boldsymbol{m}\colon \Z\Crit(f)&\to C_*(\Omega M)\otimes \Z\Crit(f),\\
    x&\mapsto \sum_y m_{x,y}\otimes y.
\end{align*}
Then the differential can be written as
\begin{align*}
    d_F=d\otimes \id +(\Phi\otimes \id)\circ (\id\otimes \boldsymbol{m})
\end{align*}
and the above sign $(-1)^{|\alpha|}$ comes from a Koszul sign.\\

For the map $\Psi\colon C_*(M,C_*(F))\to C_*(E)$, we consider the fibration $\pi'\colon E'\to M/T$ given by the pullback of $\pi\colon E\to M$ along $\theta\colon M/T\to M$. We can assume that $\theta$ is a homeomorphism. Thus $\Phi$ defines also a transitive lifting function on $\pi'\colon E'\to M/T$. This also induces a construction for $C_*(M/T,C_*(F))$. The differential on $C_*(M/T,C_*(F))$ is given by the action on $C_*(F)$ by $\Phi(\alpha\otimes m_{x,y}'):=\Phi(\alpha\otimes \theta_*(m_{x,y}'))=\Phi(\alpha\otimes  m_{x,y})$.

The map $\Psi$ is given by the composition
\begin{align*}
    C_*(M,C_*(F))\to C_*(M/T,C_*(F))\overset{\Psi'}{\to} C_*(E')\to C_*(E).
\end{align*}
The first map is an isomorphism. The third map above is given the homotopy equivalence $\theta$ and thus is a quasi-isomorphism. Thus the only map of interest is the middle map $\Psi'$. This is given by
\begin{align*}
    \Psi'(\alpha\otimes x)=\Phi_*(\alpha\otimes m'_x).
\end{align*}

Barraud, Damian, Humili\`{e}re and Oancea prove that $\Psi'$ is a quasi-isomorphism and therefore so is $\Psi$ \cite[Theorem 7.2.1]{barraud2023morse}.\\

We use an immediate corollary of this statement for relative fibrations:

\begin{Cor}
    Let $\pi\colon E'\subseteq E\to M$ be fibrations with fibres $F'\subseteq F$ with $\Phi$ a transitive lifting function of $\pi\colon E\to M$ that restricts to a transitive lifting functions on $E'$. The map
    \begin{align*}
        \Psi\colon C_*(M,C_*(F,F'))\to C_*(E,E')
    \end{align*}
    defines a quasi-isomorphism.
\end{Cor}
\begin{proof}
    The fact that $\Phi$ is a transitive lifting function for $E$ and $E'$ means that the inclusion
    \begin{align*}
        C_*(M,C_*(F'))\to C_*(M,C_*(F))
    \end{align*}
    is a chain map. Moreover, the following diagram commutes
    \begin{center}
        \begin{tikzcd}
            C_*(M,C_*(F'))\ar[r] \ar[d, "\Psi_{E'}"] & C_*(M,C_*(F)) \ar[d, "\Psi_E"]\\
            C_*(E')\ar[r] &C_*(E).
        \end{tikzcd}
    \end{center}
    The five lemma now concludes the proof.
\end{proof}

\subsubsection{Maps}

In this article, we need two different types of maps one can define in a Morse complex with DG coefficients.

In \cite[Section 9]{barraud2023morse}, a direct map is defined for a smooth map of based manifolds $\varphi \colon (M,x_0)\to (N,y_0)$ and a fibration $\pi\colon E\to N$
\begin{align*}
    C_*(M,\varphi^*C_*(F))\to C_*(N,C_*(F)).
\end{align*}
On the other hand, Riegel developed in \cite{riegel2024chain} a theory to construct maps on the Morse complex with DG coefficients for maps $C_*(F)\to C_*(F')$ which are morphisms of $\mathcal{A}_\infty$-modules. We are interested in $\mathcal{A}_\infty$-morphisms between strict $C_*(\Omega M)$-modules. We recall that $\Omega M$ denote Moore loops and thus $C_*(\Omega M)$ is indeed strictly associative.

\begin{Def}\label{def:morph of Ainfty modules}
    A \textit{morphism of $\mathcal{A}_\infty$-modules} $\boldsymbol{\eta}\colon \mathcal{F}\to \mathcal{G}$ where $(\mathcal{F}_*,d_\mathcal{F})$ and $(\mathcal{G},d_\mathcal{G})$ are strict right $C_*(\Omega M)$-modules with action $\mu_\mathcal{F}$ and $\mu_\mathcal{G}$ is given by a collection of maps for $k\geq 1$
    \begin{align*}
        \eta_k\colon \mathcal{F}_*\otimes C_*(\Omega M)^{\otimes k-1}\to \mathcal{G}_{k-1+*}
    \end{align*}
    such that for all $N\geq 0$, it holds:
    \begin{align*}
        \eta_{N+1}&\circ d_{\mathcal{F}\otimes C_*(\Omega)^{\otimes N-1}}+(-1)^{N+1}d_\mathcal{G}\circ \eta_{N+1}\\
        &=(-1)^{N+1}\eta_N(\mu_\mathcal{F}\otimes \id ^{\otimes N-1})-\mu_\mathcal{G}(\eta_N\otimes \id)+\sum_{r=1}^{N-1}(-1)^{N+1+r}\eta_N(\id^{\otimes r}\otimes \mu\otimes \id ^{\otimes {N-1-r}}).
    \end{align*}
\end{Def}

We denote $\boldsymbol{\widetilde{m}}$ for the degree $-1$ map on $T(C_*(\Omega M))\otimes \Z\Crit(f)\to T(C_*(\Omega M))\otimes \Z\Crit(f)$ defined
on $C_*(\Omega M)^{\otimes k}\otimes \Z\Crit(f)$ by $\id^{\otimes k}\otimes \boldsymbol{m}$.

\begin{Prop}[{\cite[Proposition 5.1]{riegel2024chain}}]\label{prop:Ainfty morphism DG}
    Let $\boldsymbol{\eta}\colon C_*(F_1)\to C_*(F_2)$ be a morphism of $\mathcal{A}_\infty$-modules over $C_*(\Omega M)$. The map $\widetilde{\eta}\colon C_*(M,C_*(F_1))\to C_*(M,C_*(F_2))$ defined by
    \begin{align*}
        \sum_{k\geq 0}(\eta_{k+1}\otimes \id)\circ (\id \otimes \boldsymbol{\widetilde{m}}^k)
    \end{align*}
    is a chain map.
\end{Prop}

\begin{Rem}
    The sum in the definition of $\widetilde{\eta}$ is in fact a finite sum because $\boldsymbol{\widetilde{m}}^{n+1}=0$.
\end{Rem}

The goal of this article to apply Proposition \ref{prop:Ainfty morphism DG} to a morphism of $\mathcal{A}_\infty$-modules $\boldsymbol{\nu}\colon C_*(\Omega M,x_0)\to C_{1-n+*}((\Omega M,x_0)^2)$ such that $\nu_1$ is the Goresky-Hingston coproduct on the based loop space.

\subsection{The Goresky-Hingston Coproduct}\label{subsec:GH coprod}
In this subsection, we recall the definition of the Goresky-Hingston coproduct following \cite{hingston2017product}.\\

We fix $\varepsilon>0$ such that $2\varepsilon$ is smaller than the injectivity radius of $M$. We denote the following tubular neighbourhoods
\begin{align*}
    U_\Delta&:=\{(x,y)\in M\times M\mid d_M(x,y)<\varepsilon\},&V_\Delta&:=\{(x,y)\in M\times M\mid d_M(x,y)<\tfrac{\varepsilon}{2}\}.   
\end{align*}
A \textit{Thom class of the diagonal} $\tau_\Delta\in C^n(U_\Delta,V_\Delta^c)$ is a cochain such that
\begin{align*}
    \iota^*\tau_\Delta\cap [M\times M]=\Delta_*[M]\in H_*(M\times M),
\end{align*}
where $\iota$ denotes inclusion $U_\Delta\hookrightarrow M\times M$.

We denote the map 
\begin{align*}
    e_\Lambda \colon I\times \Lambda  &\to M\times M,\\
    (t,(\gamma\colon [0,a]\to M))&\mapsto (\gamma(0), \gamma(ta)).
\end{align*}
We define the \textit{free figure 8 space} as
\begin{align*}
    \mathcal{F}\Lambda:=\mathcal{F}\Lambda M:=e_\Lambda^{-1}(\Delta(M))=\{(t,\gamma)\in I\times \Lambda\mid \gamma(0)=\gamma(ta)\}
\end{align*}
and define the map
\begin{align*}
    R\colon e_\Lambda^{-1}(U_\Delta)&\to \mathcal{F}\Lambda
\end{align*}
sending $(t,\gamma)$ with $(\gamma(0),\gamma(ta))=(x,y)\in U_\Delta$ to $(t',\gamma|_{[0,ta]}*\overline{\lambda_{x,y}}*{\lambda_{x,y}}*\gamma|_{[ta,a]})$, where $\lambda_{x,y}$ is the unit speed, geodesic path from $x$ to $y$ and $t'$ is such that $e_\Lambda$ evaluates $(t',\gamma|_{[0,ta]}*\overline{\lambda_{x,y}}*{\lambda_{x,y}}*\gamma|_{[ta,a]})$ at the tip of the two geodesics paths. This map is a homotopy inverse of the inclusion $\mathcal{F}\Lambda\hookrightarrow e_\Lambda^{-1}(U_\Delta)$.

Finally on the figure 8 space, we can define a cut map
\begin{align*}
    \mathrm{c}\colon \mathcal{F}\Lambda&\to \Lambda\times \Lambda,\\
    (t,\gamma)&\mapsto (\gamma|_{[0,ta]},\gamma|_{[ta,a]}).
\end{align*}

\begin{Def}
    A {chain-level description of the} $\nu_{GH}$ \textit{Goresky-Hingston} coproduct is given by the composition
    \begin{equation}\label{eq:def of GH coprod}
        \begin{split}
            C_*(\Lambda,M)&\overset{I\times }{\to} C_{1+*}(I\times \Lambda,\partial I\times \Lambda\cup I\times M)\\
            &\to C_{1+*}(I\times \Lambda,\partial I\times \Lambda \cup I\times M\cup e^{-1}_\Lambda(V_\Delta^c))\\
            &\simeq C_{1+*}(e^{-1}_\Lambda(U_\Delta),\partial I\times \Lambda \cup I\times M\cup e^{-1}_\Lambda(V_\Delta^c))\\
            &\overset{e_\Lambda^*\tau_\Delta\cap }{\to} C_{1-n+*}(e^{-1}_\Lambda(U_\Delta),\partial I\times \Lambda\cup I\times M)\\
            &\overset{R_*}{\to} C_{1-n+*}(\mathcal{F}\Lambda,\partial I\times \Lambda\cup I\times \Lambda)\\
            &\overset{\mathrm{c}_*}{\to}C_{1-n+*}((\Lambda,M)^2).
        \end{split}
    \end{equation}
\end{Def}
\begin{Rem}
    This definition agrees with the definition of the Goresky-Hingston coproduct in \cite[Definition 1.4]{hingston2017product} up to the sign $(-1)^k$ on $C_k(\Lambda,M)$. We choose these sign conventions because we have the convention that maps act from the left. So in particular we multiply by the interval $I$ from the left while Hingston and Wahl multiply by $\times I$ in their definition. Our convention leads to slightly neater signs (see \cite[Section 6.3]{clivio2025graded}).
    
    Moreover, we define it as a map with target $C_{1-n+*}((\Lambda,M)^2)$ rather than $C_{1-n+*}(\Lambda,M)^{\otimes 2}$. These two definitions are equivalent after post-composing our definition with a diagonal approximation. For us, it is however more convenient to consider the above map as the coproduct.
\end{Rem}

\subsubsection{Based Loop Space}

The Goresky-Hingston coproduct restricts to a coproduct on $C_*(\Omega,x_0)$. We describe it here explicitly. We denote the following open sets
\begin{align*}
    U_{x_0}&:=\{x\in M\mid d_M(x,x_0)<\varepsilon\},&V_{x_0}&:=\{x\in M\mid d_M(x,x_0)<\tfrac{\varepsilon}{2}\}.   
\end{align*}
We denote the inclusion
\begin{align*}
    \iota_{x_0}\colon (U_{x_0},V_{x_0})&\to (U_\Delta,V_\Delta),\\
    x&\mapsto (x_0,x).
\end{align*}
A \textit{Thom class of the basepoint} or simply a local orientation is given by $\tau_{\Omega}:=\iota_{x_0}^*\tau_\Delta\in C^n(U_{x_0},V_{x_0}^c)$.
We denote the map 
\begin{align*}
    e_\Omega \colon I\times \Omega  &\to M,\\
    (t,(\gamma\colon [0,a]\to M))&\mapsto  \gamma(ta)
\end{align*}
and define the \textit{based figure 8 space} as
\begin{align*}
    \mathcal{F}\Omega:=\mathcal{F}\Omega M:=e_\Omega^{-1}(\{x_0\})=\{(t,\gamma)\in I\times \Omega\mid \gamma(t)=x_0\}.
\end{align*}
A chain model of the Goresky-Hingston coproduct as in \eqref{eq:def of GH coprod} can be computed on $C_*(\Omega,x_0)$ by
\begin{align*}
    C_*(\Omega,x_0)&\overset{I\times }{\to} C_{1+*}(I\times \Omega,\partial I\times \Omega\cup I\times \{x_0\})\\
    &\to C_{1+*}(I\times \Omega,\partial I\times \Omega \cup I\times \{x_0\}\cup e^{-1}_\Omega(V_{x_0}^c))\\
    &\simeq C_{1+*}(e^{-1}_\Omega(U_{x_0}),\partial I\times \Omega \cup I\times \{x_0\}\cup e^{-1}_\Omega(V_{x_0}^c))\\
    &\overset{e_\Omega^*\tau_{x_0}\cap }{\to} C_{1-n+*}(e^{-1}_\Omega(U_{x_0}),\partial I\times \Omega\cup I\times \{x_0\})\\
    &\overset{R_*}{\to} C_{1-n+*}(\mathcal{F}\Omega,\partial I\times \Omega\cup I\times \Omega)\\
    &\overset{\mathrm{c}_*}{\to}C_{1-n+*}((\Omega,x_0)^2).
\end{align*}

\subsection{Sign Conventions}

In this article, we heavily use maps of chain complexes that do not preserve degree. For such maps, the signs are delicate. In this subsection, we specify what convention we use.\\

Let $C_*$ denote a chain complex with differential $d$. Then $C_{k+*}:=\Sigma^d C_*$ denotes the chain complex $C_*$ with the differential $\Sigma^kd:=(-1)^kd$.

Thus a degree-shifting map $f\colon C_*\to D_{|f|+*}$ is a chain map if on $\alpha\in C_*$, we have
\begin{align*}
    f(d \alpha)=(-1)^{|f|}d f(\alpha).
\end{align*}

\begin{Exmp}
    The chain-level Goresky-Hingston coproduct as in \eqref{eq:def of GH coprod} is a chain map with this sign convention. Indeed, in the composition \eqref{eq:def of GH coprod} the maps 
    \begin{align*}
        C_{1+*}(I\times \Lambda,\partial I\times \Lambda\cup I\times M)&\to C_{1+*}(I\times \Lambda,\partial I\times \Lambda \cup I\times M\cup e^{-1}_\Lambda(V_\Delta^c))\\
        &\simeq C_{1+*}(e^{-1}_\Lambda(U_\Delta),\partial I\times \Lambda \cup I\times M\cup e^{-1}_\Lambda(V_\Delta^c))
    \end{align*}
    and 
    \begin{align*}
         C_{1-n+*}(e^{-1}_\Lambda(U_\Delta),\partial I\times \Lambda\cup I\times M)&\overset{R_*}{\to} C_{1-n+*}(\mathcal{F}\Lambda,\partial I\times \Lambda\cup I\times \Lambda)\\
        &\overset{\mathrm{c}_*}{\to}C_{1-n+*}((\Lambda,M)^2)
    \end{align*}
    are degree-preserving chain maps.
    
    We compute for $\alpha\in C_*(\Lambda,M)$:
    \begin{align}\label{eq:differential sign Ix}
        d  (I\times\alpha)&= (d I)\times \alpha- I\times (d \alpha)=0-I\times (d \alpha)
    \end{align}
    For $\alpha\in C_{1+*}(e^{-1}_\Lambda(U_\Delta),\partial I \times \Lambda\cup I\times M\cup e^{-1}_\Lambda(V_\Lambda^c))$, we compute:\footnote{The signs in \cite[Proposition VI.5.1. (v)]{bredon2013topology} are computed as $d(f\cap c)=\delta (f)\cap c+(-1)^{|f|}f\cap dc$. This is a sign mistake in \cite[Page 335, Line -2]{bredon2013topology}. The signs that we use are $d(f\cap c)=(-1)^{|f|}(f\cap dc-\delta (f)\cap c)$.\\    
    The signs in \cite[Section 3.3]{hatcher2002algebraic} are correct. However, we have the opposite convention that in the cap product the cochain acts from the left on the chain.}
    \begin{align}\label{eq:differential sign tau cap}
        d (e^*_\Lambda\tau_\Delta \cap \alpha)=(-1)^{n}(e^*_\Lambda\tau_\Delta \cap d \alpha- \overset{=0}{\overbrace{\delta (e^*_\Lambda\tau_\Delta)}} \cap \alpha)=(-1)^ne^*_\Lambda\tau_\Delta \cap d \alpha.
    \end{align}
    Combining \eqref{eq:differential sign Ix} and \eqref{eq:differential sign tau cap}, we find 
    \begin{align*}
        \nu_\mathrm{GH}(d \alpha)=(-1)^{1-n}d \nu_{\mathrm{GH}}(\alpha)
    \end{align*}
    and thus $\nu_\mathrm{GH}$ is a chain map with our sign conventions.
\end{Exmp}

This sign convention has moreover the following consequence for DG coefficients: let $F,G$ be fibres of fibrations over $M$. We write the differential on $C_*(M,C_*(F))$ and $C_*(M,C_*(G))$ as 
\begin{align*}
    d_F&=d\otimes \id +(\Phi_F\otimes \id )\circ (\id \otimes \boldsymbol{m}),\\
    d_G&=d\otimes \id +(\Phi_G\otimes \id )\circ (\id \otimes \boldsymbol{m}).
\end{align*}
Let $f\colon C_*(F)\to C_{|f|+*}(G)$ be a chain map which is equivariant under the $\Omega$-action. Then $f\otimes \id$ is a chain map with respect to $d_F$ and $\Sigma^{|f|}d_G$
\begin{align*}
    (f\otimes \id)\circ d_F&=(f\otimes \id)\circ (d\otimes \id)+(f\otimes \id)\circ (\Phi_F\otimes \id )\circ (\id \otimes \boldsymbol{m})\\
    &=(fd\otimes \id)+(\Phi_G\otimes \id)\circ(f \otimes\boldsymbol{m})\\
    &= (-1)^{|f|}d f\otimes \id+(-1)^{|f|}(\Phi_G\otimes \id)\circ(\id\otimes \boldsymbol{m})\circ (f\otimes \id)\\
    &=\Sigma^{|f|}d_G\circ (f\otimes \id).
\end{align*}

\section{Transitive Diffeomorphism Lifts}\label{sec:trans diff lifts}
In this section, we turn our attention to constructing particularly nice transitive lifting functions on the fibration $e_0\colon \Lambda M \to M$. These transitive lifting functions are constructed such that they send the spaces $e_\Lambda^{-1}(U_\Delta), e_\Lambda^{-1}(V_\Delta^c), \mathcal{F}\Lambda M$ to themselves. These are the spaces that appear in a chain-level description of the Goresky-Hingston coproduct as in \eqref{eq:def of GH coprod}. We construct these better-behaved transitive lifting functions by defining and constructing transitive diffeomorphism lifts.

\subsection{Definition}

\begin{Def}
    A \textit{transitive diffeomorphism lift} is a map
    \begin{align*}
        \Theta\colon \mathcal{P}M\to \Diff(M)
    \end{align*}
    such that
    \begin{enumerate}[(i)]
        \item\label{Theta property(i)} for $c_x$ the constant path at $x\in M$, $\Theta(c_x)=\id$;
        \item\label{Theta property(ii)} for composable $\lambda,\delta\in \mathcal{P}M$, it holds
        \begin{align*}
            \Theta(\lambda\cdot \delta)=\Theta(\delta)\circ \Theta(\lambda);
        \end{align*}
        \item\label{Theta property(iii)} for all $x,y\in M$ and $\lambda\in \mathcal{P}_{x\to y}M$, it holds
        \begin{align*}
            \Theta(\lambda)(x)=y;
        \end{align*}
        \item\label{Theta property(iv)} for all $x,y,z\in M$ and $\lambda\in \mathcal{P}_{x\to y}M$ with $d_M(x,z)<\varepsilon$, it holds
        \begin{align*}
            d_M(\Theta(\lambda)(z),y)<\varepsilon;
        \end{align*}
        \item\label{Theta property(v)} for all $x,y,z\in M$ and $\lambda\in \mathcal{P}_{x\to y}M$ with $d_M(x,z)>\tfrac{\varepsilon}{2}$, it holds 
        \begin{align*}
            d_M(\Theta(\lambda)(z),y)>\tfrac{\varepsilon}{2}.
        \end{align*}
    \end{enumerate}
\end{Def}

Before we prove that such a transitive diffeomorphism lift always exists in Proposition \ref{prop:constr of tdl}, we prove that, given such a function, we can construct a particularly nice transitive lifting function.

\begin{Def}
    Let $\Theta\colon \mathcal{P}M\to \Diff(M)$ be a transitive diffeomorphism lift. The \textit{associated transitive lifting function}
    is defined by
    \begin{align*}
        \Phi\colon \Lambda M\times\mathcal{P}M&\to  \Lambda M,\\
        (\gamma,\lambda)&\mapsto \Theta(\lambda)_*(\gamma).
    \end{align*}
\end{Def}

\begin{Not}
    We denote $e_t\colon I\times \Lambda M\to M$ for the map that sends $(t,\gamma\colon[0,a]\to M)$ to $\gamma(ta)$.
\end{Not}

\begin{Lem}\label{lem:ass trans lift fun}
    Let $\Theta\colon \mathcal{P}M\to \Diff(M)$ be a transitive diffeomorphism lift and $\Phi \colon \Lambda M\times\mathcal{P}M\to  \Lambda M$ its associated transitive lifting function. Then the following holds:
    \begin{enumerate}[(a)]
        \item\label{item:ass trans lift fun (a)} the restriction of $\Phi$ to $\Lambda M\,{}_{e_0}\!\times_{e_0}\mathcal{P}M$ and the map
        \begin{align}\label{eq:trans lift on IxLambda}
            \begin{split}
                \Phi\colon (I\times \Lambda M)\, {}_{e_0}\!\times_{e_0}\mathcal{P}M&\to I\times \Lambda M,\\
                ((t,\gamma),\lambda)&\mapsto (t,\Phi(\gamma,\lambda))
            \end{split}
        \end{align}
        are transitive lifting functions;
        \item\label{item:ass trans lift fun (b)} as maps $I\times \Lambda M\times \mathcal{P}M\to M$
        \begin{align*}
            e_t\circ \Phi=\Theta(\lambda);
        \end{align*}
        \item\label{item:ass trans lift fun (c)} the map \eqref{eq:trans lift on IxLambda} restricts to a transitive lifting function on $e_\Lambda^{-1}(U_\Lambda), e_\Lambda^{-1}(V_\Lambda^c),\mathcal{F}\Lambda\subseteq I\times \Lambda$.
    \end{enumerate}
\end{Lem}
\begin{proof}
    We first prove \eqref{item:ass trans lift fun (a)} for the fibration $e_0\colon \Lambda M\to M$. Property \eqref{Theta property(i)} of $\Theta$ ensures that $\Phi(\gamma,c_{\gamma(0)})=\gamma$. Property \eqref{Theta property(ii)} gives that for composable paths $\lambda,\delta\in \mathcal{P}M$, we have 
    \begin{align*}
        \Phi(\gamma,\lambda\cdot \delta)=\Phi(\Phi(\gamma,\lambda),\delta).
    \end{align*}
    Finally, property \eqref{Theta property(iii)} ensures that $e_0\circ \Phi=e_1\circ \mathrm{pr}_2$.
    
    It is immediate that for any transitive lifting function $\Phi\colon \Lambda M\,{}_{e_0}\!\times_{e_0}\mathcal{P}M\to \Lambda$ the map $\id\times \Phi\colon (I\times \Lambda M)\,{}_{e_0}\!\times_{e_0}\mathcal{P}M\to I\times \Lambda M$ is a transitive lifting function. This shows \eqref{item:ass trans lift fun (a)}.
    
    We compute for $(t,\gamma)\in I\times \Lambda$ and $\lambda\in \mathcal{P}$
    \begin{align*}
        e_t\circ \Phi((t,\gamma),\lambda)=e_t(t,\Phi(\gamma,\lambda))=e_t(t,\Theta(\lambda)_*(\gamma))=\Theta(\lambda)(e_t(t,\gamma)).
    \end{align*}
    
    For \eqref{item:ass trans lift fun (c)}, we need to check that, for $((t,\gamma),\lambda)\in (I\times \Lambda M)\,{}_{e_0}\!\times_{e_0}\mathcal{P}M$, it holds
    \begin{align*}
        \Phi((t,\gamma),\lambda)\in \begin{cases}
            e_\Lambda^{-1}(U_\Lambda) &\text{if } (t,\gamma) \in e^{-1}_\Lambda(U_\Lambda),\\
            e_\Lambda^{-1}(V_\Lambda^c) &\text{if } (t,\gamma) \in e^{-1}_\Lambda(V_\Lambda^c),\\
            \mathcal{F}\Lambda &\text{if } (t,\gamma) \in \mathcal{F}\Lambda;
        \end{cases}
    \end{align*} 
    Let $\lambda\in \mathcal{P}_{x\to y} M$ be a path and $(t,\gamma)\in I\times \Lambda$ with $\gamma(0)=x$. If $(t,\gamma)\in e^{-1}_\Lambda(U_\Lambda)$, then 
    \begin{align*}
        d_M(x,e_t(t,\gamma))<\varepsilon.
    \end{align*}
    We then compute using \eqref{item:ass trans lift fun (b)} and Property \eqref{Theta property(iv)} of $\Theta$
    \begin{align*}
        d_M(y,e_t\circ \Phi((t,\gamma),\lambda))\overset{\eqref{item:ass trans lift fun (b)}}{=}d_M(y,\Theta(\lambda)(e_t(t,\gamma)))\overset{\eqref{Theta property(iv)}}{<}\varepsilon.
    \end{align*}
    Therefore we have $\Phi(\lambda)(t,\gamma)\in e_\Lambda^{-1}(U_\Delta)$. A similar computation shows that if $(t,\gamma)\in e_\Lambda^{-1}(V_\Delta^c)$, then $\Phi(\lambda)(t,\gamma)\in e_\Lambda^{-1}(V^c_\Delta)$ by Property \eqref{Theta property(v)} of $\Theta$.
    
    If $(t,\gamma)\in \mathcal{F}\Lambda$ holds, then $e_\Lambda(t,\gamma)=x$ and thus by \eqref{item:ass trans lift fun (b)} and Property \eqref{Theta property(iii)}, we compute
    \begin{align*}
        e_t\circ \Phi((t,\gamma),\lambda)\overset{\eqref{item:ass trans lift fun (b)}}{=}\Theta(\lambda)(e_t(t,\gamma))=\Theta(\lambda)(x)\overset{\eqref{Theta property(iii)}}{=}y.
    \end{align*}
    Therefore, we have $\Phi(\lambda)(t,\gamma)\in\mathcal{F}\Lambda$.
\end{proof}

\subsection{Construction}
In this subsection, we construct a transitive diffeomorphism lift for $M$. We first construct for each path $\lambda\in \mathcal{P}M$ a time-dependent vector field on $M$. Then the flow along this vector field defines the diffeomorphism $\Theta(\lambda)$.\\

We first construct $\Theta(\lambda)$ for $\lambda\colon [0,a]\to M$ smooth. We define $\Theta(\lambda)$ as the flow along a time dependent vector field:
\begin{align*}
    V_\lambda\colon [0,a]\times M\to TM.
\end{align*}

To ensure Property \eqref{Theta property(iii)} that $\Theta(\lambda)$ sends its start point to its end point, we require $V_\lambda(t,\lambda(t)):=\dot{\lambda}(t)$.

We recall that $\varepsilon>0$ was chosen such $2\varepsilon$ is smaller than the injectivity radius of $M$. The vector field $V_\lambda$ is such that for all $t\in [0,a]$ the vector field is zero outside of $B_{2 \varepsilon}(\lambda(t))$.

\begin{figure}[H]
    \centering
    \begin{minipage}{0.45\textwidth}
        \centering
        \captionsetup{width=0.85\textwidth}
        \includegraphics[width=\textwidth]{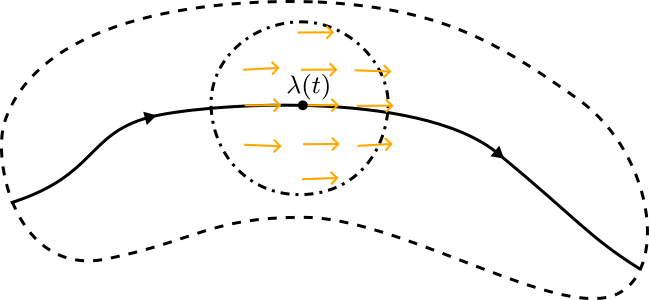} 
        \caption{The vector field $v_{\lambda,t}$ is defined on the ball $B_{2\varepsilon}(\lambda(t))$ and points in the direction of $\dot{\lambda}(t)$.}\label{fig:v lambda}
    \end{minipage}\hfill
    \begin{minipage}{0.45\textwidth}
        \centering
        \captionsetup{width=0.85\textwidth}
        \includegraphics[width=\textwidth]{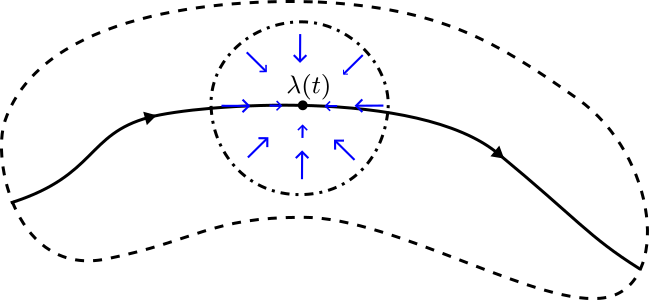} 
        \caption{The vector field $w_{\lambda,t}$ is defined on the ball $B_{2\varepsilon}(\lambda(t))$ and points inwards.}\label{fig:w lambda}
    \end{minipage}
\end{figure}

It thus remains to explain how $V_\lambda(t,-)$ is defined on $B_{2\varepsilon}(\lambda(t))$. Independently of $\lambda$, we fix a smooth functions $\chi_1\colon \R^{\geq 0}\to \R$ such that
\begin{itemize}
    \item $\chi_1|_{[0,\varepsilon]}\equiv 1$;
    \item $\chi_1|_{[2\varepsilon,\infty)}\equiv 0$;
\end{itemize}
For $x\in B_{2\varepsilon}(\gamma(t))$, let $\sigma_x\colon I\to B_{2\varepsilon}(\lambda(t))$ be the geodesic from $\gamma(t)$ to $x$ and let $P_{x}\colon T_{\lambda(t)} M\to T_{x}M$ the parallel transport along it. We set 
\begin{align*}
    v_{\lambda,t}(x):=P_{x}(\dot{\lambda}(t))\in T_xM.
\end{align*}
See Figure \ref{fig:v lambda}. We have in particular $v_{\lambda,t}(\lambda(t))=\dot{\lambda}(t)$.

We now set 
\begin{align*}
    V_{\lambda,1}\colon [0,a]\times M&\to TM,\\
    (t,x)&\mapsto \begin{cases}
        \chi_1(d_M(x,\lambda(t)))\cdot v_{\lambda,t}(x) &\text{if } x\in B_{2\varepsilon}(\lambda(t));\\
        0 &\text{else.}
    \end{cases}
\end{align*}
Taking the flow of this vector field guarantees the identity property as in \eqref{Theta property(i)}, the composability as in \eqref{Theta property(ii)} and that it sends the start point to the end point as in \eqref{Theta property(iii)}. (We argue this in further detail in the proof of Proposition \ref{prop:constr of tdl}). However to achieve Properties \eqref{Theta property(iv)} and \eqref{Theta property(v)}, we need a correction term.

\begin{figure}[H]
    \centering
    \includegraphics[width=0.65\linewidth]{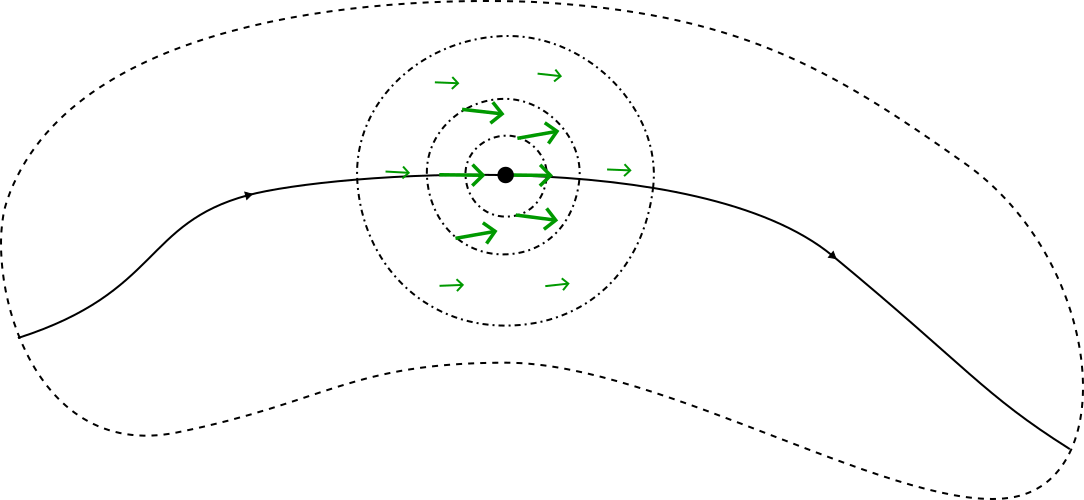}
    \caption{The vector field $V_\lambda$ is a linear combination of $v_{\lambda,t}(x)$ and $w_{\lambda,t}(x)$. In particular, it is zero for $x$ outside $B_{2\varepsilon}(\lambda(t))$ and $v_{\lambda,t}(x)$ on $B_{\frac{\varepsilon}{4}}(\gamma(t))$. The $w_{\lambda,t}(x)$-term is such that $V_\lambda$ moves the points with $d_M(x,\lambda(t))=\varepsilon$ closer to $\lambda(t)$ and the points with $d_M(x,\lambda(t))=\frac{\varepsilon}{2}$ further away from $\lambda(t)$.}
    \label{fig:V_lambda}
\end{figure}

Our strategy is to define a second vector field $V_{\lambda,2}\colon[0,a]\times M\to TM$ such that $V_{\lambda}:=V_{\lambda,1}+V_{\lambda,2}$ \enquote{points inward} on the sphere around $\lambda(t)$ with radius $\varepsilon$ and \enquote{points outward} on the sphere with radius $\frac{\varepsilon}{2}$ (see Figure \ref{fig:V_lambda}). To make precise what \enquote{points inward/outward} means in our context, we recall the following formula for the derivative of the Riemannian distanc:
\begin{Lem}\label{lem:derivative of distance}
    Let $x, y\in M$ with $d_M(x,y)<2\varepsilon$ and let $v_x\in T_xM$ and $v_y\in T_yM$ be tangent vectors. It holds that
    \begin{align*}
        D_{v_x\oplus v_y}(d_M(x,y)^2)=-2\langle v_x,\exp_x^{-1}(y)\rangle_x-2\langle v_y,\exp_y^{-1}(x)\rangle_y
    \end{align*}
    where $\langle-,-\rangle_x$ and $\langle-,-\rangle_y$ denote the Riemannian metric at $x$ and $y$, respectively.
\end{Lem}
\begin{proof}
    We note that $d_M(x,y)<2\varepsilon$ guarantees that their distance is smaller than the injectivity radius of $M$ and therefore $\exp_x^{-1}(y)$ and $\exp_y^{-1}(x)$ are well-defined.
    
    We have
    \begin{align*}
        D_{v_x\oplus v_y}(d_M(x,y)^2)&=D_{v_x}(d_M(x,y)^2)+D_{v_y}(d_M(x,y)^2).
    \end{align*}
    To compute $D_{v_x}(d_M(x,y)^2)$, we fix a path $p\colon (-a,a)\to M$ such that $p(0)=x$ and $\dot{p}(0)=v_x$. We define the variation $\Gamma(s,t)=\gamma_s(t)$ where $\gamma_s\colon[0,1]\to M$ is the constant-speed minimizing path with $\gamma_s(0)=y$ and $\gamma_s(1)=p(s)$. We consider the energy functional $E(\gamma_s)=\frac{1}{2}\int \langle \dot\gamma_s(t), \dot\gamma_s(t)\rangle dt$. We thus have
    \begin{align*}
        D_{v_x}(d_M(x,y)^2)=\partial_s|_{s=0} d_M(y,p(s))^2=2\partial_s|_{s=0}E(\gamma_s).
    \end{align*}
    The first variation formula (e.g.\ see \cite[Theorem 3.31]{gallot1990riemannian}) then gives
    \begin{align*}
        \partial_s|_{s=0}E(\gamma_s)&=\langle \partial_s|_{s=0}\gamma_s(1), \dot{\gamma}_0(1)\rangle-\langle\partial_s|_{s=0} \gamma_s(0),\dot{\gamma}(0)\rangle-\int_0^1\langle \partial_s|_{s=0} \gamma_s(t),\partial_t\dot{\gamma}(t)\rangle dt\\
        &=\langle v_x, \dot{\gamma}_0(1)\rangle-\langle0,\dot{\gamma}(0)\rangle-\int_0^1\langle \partial_s|_{s=0} \gamma_s(t),0\rangle dt=\langle v_x, \dot{\gamma}_0(1)\rangle
    \end{align*}
    We note that $\dot{\gamma}_0(1)=-\exp_x(y)$ because $\gamma_0(t)=\exp_x(\exp^{-1}_x (y)(1-t))$.
    
    We deduce that $D_{v_x}(d_M(x,y)^2)=-2\langle v_x,\exp_x^{-1}(y)\rangle_x$ and, by symmetry, $D_{v_y}(d_M(x,y)^2)=-2\langle v_y,\exp_y^{-1}(x)\rangle_y$. This concludes the proof.
\end{proof}

For $z\in B_{2\varepsilon}(\lambda(t))$, we denote by $w_{\lambda,t}(z):=\exp_z^{-1}(\lambda(t))\in T_zM$. This is well-defined because $d_M(z,\lambda(t))$ is smaller than the injectivity radius. This is a vector field on $B_{2\varepsilon}(\lambda(t))$ pointing towards $\lambda(t)$ (see Figure \ref{fig:w lambda}). We also define $u_{\lambda,t}(z):=\exp_{\lambda(t)}^{-1}(z)$.

Independently of $\lambda$, we fix some $\delta>0$ and a smooth function $\chi_2\colon\R^2\times \R^{\geq 0}\to \R$ such that 
\begin{itemize}
    \item $\chi_2(x,y,d)= 0$ for $d<\tfrac{\varepsilon}{4}$ and $d>2\varepsilon$;
    \item $\chi_2(x,y,\varepsilon)=x$;
    \item $\chi_2(x,y,\tfrac{\varepsilon}{2})=y$.
\end{itemize}
Let $\delta_1^\lambda(t)\in \R$ be the infimum on all $\delta_1$ such that, for all $z\in M$ with $d_M(z,\lambda(t))=\varepsilon$, we have
\begin{align}\label{eq:def of delta1}
    \left\langle v_{\lambda,t}(z)+\delta_1w_{\lambda,t}(z),w_{\lambda,t}(z)\right\rangle_z +\left\langle\dot{\lambda}(t),u_{\lambda,t}(z)\right\rangle_{\lambda(t)}\geq \delta.
\end{align}
Similarly, let $\delta_2^\lambda(t)$ be the supremum on all $\delta_2$ such that, for all $z\in M$ with $d_M(z,\lambda(t))=\tfrac{\varepsilon}{2}$, we have
\begin{align*}
    \left\langle v_{\lambda,t}(z)+\delta_2w_{\lambda,t}(z),w_{\lambda,t}(z)\right\rangle_z +\left\langle\dot{\lambda}(t),u_{\lambda,t}(z)\right\rangle_{\lambda(t)}\leq-\delta.
\end{align*}
Both $\delta_1(t)$ and $\delta_2(t)$ are piecewise smooth in $t$.

We then define 
\begin{align*}
    V_{\lambda,2}\colon [0,a]\times M &\to TM,\\
    (t,z)&\mapsto \begin{cases}
        \chi_2(\delta_1^\lambda(t),\delta_2^\lambda(t),d_M(z,\lambda(t)))\cdot  w_{\lambda,t}(z) & \text{if }\\
        z\in B_{2\varepsilon}(\lambda(t)),
        0& \text{else.}
    \end{cases}
\end{align*}
By our assumption on $\chi_2$, this is a vector field that satisfies $V_{\lambda,2}(t,\lambda(t))=0$, $V_{\lambda,2}(t,z)= \delta_1^\lambda(t)w_{\lambda,t}(z)$ if $d_M(z,\lambda(t))=\varepsilon$ and $V_{\lambda,2}(t,z))=\delta_2^\lambda(t)w_{\lambda,t}(z)$ if $d_M(z,\lambda(t))=\frac{\varepsilon}{2}$.

We define $V_\lambda:=V_{\lambda,1}+V_{\lambda,2}$. This vector field is smooth for a fixed $t$ and piecewise smooth in $t$ because both $V_{\lambda,1}$ and $V_{\lambda,2}$ are. Therefore, the flow
\begin{align*}
    \theta^\lambda_t\colon M\to M
\end{align*}
is a well-defined diffeomorphism for all $t\in [0,a]$. We define $\Theta(\lambda):=\theta^\lambda_a$.

For a general piecewise smooth $(\lambda\colon [0,a]\to M)\in \mathcal{P}M$, we fix $0=a_0<\dots <a_k=a$ for minimal $k\geq 0$ such that $\lambda|_{[a_i,a_{i+1}]}$ is smooth. We then define
\begin{align}\label{eq:Theta def compo}
    \Theta(\lambda):=\Theta(\lambda|_{[a_{k-1},a_k]})\circ \dots \circ \Theta(\lambda|_{[a_0,a_1]}).
\end{align}

\begin{Prop}\label{prop:constr of tdl}
    The map $\Theta\colon \mathcal{P}M\to \Diff(M)$ as constructed above is a transitive diffeomorphism lift.
\end{Prop}
\begin{proof}
    The first condition \eqref{Theta property(i)} of a transitive diffeomorphism lift states that $\Theta(c_x)$ is the identity for $c_x$ a constant path. This is true because in this case, we flow for time $0$ and thus $\Theta(c_x)=\theta^\lambda_0=\id$.
    
    The second condition \eqref{Theta property(ii)} states that $\Theta$ commutes with composition. That means that for $\lambda_1$ and $\lambda_2$ composable paths, we have $\Theta(\lambda_1\cdot \lambda_2)=\Theta(\lambda_2)\circ\Theta(\lambda_1)$. If the composition is not smooth, then this property holds directly from the definition as in \eqref{eq:Theta def compo}. If the composition is smooth, this follows from the composition property of the flow and the fact that $V_\lambda(t,-)$ only depends on $\lambda(t)$ and $\dot{\lambda}(t)$.
    
    The third condition \eqref{Theta property(iii)} states that if $\lambda$ is a path from $x$ to $y$, then $\Theta(\lambda)(x)=y$. By condition \eqref{Theta property(ii)}, we can assume that $\lambda$ is smooth. We have
    \begin{align*}
        V_\lambda(t,\lambda(t))
        &= V_{\lambda,1}(t,\lambda(t))+V_{\lambda,2}(t,\lambda(t))=\dot{\lambda}(t)+0.
    \end{align*}
    Therefore $\lambda(t)=\theta^\lambda_t(x)$ and in particular $\Theta(\lambda)(x)=\theta^\lambda_a(x)=\lambda(a)=y$.
    
    The fourth condition \eqref{Theta property(iv)} states that for $x,y,z\in M$ with $\lambda$ a path from $x$ to $y$ and $d_M(x,z)<\varepsilon$, we have $d_M(y,\Theta(\lambda)(z))<\varepsilon$. Again by condition \eqref{Theta property(ii)}, we can assume that $\lambda$ is smooth. By continuity it is enough to show that for $z$ with $d_M(z,x)=\varepsilon$ that we have 
    \begin{align*}
        \left.\partial_t\right|_{t=0} \left[d_M(\theta^\lambda_{t}(z),\lambda(t))^2\right]<0.
    \end{align*}
    We compute
    \begin{align*}
        \left.\partial_t\right|_{t=0} \left[d_M(\theta^\lambda_{t}(z),\lambda(t))^2\right]&=\left.D_{\partial_t \theta^\lambda_{t}(z)\oplus \partial_t \lambda(t)}\left[d_M(\theta^\lambda_{t}(z),\lambda(t)\right]\right|_{t=0}=D_{V_{\lambda}(0,z)\oplus \dot{\lambda}(t)}\left[d_M(z,\lambda(t))^2\right].
    \end{align*}
    By Lemma \ref{lem:derivative of distance}, we find
    \begin{align*}
        D_{V_{\lambda}(0,z)\oplus \dot{\lambda}(t)}\left[d_M(z,x)^2\right]&=-2\langle V_{\lambda}(0,z),\exp_z^{-1}(x)\rangle_{z}-2\langle \dot{\lambda}(0),\exp_{x}^{-1}(z)\rangle_{x}.
    \end{align*}
    In our notation, we have $\exp_{x}^{-1}(z)=u_{\lambda,0}(z)$ and $\exp_z^{-1}(x)=w_{\lambda,0}(z)$. Because $d_M(z,x)=\varepsilon$ and the definition of $\chi_1$ and $\chi_2$, we have $V_{\lambda}(0,z)=v_{\lambda,0}(z)+\delta_1^\lambda(0)w_{\lambda,0}(z)$. We therefore have
    \begin{align*}
        \left.\partial_t\right|_{t=0} \left[d_M(\theta^\lambda_{t}(z),\lambda(t))^2\right]&=-2\left\langle v_{\lambda,0}(z)+\delta_1^\lambda(0)w_{\lambda,0}(z),w_{\lambda,0}(z)\right\rangle_z -2\left\langle\dot{\lambda}(0),u_{\lambda,0}(z)\right\rangle_{x}
    \end{align*}
    which is smaller than $-2\delta<0$ by definition of $\delta_1^\lambda(0)$ as in \eqref{eq:def of delta1}. This shows condition \eqref{Theta property(iv)}.
    
    Condition \eqref{Theta property(v)} follows analogously from Lemma \ref{lem:derivative of distance} and the definition of $\delta_2(t)$.
    
    This shows that $\Theta$ as constructed above is a transitive diffeomorphism lift.
\end{proof}

\section{Goresky-Hingston Coproduct as $\mathcal{A}_\infty$-Morphism}\label{sec:ainfty struc}
In this section, we construct an $\mathcal{A}_\infty$-morphism $\boldsymbol{\nu}\colon  C_*(\Omega M,x_0)\to C_{1-n+*}((\Omega M,x_0)^2)$. The map $\nu_1$ will be the Goresky-Hingston coproduct on the based loop space. This then lets us define a map $\widetilde{\nu}\colon C_*(M,C_*(\Omega M,x_0))\to C_*(M ,C_{1-n+*}((\Omega M,x_0)^2)$ on Morse homology with DG coefficients.\\

We fix some transitive diffeomorphism lift $\Theta\colon \mathcal{P} M\to \Diff(M)$ and its associated transitive lifting functions $\Phi\colon \Lambda M \,{}_{e_0}\!\times_{e_0}  \mathcal{P}M\to \Lambda M$ and $\Phi\colon  (I\times\Lambda M) \,{}_{e_0}\!\times_{e_0}  \mathcal{P}M\to I\times \Lambda M$. By Lemma \ref{lem:ass trans lift fun}, this restricts to a right $\Omega M$-action on $\Omega M$, $I\times \Omega M$, $e_\Omega^{-1}(U_{x_0})$, $e_\Omega^{-1}(V_{x_0}^c)$, $\mathcal{F}\Omega$ and $\mathcal{H}:=\partial I\times \Omega\cup I\times \{x_0\}\subseteq I\times \Omega$. It also acts on $\Omega M\times \Omega M$ via the diagonal action.

The based loop space $\Omega M$ thus plays a double role: it is the monoid acting but also the space that is acted on. To minimize the resulting confusion, we use the convention that we denote a loop $\lambda$ or $\mu$ if it is acting and $\gamma$ if it is acted on.

In this section, we study the (failure of) commutativity of the following diagram
\begin{center}
    \begin{tikzcd}
        C_*(\Omega,x_0)  \otimes C_*(\Omega) \ar[d, "(I\times) \otimes \id"]\ar[r,"\Phi_*"] \ar[rd,phantom ,"\eqref{sq:I}"]& C_*(\Omega, x_0)\ar[d, "I\times "]\\
        C_{1+*}(I\times \Omega, \mathcal{H})\otimes C_*(\Omega) \ar[d] \ar[r, "\Phi_*"] \ar[rd,phantom ,"\eqref{sq:II}"] & C_{1+*}(I\times \Omega,\mathcal{H})\ar[d]\\
        C_{1+*}(I\times \Omega, \mathcal{H}\cup e^{-1}_\Omega(V_{x_0}^c)) \otimes C_*(\Omega) \ar[d, "\simeq"]  \ar[rd,phantom ,"\eqref{sq:III}"]\ar[r, "\Phi_*"] & C_{1+*}(I\times \Omega, \mathcal{H} \cup e^{-1}_\Omega(V_{x_0}^c)) \ar[d, "\simeq"]\\
        C_{1+*}(e^{-1}_\Omega(U_{x_0}),\mathcal{H} \cup e^{-1}_\Omega(V_{x_0}^c))\otimes C_*(\Omega) \ar[d, "e_\Omega^*\tau_{x_0}\cap \otimes \id "] \ar[r, "\Phi_*"]  \ar[rd,phantom ,"\eqref{sq:IV}"] & C_{1+*}(e^{-1}_\Omega(U_{x_0}),\mathcal{H} \cup e^{-1}_\Omega(V_{x_0}^c)) \ar[d, "e_\Omega^*\tau_{x_0}\cap "] \\
        C_{1-n+*}(e^{-1}_\Omega(U_{x_0}),\mathcal{H})\otimes C_*(\Omega) \ar[d, "R_*\otimes \id"] \ar[r, "\Phi_*"]  \ar[rd, phantom, "\eqref{sq:V}"]&C_{1-n+*}(e^{-1}_\Omega(U_{x_0}),\mathcal{H})\ar[d, "R_*"]\\
        C_{1-n+*}(\mathcal{F}\Omega,\mathcal{H}) \otimes C_*(\Omega)\ar[d, "\mathrm{c}\otimes \id"]\ar[r, "\Phi_*"]  \ar[rd, phantom, "\eqref{sq:VI}"] & C_{1-n+*}(\mathcal{F}\Omega,\mathcal{H}) \ar[d, "\mathrm{c}"]\\
        C_{1-n+*}((\Omega,x_0)^2)\otimes C_*(\Omega)\ar[r, "\Phi_*"] & C_{1-n+*}((\Omega,x_0)^2).
    \end{tikzcd}
\end{center}

\begin{enumerate}[(1)]
    \item\label{sq:I} This square commutes because the transitive lifting function on $I\times \Omega$ is defined by $\id\times \Phi$.
    \item\label{sq:II} This square commutes because it is given by a map of pairs of topological spaces with a compatible action.
    \item\label{sq:III} In Lemma \ref{lem:bary subd}, we prove that there exists a morphism of $\mathcal{A}_\infty$-modules making this square commute.
    \item\label{sq:IV} We construct a morphism of $\mathcal{A}_\infty$-modules making this square commute in Theorem \ref{thm:q is an Ainfty morphism}. This is the main goal of this section.
    \item\label{sq:V} There exists a morphism of $\mathcal{A}_\infty$-modules making this square commute by \cite[Theorem 5.8]{riegel2024chain} and the fact that $R\colon e_\Lambda^{-1}(U_\Delta)\to \mathcal{F}\Lambda$ is a map of fibrations.
    \item\label{sq:VI} This square commutes on the level of spaces.
\end{enumerate}

Summa summarum, each square has a (possibly strict) morphism of $\mathcal{A}_\infty$-modules making it commute. We denote their composition by
\begin{align}\label{eq:Ainfty nu}
    \boldsymbol{\nu}\colon C_*(\Omega, x_0)\otimes \to  C_{1-n+*}((\Omega,x_0)^2).
\end{align}
This morphism $\boldsymbol{\nu}$ is such that $\nu_1$ computes the Goresky-Hingston coproduct on the based loop space.

In this section, we construct the remaining morphisms of $\mathcal{A}_\infty$ for square \eqref{sq:III} and \eqref{sq:IV}. In fact giving a morphism for \eqref{sq:III} is done by the following lemma. The proof of which uses some general $\mathcal{A}_\infty$-theory we produce in Appendix \ref{app:inv qi of Ainfty modules}.

\begin{Lem}\label{lem:bary subd}
    For field coefficients, there exists a morphism of $\mathcal{A}_\infty$-modules 
    \begin{align*}
        \boldsymbol{\rho}\colon C_*(I\times \Omega, \mathcal{H}\cup e^{-1}_\Omega(V_{x_0}^c))\to C_*(I\times \Omega, \mathcal{H})
    \end{align*}
    such that $\rho_1$ is an inverse to the inclusion $C_*(I\times \Omega, \mathcal{H})\to C_*(I\times \Omega, \mathcal{H}\cup e^{-1}_\Omega(V_{x_0}^c))$ on homology.
\end{Lem}
\begin{proof}
    The inclusion $C_*(I\times \Omega, \mathcal{H})\to C_*(I\times \Omega, \mathcal{H}\cup e^{-1}_\Omega(V_{x_0}^c))$ is a quasi-isomorphism. Moreover, it is $C_*(\Omega)$-equivariant. Therefore, the inclusion defines an $\infty$-quasi-isomorphism (see Definition \ref{def:infty iso}). In the appendix in Corollary \ref{cor:inverting infty-quasi over field}, we show that over a field we can find for $\infty$-quasi-isomorphism an $\infty$-quasi-isomorphism which is the inverse in homology.
\end{proof}

It thus remains to produce the morphism of $\mathcal{A}_\infty$-modules
\begin{align*}
    \boldsymbol{q}\colon C_*(e_\Omega^{-1}(U_{x_0}),\mathcal{H}\cup e^{-1}_\Omega(V_{x_0}^c))\to C_{-n+*} (e_\Omega^{-1}(U_{x_0}),\mathcal{H})
\end{align*}
such that $q_1=e_\Omega^*\tau_{x_0}\cap$.\\

We first focus on the construction of $q_2$. The higher morphisms $q_k$ for $k\geq 3$ are given by appropriately composing the construction of $q_2$.

The map $q_2$ is a chain homotopy between $q_1\circ \Phi$ and $\Phi\circ (q_1\times \id)$.

\subsection{Rewriting the Cap Croduct}
To produce the chain homotopy $q_2$ it is convenient to rewrite $q_1\circ \Phi$ and $\Phi\circ (q_1\times \id)$ as in Lemma \ref{lem:rewritten q1 Phi}. To appropriately rewrite these two terms we need to introduce some notation.

\begin{Not}
    We note that by Properties \eqref{Theta property(iv)} and \eqref{Theta property(v)} of $\Theta$, the restriction $\Theta|_{U_{x_0}}$ defines a pointed, orientation-preserving embeddings $(U_{x_0},V_{x_0}^c)\hookrightarrow (U_{x_0},V_{x_0}^c)$. 
    
    We denote by $\Emb_*^+((U_{x_0},V_{x_0}^c),(U_{x_0},V_{x_0}^c))$ the space of all such embeddings and denote the maps
    \begin{align*}
        \Theta_{0}\colon \Omega &\to \Emb_*^+((U_{x_0},V_{x_0}^c),(U_{x_0},V_{x_0}^c)),\\
        \lambda&\mapsto \Theta(\lambda)|_{U_{x_0}}
    \end{align*}
    and
    \begin{align*}
        \eta\colon U_{x_0}\times \Emb^+_*((U_{x_0},V_{x_0}^c),(U_{x_0},V_{x_0}^c))&\to U_{x_0},\\
        (x,\varphi)&\mapsto \varphi(x),
    \end{align*}
    We denote the composition
    \begin{align*}
        \mu\colon \Emb_*^+(U_{x_0},U_{x_0})\times \Emb_*^+(U_{x_0},U_{x_0})&\to\Emb_*^+(U_{x_0},U_{x_0}),\\
        (\varphi_1,\varphi_2)&\mapsto \varphi_2\circ \varphi_1.
    \end{align*}
    We recall that we also denote $\mu\colon \Omega\times \Omega\to \Omega$ for the composition map. 
\end{Not}

This notation is chosen such that the following holds:
\begin{align*}
    \mu\circ (\Theta_0\times \Theta_0)&=\Theta_0\circ \mu,\\
    \eta \circ (\id\times \mu)&=\eta \circ (\eta \times \id),\\
    \Phi\circ (\id\times \mu)&=\Phi\circ (\Phi\times \id),
\end{align*}
and, by Lemma \ref{lem:ass trans lift fun} \eqref{item:ass trans lift fun (b)}, we have on $C_*((e^{-1}_\Omega(U_{x_0}),\mathcal{H}\cup e_\Omega^{-1}( V_{x_0}^c))\times \Omega)$
\begin{align}\label{eq:Theta0 and Phi commute}
    e_\Omega\circ\Phi=\eta\circ (e_\Omega\times \Theta_0).
\end{align}

We recall that the map $e_\Omega^*\tau_{x_0}\cap $ on $\alpha\in C_*(e^{-1}_\Omega(U_{x_0}),\mathcal{H}\cup e_\Omega^{-1}( V_{x_0}^c))$ is defined as the composition
\begin{align*}
    ((\tau_{x_0}\circ e_\Omega)\otimes \id)\circ \Delta_*(\alpha).
\end{align*}
Here $\tau_{x_0}$ is interpreted as a degree $-n$ map $C^n(U_{x_0},V_{x_0}^c)\to \R$ and $\Delta_*$ is a diagonal approximation. We fix the \textit{Serre diagonal} as defined in \cite{serre1951homologie} given by
\begin{align*}
    \Delta_*\colon C_n(X)&\to C_*(X)\otimes C_*(X),\\
    \sigma &\mapsto \sum_{J\sqcup K=\{1,\dots, n\}} \sgn(J,K) \sigma|_{I^J\times \{0\}^K}\otimes \sigma|_{\{1\}^J\times I^K}
\end{align*}
for $\sgn(J,K)=(-1)^{|\{j>k\mid j\in J,k\in K\}|}$.

The last piece of notation is the following:
using the exponential map at $x_0$, we get an inner product on $U_{x_0}$. We denote by
\begin{align*}
    \SOEmb(U_{x_0},U_{x_0})\subseteq \Emb_*^+((U_{x_0},V_{x_0}^c),(U_{x_0},V_{x_0}^c))
\end{align*}
the subspace of embeddings which leave this scalar product invariant. This space is a deformation retract of $\Emb_*^+((U_{x_0},V_{x_0}^c),(U_{x_0},V_{x_0}^c))$. We denote a retract by 
\begin{align*}
    r\colon \Emb_*^+((U_{x_0},V_{x_0}^c),(U_{x_0},V_{x_0}^c)) \to \SOEmb((U_{x_0},V_{x_0}^c),(U_{x_0},V_{x_0}^c)).
\end{align*}
In Lemma \ref{lem:H1 for Omega}, we construct a particularly nice retract.

We are now in a position to rewrite the two maps that we want to give a homotopy between:

\begin{Lem}\label{lem:rewritten q1 Phi}
    Let $\tau_{x_0}\in C^n(U_{x_0},V_{x_0}^c)$ be a $\SOEmb(U_{x_0},U_{x_0})$-invariant Thom class, that is, it holds as maps $C_*(U_{x_0},V_{x_0}^c)\otimes C_*(\SOEmb(U_{x_0},U_{x_0}))\to \R$ that 
    \begin{align}\label{eq:inv taux0}
        \tau_{x_0}\circ \eta\circ \times =\tau_{x_0}\circ  \mathrm{pr}_1\circ \times.
    \end{align}
    
    On $\alpha\otimes \lambda\in C_{*}(e_\Omega^{-1}(U_{x_0}),\mathcal{H} \cup e_\Omega^{-1}(x_0)^c)\otimes C_*(\Omega)$, 
    \begin{enumerate}[(a)]
        \item\label{item:rewritten a} the map $q_1\circ \Phi$ is computed by
        \begin{align*}
            ((\tau_{x_0}\circ \eta\circ (e_\Omega\times \Theta_0))\otimes\Phi)\circ \Delta_*(\alpha\times \lambda);
        \end{align*}
        \item\label{item:rewritten b} the map $\Phi\circ (q_1\otimes \id)$ is computed by
        \begin{align*}
            (\left(\tau_{x_0}\circ \eta\circ (e_\Omega\times (r\circ\Theta_0))\right)\otimes \Phi) \circ \Delta_* (\alpha\times \lambda).
        \end{align*}
    \end{enumerate}
\end{Lem}
\begin{Rem}
    In Corollary \ref{cor:invariant taux0}, we construct a Thom class with real coefficients such that \eqref{eq:inv taux0} holds. We postpone the construction to Section \ref{sec:inv Thom class}, because for our purposes it is helpful to create a Thom class $\tau_\Delta\in C^n(U_\Delta,V_\Delta^c)$ with an analogous property. The class $\tau_{x_0}$ as in the previous lemma is then given by the restriction of $\tau_\Delta$ to $U_{x_0}$.
\end{Rem}

\begin{Rem}\label{rem:reduced to homotopies of r}
    Lemma \ref{lem:rewritten q1 Phi} is key because it shows that constructing $q_2$, a homotopy between the two maps $q_1\circ \Phi$ and $\Phi\circ (q_1\otimes \id)$, comes down to giving a homotopy between $\id$ and $r$ on $\Emb_*^+((U_{x_0},V_{x_0}^c),(U_{x_0},V_{x_0}^c))$.
\end{Rem}

\begin{proof}[Proof of Lemma \ref{lem:rewritten q1 Phi}]
    We first prove \eqref{item:rewritten a}:
    \begin{align*}
        q_1\circ \Phi(\alpha\otimes \lambda)&= ((\tau_{x_0}\circ e_\Omega)\otimes \id)\circ \Delta_*(\Phi(\alpha\otimes \lambda))\\
        &=((\tau_{x_0}\circ e_\Omega\circ \Phi)\otimes \Phi)\circ \Delta_*(\alpha\times \lambda)\\
        &\overset{\eqref{eq:Theta0 and Phi commute}}{=}((\tau_{x_0}\circ \eta\circ (e_\Omega\times \Theta_0))\otimes\Phi)\circ \Delta_*(\alpha\times \lambda).
    \end{align*}
    For \eqref{item:rewritten b}, we set $a:=|\alpha|$ and $l:=|\lambda|$. We compute:
    \begin{align*}
        \Phi\circ (q_1\otimes \id)(\alpha\otimes \lambda)&= \Phi(((\tau_{x_0} \circ e_\Omega )\otimes \id) \Delta_*)\otimes \id)(\alpha\otimes \lambda)\\
        &=\Phi\left(\left(\sum_{J\sqcup K=\{1,\dots a\}}  \sgn(J,K) \tau_{x_0} (e_\Omega(\alpha|_{I^J\times \{0\}^K}))\alpha|_{\{1\}^J\times I^K}\right)\otimes \lambda\right)\\
        &=\sum_{\substack{J\sqcup K=\{1,\dots a+l\}\\J\subseteq \{1,\dots,a\}}}  \sgn(J,K) \tau_{x_0} (e_\Omega(\mathrm{pr}_1((\alpha\times \lambda)|_{I^J\times \{0\}^K})))(\Phi(\alpha\otimes\lambda))|_{\{1\}^J\times I^K}.
    \end{align*}
    We now note that $\mathrm{pr}_1((\alpha\times \lambda)|_{I^J\times \{0\}^K})$ is a degenerate chain if $J\cap \{a+1,\dots,a+l\}\neq \varnothing$. We can thus drop the condition $J\subseteq\{1,\dots,a\}$ in the sum because the added terms are zero. We therefore find
    \begin{align*}
        \Phi\circ (q_1\otimes \id)(\alpha\otimes \lambda)&=\sum_{\substack{J\sqcup K=\{1,\dots a+l\}}}   \sgn(J,K) \tau_{x_0} (e_\Omega(\mathrm{pr}_1((\alpha\times \lambda)|_{I^J\times \{0\}^K})))(\Phi(\alpha\otimes\lambda))|_{\{1\}^J\times I^K}\\
        &=(\left(\tau_{x_0}\circ e_\Omega\circ \mathrm{pr}_1\right)\otimes \Phi )\circ \Delta_* (\alpha\times \lambda)
    \end{align*}
    We moreover have
    \begin{align*}
        (\left(\tau_{x_0}\circ e_\Omega\circ \mathrm{pr}_1\right)\otimes \Phi )\circ \Delta_* (\alpha\times \lambda)&=(\left(\tau_{x_0}\circ \mathrm{pr}_1\circ (e_\Omega\times r\Theta_0)\right)\otimes \Phi )\circ \Delta_* (\alpha\times \lambda).
    \end{align*}
    By our assumption on $\tau_{x_0}$ as in \eqref{eq:inv taux0} and the fact that $r\Theta_0$ lands in $\SOEmb(U_{x_0},U_{x_0})$, we find
    \begin{align*}
        (\left(\tau_{x_0}\circ \mathrm{pr}_1\circ (e_\Omega\times r\Theta_0)\right)\otimes \Phi )\circ \Delta_* (\alpha\times \lambda)
        &=(\left(\tau_{x_0}\circ \eta\circ (e_\Omega\times r \Theta_0)\right)\otimes \Phi) \circ \Delta_* (\alpha\times \lambda)
    \end{align*}
    which is what we wanted to show.
\end{proof}

\subsection{Higher Homotopies of Embeddings}\label{subsec:first hpty}
We now construct a homotopy $H_1$ between the identity and the retraction $r\colon \Emb_*^+((U_{x_0},V_{x_0}^c),(U_{x_0},V_{x_0}^c))\to \SOEmb(U_{x_0},U_{x_0})$. As explained in Remark \ref{rem:reduced to homotopies of r}, this lets us construct $q_2$. Finally, composing the homotopy $H_1$ produces the higher maps $q_k$ for $k\geq 3$.

\begin{Lem}\label{lem:H1 for Omega}
    There exists a homotopy
    \begin{align*}
        H_1\colon I\times \Emb_*^+((U_{x_0},V_{x_0}^c),(U_{x_0},V_{x_0}^c))\to \Emb_*^+((U_{x_0},V_{x_0}^c),(U_{x_0},V_{x_0}^c))
    \end{align*}
    such that, for all $\varphi,\varphi_1,\varphi_2\in \Emb_*^+((U_{x_0},V_{x_0}^c),(U_{x_0},V_{x_0}^c))$ and $t\in I$, it holds
    \begin{enumerate}[(i)]
        \item\label{item:H1 i} $H_1(0,\varphi)=:r(\varphi)\in \SOEmb(U_{x_0},U_{x_0})$;
        \item $H_1(1,\varphi)=\varphi$;
        \item\label{item:H1 iii} $H_1(t, r(\varphi_1)\circ \varphi_2)= r(\varphi_1)\circ H_1(t,\varphi_2)$.
    \end{enumerate}
\end{Lem}
\begin{proof}
    We fix some orthonormal coordinates on $U_{x_0}$. We may thus identify $U_{x_0}\cong B_\varepsilon(0)\subseteq \R^n$ such that $\|x\|=d_M(x,x_0)$. We then define
    \begin{align*}
        \alpha\colon (U_{x_0},V_{x_0}^c)&\cong (\R^n,B_1(0)^c),\\
        x&\mapsto \frac{x}{\varepsilon-\|x\|}.
    \end{align*}
    This in particular identifies $\Emb_*^+((U_{x_0},V_{x_0}^c),(U_{x_0},V_{x_0}^c))$ with $\Emb_*^+((\R^n,B_1(0)^c),((\R^n,B_1(0)^c))$ and their respective subspaces $\SOEmb(U_{x_0},U_{x_0})$ and $\SO(n)$.
    
    In Proposition \ref{prop:hpt of rel emb}, we define a homotopy $H$ that exhibits $\SO(n)$ as a deformation retract of $\Emb_*^+((\R^n,B_1(0)^c),((\R^n,B_1(0)^c))$ and commutes with left multiplication of $\SO(n)$. The map
    \begin{align*}
        H_1(t,\varphi):=\alpha^{-1}\circ H(t,\varphi\circ \alpha)
    \end{align*}
    defines the desired homotopy. The invariance of $H$ under left $\SO(n)$-multiplication moreover shows that $H_1$ is independent of our choice of $\alpha$, because two different choices of $\alpha$ differ by an element in $\SO(n)$.
\end{proof}

We denote the corresponding chain homotopy by
\begin{align*}
    h_1\colon C_*(\Emb_*^+((U_{x_0},V_{x_0}^c),(U_{x_0},V_{x_0}^c)))&\to C_{1+*}(\Emb_*^+((U_{x_0},V_{x_0}^c),(U_{x_0},V_{x_0}^c))),\\
    \varphi&\mapsto H_1(I\times \varphi).
\end{align*}
We define inductively for $k\geq 2$
\begin{align*}
    H_k\colon I^k\times \Emb^+_*((U_{x_0},V_{x_0}^c),(U_{x_0},V_{x_0}^c))^{\times k}&\to \Emb_*^+((U_{x_0},V_{x_0}^c),(U_{x_0},V_{x_0}^c)),\\
    (t_1,\dots, t_k,\varphi_1,\dots,\varphi_k)& \mapsto H_1(t_1,\mu(\varphi_1,H_{k-1}(t_2,\dots,t_k,\varphi_2,\dots,\varphi_k)))
\end{align*}
and
\begin{align*}
    h_k\colon C_*(\Emb^+_*((U_{x_0},V_{x_0}^c),(U_{x_0},V_{x_0}^c))^{\times k})&\to C_{k+*}(\Emb^+_*((U_{x_0},V_{x_0}^c),(U_{x_0},V_{x_0}^c))),\\
    \varphi&\mapsto H_k(I^k\times \varphi).
\end{align*}
The higher homotopies are defined in a coherent way such that they satisfy the following useful formulas:
\begin{Lem}\label{lem:h_N is homotopy}
    The map $h_k$ commutes with $r$ in the following sense
    \begin{align}\label{eq:h_k and r commute}    
        h_k\circ (\id^{\times k-1}\times \mu)\circ (\id^{\times k}\times r)&=\mu\circ (h_k\times r).
    \end{align}
    The maps $h_k$ are such that
    \begin{align*}
        h_1d+d h_1=\id -r
    \end{align*}
    and for $k\geq 2$
    \begin{align*}
        h_kd +(-1)^{k+1}d h_k=&\ (-1)^{k+1}(\mu\circ (\id \times h_{k-1})-r\circ \mu\circ (\id \times h_{k-1}))\\
       &+ (-1)^{k+1}\sum_{i=2}^{k}(-1)^{i+1}h_{k-1}\circ (\id^{i-2}\times \mu\times \id^{k-i})\\
       &+(-1)^{k+1}\sum_{i=2}^{k-1}(-1)^{i}\mu\circ (h_{i-1}\times (r\circ\mu \circ (\id \times  h_{k-i})))\\
       &-\mu \circ (h_{k-1}\times r)).
    \end{align*}
\end{Lem}
\begin{proof}
    Equation \eqref{eq:h_k and r commute} holds for $k=1$ by Lemma \ref{lem:H1 for Omega} \eqref{item:H1 iii}. It follows by induction for $k\geq 2$.\\
    
    We evaluate on $\varphi_1\times \dots \times \varphi_k\in C_*( \Emb_*^+((U_{x_0},V_{x_0}^c),(U_{x_0},V_{x_0}^c))^{\times k})$:
    \begin{align*}
        ( h_kd+(-1)^{k+1}d& h_k)(\varphi_1\times \dots \times \varphi_k)=(-1)^{k+1}  H_k(d(I^k)\times \varphi_1\times \dots \times \varphi_k).
    \end{align*}
    For $k=1$, we compute:
    \begin{align*}
        h_1d&(\varphi_1) +d h_1(\varphi_1)=H_1(1,\varphi_1)-H_1(0,\varphi_1)=\varphi_1-r(\varphi_1).
    \end{align*}
    For $k\geq 2$, we find:
    \begin{align}\label{eq:hkd+dhk}
        \begin{split}
            (-1)^{k+1} & H_k(d(I^k)\times \varphi_1\times \dots \times \varphi_k)\\
            =&\ (-1)^{k+1}\sum_{i=1}^k (-1)^{i-1}H_k(I^{i-1}\times \{1\}\times I^{k-i}\times (\varphi_1\times \dots \times \varphi_k))\\
            &-(-1)^{k+1}\sum_{i=1}^k (-1)^{i-1}H_k(I^{i-1}\times \{0\}\times I^{k-i}\times (\varphi_1\times \dots \times \varphi_k)).
        \end{split}
    \end{align}
    For $i=1$, we compute:
    \begin{align*}
        H_k(\{1\}\times I^{k-1} \times\varphi_1\times \dots \times \varphi_k)&=\mu\circ (\id \times h_{k-1})(\varphi_1\times \dots \times \varphi_k)
    \end{align*}
    and
    \begin{align*}
        H_k(\{0\}\times I^{k-1} \times\varphi_1\times \dots \times \varphi_k)&=r\circ \mu\circ (\id \times h_{k-1})(\varphi_1\times \dots \times \varphi_k).
    \end{align*}
    For $i=k$, we have
    \begin{align*}
        H_k(I^{k-1}\times \{1\}\times \varphi_1\times \dots \times \varphi_k)&=H_{k-1}(I^{k-1} \times \varphi_{1}\times \dots \times \varphi_{k-2}\times (\varphi_k\circ \varphi_{k-1}))\\
        &=h_{k-1}\circ (\id^{k-2}\times \mu)(\varphi_1\times \dots \times \varphi_k)
    \end{align*}
    and 
    \begin{align*}
        H_k(I^{k-1}\times \{0\}\times \varphi_1\times \dots \times \varphi_k)&=H_{k-1}(I^{k-1} \times \varphi_{1}\times \dots \times \varphi_{k-2}\times (r(\varphi_k)\circ \varphi_{k-1}))\\
        &=h_{k-1}\circ (\id^{k-2}\times \mu)\circ (\id^{k-1}\times r)(\varphi_1\times \dots \times \varphi_k)\\
        &\overset{\eqref{eq:h_k and r commute}}{=}\mu \circ (h_{k-1}\times r)(\varphi_1\times \dots \times \varphi_k).
    \end{align*}
    For $1<i<k$, we have
    \begin{align*}
        H_k(I^{i-1}\times \{1\}&\times I^{k-i}\times \varphi_1\times \dots \times \varphi_k)\\
        &=H_{i-1}(I^{i-1}\times \varphi_1\times \dots \times \varphi_{i-2 } \times (H_{k-i}(I^{k-i}\times \varphi_{i+1}\times \dots \times \varphi_k)\circ \varphi_i\circ \varphi_{i-1}))\\
        &=H_{k-1}(I^{k-1}\times \varphi_1\times \dots \times \varphi_{i-2}\times (\varphi_i\circ \varphi_{i-1})\times \varphi_{i+1}\times \dots\times\varphi_k)\\
        &=h_{k-1}\circ (\id^{i-2}\times \mu\times \id^{k-i})(\varphi_1\times \dots \times \varphi_k)
    \end{align*}
    and 
    \begin{align*}
         H_k(I^{i-1}&\times \{0\}\times I^{k-i}\times \varphi_1\times \dots \times \varphi_k)\\
         &=H_{i-1}(I^{i-1}\times \varphi_1\times \dots \times \varphi_{i-2 } \times (r(H_{k-i}(I^{k-i}\times \varphi_{i+1}\times \dots \times \varphi_k)\circ \varphi_i)\circ \varphi_{i-1}))\\
         &=h_{i-1}(\id^{i-2}\times (\mu\circ (\id \times (r\circ \mu \circ (\id \times h_{k-i})))))(\varphi_1\times \dots \times \varphi_k)\\
         &\overset{\eqref{eq:h_k and r commute}}{=}\mu\circ (h_{i-1}\times (r\circ\mu \circ (\id \times  h_{k-i})))(\varphi_1\times \dots \times \varphi_k)
    \end{align*}
    Putting this into \eqref{eq:hkd+dhk} then shows the lemma.
\end{proof}

\subsection{The Higher Morphisms}
We now define the morphisms $q_k$ for $k\geq 2$. Roughly speaking, they arise from plugging the chain homotopies $h_k$ into Lemma \ref{lem:rewritten q1 Phi}.\\

We define for $k\geq 1$
\begin{align}\label{eq:def of qk}
    q_{k+1}(\alpha\otimes \lambda_1\otimes \dots \otimes  \lambda_k):=((\tau_{x_0} \circ \eta  \circ (e_\Omega\times ( h_k\circ \Theta_0^{\times k})))\otimes \Phi)\circ \Delta_*(\alpha\times \lambda_1\times \dots\times \lambda_k).
\end{align}
Here $\Phi\colon e_\Omega^{-1}(U_{x_0})\times \Omega\to e_\Omega^{-1}(U_{x_0})$ is extended to a map 
\begin{align*}
    e_\Omega^{-1}(U_{x_0})\times  \Omega^{\times k}&\to e_\Omega^{-1}(U_{x_0}),\\
    (x,\lambda_1,\dots,\lambda_k)& \mapsto \Phi(x,\lambda_1\cdots \lambda_k).
\end{align*}

\begin{Thm}\label{thm:q is an Ainfty morphism}
    The tuple $q_k$ for $k\geq 1$ defines a morphism of $\mathcal{A}_\infty$-modules
    \begin{align*}
        \boldsymbol{q}\colon C_*(e_\Omega^{-1}(U_{x_0}),\mathcal{H} \cup e_\Omega^{-1}(V_{x_0})^c)\to C_{-n+*}(e_\Omega^{-1}(U_{x_0}),\mathcal{H})
    \end{align*}
    with respect to the action given by $\Phi$ a transitive lifting function induced by a transitive diffeomorphism lift.
\end{Thm}
\begin{proof}
    We note that $\Phi$ defines strictly associative $C_*(\Omega)$-module structure on source and target. Moreover, since the map is not degree-preserving we need the differential $(-1)^nd$ on the target. Therefore the condition for $\boldsymbol{q}$ being a morphism of $\mathcal{A}_\infty$-modules reads for $N\geq 0$ as:
    \begin{equation}\label{eq:Ainfty-eq}
        \begin{split}
            q_{N+1}&\circ d+(-1)^{N+1+n}d \circ q_{N+1}\\
            &=(-1)^{N+1}q_N\circ (\Phi\otimes \id ^{\otimes N-1})-\Phi\circ (q_N\otimes \id)+\sum_{r=1}^{N-1}(-1)^{N+1+r}q_N\circ (\id^{\otimes r}\otimes \mu\otimes \id ^{\otimes {N-1-r}})
        \end{split} 
    \end{equation}
    For $N=0$, this states that $q_1$ is a chain map, which is already proven.

    For $N\geq 1$, we recall that 
    \begin{align*}
        q_{N+1}
        =&\ ((\tau_{x_0} \circ \eta  \circ (e_\Omega\times h_N\Theta_0^{\times N}))\otimes \Phi)\circ \Delta_*\circ\times^N.
    \end{align*}
    Except $h_N$, all maps in this composition are chain maps and only $\tau_{x_0}$ and $h_N$ are not degree-preserving. We thus have
    \begin{equation}\label{eq:get homotopy inside}
        \begin{split}
            q_{N+1}d+(-1)^{N+1+n}d q_{N+1}=((\tau_{x_0} \circ \eta  \circ (e_\Omega\times (h_Nd +(-1)^{N+1}d h_N)\Theta_0^{\times N}))\otimes \Phi)\circ \Delta_*\circ \times^N.
        \end{split}
    \end{equation}
    We can thus apply Lemma \ref{lem:h_N is homotopy} to $h_Nd +(-1)^{N+1}d h_N$.
    
    For $N=1$, we apply Lemma \ref{lem:rewritten q1 Phi} and find:
    \begin{align*}
        q_{2}(d (\alpha\otimes \lambda_1))+(-1)^{n} d q_{2}(\alpha\otimes \lambda_1)=&\ ((\tau_{x_0} \circ \eta  \circ (e_\Omega\times \Theta_0))\otimes \Phi)\circ \Delta_*(\alpha\times \lambda_1)\\
        &- ((\tau_{x_0} \circ \eta  \circ (e_\Omega\times (r\circ \Theta_0)))\otimes \Phi)\circ \Delta_*(\alpha\times \lambda_1)\\
        =&\  q_1(\Phi(\alpha\otimes \lambda))-\Phi(q_1(\alpha)\otimes \lambda).
    \end{align*}
    This shows the Equation \eqref{eq:Ainfty-eq} for $N=1$.\\
    
    For $N\geq 2$, we find
    \begin{align*}
        q_{N+1}&\circ d+(-1)^{N+1+n} d\circ  q_{N+1}\\
        =&\ (-1)^{N+1}((\tau_{x_0} \circ \eta  \circ (e_\Omega\times (\mu\circ (\id \times h_{N-1})\circ \Theta_0^{\times N})))\otimes \Phi)\circ \Delta_*\circ \times^N\\
        &+(-1)^{N}((\tau_{x_0} \circ \eta  \circ (e_\Omega\times (r\circ \mu\circ (\id \times h_{N-1})\circ \Theta_0^{\times N})))\otimes \Phi)\circ \Delta_*\circ \times^N)\\
        &+ (-1)^{N+1}\sum_{i=2}^{N}(-1)^{i+1}((\tau_{x_0} \circ \eta  \circ (e_\Omega\times (h_{N-1}\circ (\id^{i-2}\times \mu\times \id^{N-i})\circ \Theta_0^{\times N})))\otimes \Phi)\circ \Delta_*\circ \times^N\\
        &+(-1)^{N+1}\sum_{i=2}^{N-1}(-1)^{i}((\tau_{x_0} \circ \eta  \circ (e_\Omega\times (\mu\circ (h_{i-1}\times (r\mu  (\id \times  h_{N-i})))\circ \Theta_0^{\times N})))\otimes \Phi)\circ \Delta_*\circ \times^N\\
        &-((\tau_{x_0} \circ \eta  \circ (e_\Omega\times ( (\mu \circ (h_{N-1}\times r)))\circ \Theta_0^{\times N}))\otimes \Phi)\circ \Delta_*\circ \times^N.
    \end{align*}
    The following lemma thus proves \eqref{eq:Ainfty-eq} for $N\geq 2$:
    \begin{Lem}\label{lem:five lemma}
        The following five equations hold:
        \begin{enumerate}[(a)]
            \item\label{item:five lemma a} \begin{align*}
                ((\tau_{x_0} & \circ \eta  \circ (e_\Omega\times (\mu \circ (\id \times h_{N-1})\circ \Theta_0^{\times N}))))\otimes \Phi)\circ \Delta_*\circ \times^N\\
                &=q_N(\Phi\otimes \id^{\otimes N-1});
            \end{align*}
            \item\label{item:five lemma b} \begin{align*}
                (\tau_{x_0} \circ \eta  \circ (e_\Omega\times( r\circ  \mu \circ (\id \times h_{N-1})\circ \Theta_0^{\times N})))\otimes \Phi)\circ \Delta_*\circ \times^N=0;
            \end{align*}
            \item\label{item:five lemma c} \begin{align*}
                \sum_{i=2}^{N}& (-1)^{i+1}((\tau_{x_0} \circ \eta  \circ (e_\Omega\times( h_{N-1}(\id^{\otimes i-2}\otimes  \mu \otimes \id^{\otimes N-i})\circ \Theta_0^{\times N})))\otimes \Phi)\circ \Delta_*\circ \times^N\\
                &=\sum_{r=1}^{N-1}(-1)^{r} q_N\circ (\id^{\otimes r}\otimes \mu \otimes \id^{\otimes N-1-r});
            \end{align*}
            \item\label{item:five lemma d} for $2\leq i\leq N-1$
            \begin{align*}
                ((\tau_{x_0} \circ \eta  \circ (e_\Omega\times (\mu(h_{i-1}\otimes (r\circ \mu\circ (\id \otimes h_{N-i})))\circ \Theta_0^{\times N})))\otimes \Phi)\circ \Delta_*\circ \times^N=0;
            \end{align*}
            \item\label{item:five lemma e} \begin{align*}
                ((\tau_{x_0} &\circ \eta  \circ (e_\Omega\times( \mu \circ (h_{N-1}\times r)\circ \Theta_0^{\times N}))))\otimes \Phi)\circ \Delta_*\circ \times^N=\Phi(q_N\otimes \id).
            \end{align*}
        \end{enumerate}
    \end{Lem}
    \begin{proof}[Proof of Lemma \ref{lem:five lemma}]
        For \eqref{item:five lemma a}, we use that $\eta $ and $\mu$ commute:
        \begin{align*}
            ((\tau_{x_0}& \circ \eta  \circ (e_\Omega\times (\mu \circ (\id \times h_{N-1})\circ \Theta_0^{\times N}))))\otimes \Phi)\circ \Delta_*\circ \times^N\\
            &=((\tau_{x_0} \circ \eta  \circ ((\eta \times \id)\circ   (\id\times \id \times h_{N-1})\circ (e_\Omega\times\Theta_0^{\times N}))\otimes \Phi)\circ \Delta_*\circ \times^N\\
            &=((\tau_{x_0} \circ \eta  \circ ( e_\Omega \times h_{N-1})\circ (\Phi\times\Theta_0^{\times N-1}))\otimes \Phi)\circ \Delta_*\circ \times^N\\
            &=((\tau_{x_0} \circ \eta  \circ ( e_\Omega \times h_{N-1})\circ (\id\times\Theta_0^{\times N-1}))\otimes \Phi)\circ \Delta_*\circ (\Phi \times\id^{\times N-1})\circ  \times^N\\
            &=((\tau_{x_0} \circ \eta  \circ ( e_\Omega \times (h_{N-1}\circ \Theta_0^{\times N-1})))\otimes \Phi)\circ \Delta_*\circ  \times^{N-1}\circ (\Phi\otimes \id^{\otimes N-1})\\
            &= q_N(\Phi\otimes \id^{\otimes N-1}).
        \end{align*}
        For \eqref{item:five lemma b}, we use $\SOEmb(U_{x_0},U_{x_0})$-invariance of $\tau_{x_0}$:
        \begin{align*}
            ((\tau_{x_0} &\circ \eta  \circ (e_\Omega\times( r\circ  \mu \circ (\id \times h_{N-1})\circ \Theta_0^{\times N})))\otimes \Phi)\circ \Delta_*\circ \times^N\\
            &=(((\tau_{x_0}   \circ \mathrm{pr}_1\circ (e_\Omega\times( r\circ  \mu \circ (\id \times h_{N-1})\circ \Theta_0^{\times N}))))\otimes \Phi)\circ \Delta_*\circ \times^N.
        \end{align*}
        The map $\mathrm{pr}_1\circ (e_\Omega\times( r\circ  \mu \circ (\id \times h_{N-1})))$ has as image only degenerate chains as it is independent of the coordinates in $h_{N-1}$ and $N-1>0$. Because degenerate chains are zero this proves \eqref{item:five lemma b}.
        
        For \eqref{item:five lemma c}, we compute for $2\leq i\leq N$:
        \begin{align*}
            ((\tau_{x_0} &\circ \eta  \circ (e_\Omega\times h_{N-1}((\id^{\otimes i-2}\otimes \mu \otimes \id^{\otimes N-i})\circ \Theta_0^{\times N})))\otimes \Phi)\circ \Delta_*\circ\times^N\\
            &=((\tau_{x_0} \circ \eta  \circ (e_\Omega\times h_{N-1}( \Theta_0^{\times N}\circ (\id^{\otimes i-2}\otimes \mu \otimes \id^{\otimes N-i}))))\otimes \Phi)\circ \Delta_*\circ \times^N\\
            &=((\tau_{x_0} \circ \eta  \circ (e_\Omega\times (h_{N-1} \Theta_0^{\times N})))\otimes \Phi)\circ \Delta_*\circ (\id^{\otimes i-1}\otimes \mu \otimes \id^{\otimes N-i})\circ \times^N\\
            &=q_{N}\circ (\id^{\otimes i-1}\otimes \mu \otimes \id^{\otimes N-i}).
        \end{align*}
        The variable change $r=i+1$ then shows \eqref{item:five lemma c}.
        
        For \eqref{item:five lemma d}, we compute for $2\leq i\leq N-1$ using that $\mu$ and $\eta $ commute by associativity of the evaluation:
        \begin{align*}
            ((\tau_{x_0} &\circ \eta  \circ (e_\Omega\times (\mu(h_{i-1}\otimes (r \mu (\id \otimes h_{N-i})))\circ \Theta_0^{\times N})))\otimes \Phi)\circ \Delta_*\circ \times^N \\
            &=((\tau_{x_0} \circ \eta  \circ (\eta \times r)\circ (\id\times \id\times \mu )\circ  ( \id \times h_{i-1} \times \id  \times h_{N-i})\circ (  e_\Omega\times \Theta_0^{\times N}))\otimes \Phi)\circ \Delta_*\circ \times^N.
        \end{align*}
        By $\SOEmb(U_{x_0},U_{x_0})$-invariance of $\tau_{x_0}$, we have
        \begin{align*}
            \tau_{x_0} \circ \eta  \circ (\eta \times r)=\tau_{x_0}\circ \eta\circ \mathrm{pr}_1\circ (\eta\times r).
        \end{align*}
        This is zero because the image of $\mathrm{pr}_1\circ  (\eta \times \id)\circ (\id\times \id\times \mu )\circ  ( \id \times h_{i-1} \times \id  \times h_{N-i})$ is independent of the coordinates of $\mu\circ (\id \otimes h_{N-i})$ and is thus a degenerate chain.
        
        For \eqref{item:five lemma e}, we use again that $\eta $ and $\mu$ commute:
        \begin{align*}
            ((\tau_{x_0} &\circ \eta  \circ (e_\Omega\times( \mu \circ (h_{N-1}\times r)\circ \Theta_0^{\times N}))))\otimes \Phi)\circ \Delta_*\circ \times^N\\
            &=((\tau_{x_0} \circ \eta  \circ (\eta \times r)\circ ( \id\times h_{N-1}\times \id) \circ (e_\Omega\times \Theta_0^{\times N}))\otimes \Phi)\circ \Delta_*\circ \times^N\\
            &=((\tau_{x_0} \circ \eta  \circ \mathrm{pr}_1\circ  (\eta \times r)\circ ( \id\times h_{N-1}\times \id) \circ (e_\Omega\times \Theta_0^{\times N}))\otimes \Phi)\circ \Delta_*\circ \times^N.
        \end{align*}
        An analogous computation as in the proof of Lemma \ref{lem:rewritten q1 Phi} \eqref{item:rewritten b} shows that
        \begin{align*}
            ((\tau_{x_0}& \circ \eta  \circ  \mathrm{pr}_1\circ  (\eta \times r)\circ ( \id\times h_{N-1}\times \id) \circ (e_\Omega\times \Theta_0^{\times N}))\otimes \Phi)\circ \Delta_*\circ \times^N=\Phi(q_N\otimes \id).
        \end{align*}
        This proves the lemma and thus Theorem \ref{thm:q is an Ainfty morphism}.
    \end{proof}
\end{proof}

\section{Coherent Chain Homotopies}\label{sec:cch}
In this section, we give a proof strategy how one can prove that a map $\widetilde{\eta}\colon C_*(M,C_*(F_1,F_1'))\to C_*(M,C_{m+*}(F_2, F_2'))$ models a certain map $\eta\colon C_*(E_1, E_1')\to C_*(E_2,E_2')$. We follow the proof strategy from \cite[Theorem 5.8]{riegel2024chain}. Riegel constructs a collection of maps called coherent homotopies (\cite[Definition 5.10]{riegel2024chain}) that witness that a map $\widetilde{\varphi} \colon C_*(M,C_*(F_1))\to C_*(M,C_*(F_2))$ coming from a morphism of fibration $\varphi\colon E_1\to E_2$ models $\varphi_*\colon C_*(E_1)\to C_*(E_2)$. This is in contrast to our map where we are interested in maps that do not come from maps of fibrations. We thus isolate the algebraic part of his argument: we define in Definition \ref{def:coherent chain hpty} coherent chain homotopies as an algebraic version of Riegel's coherent homotopies. In Proposition \ref{prop:coh chain hpty witness}, we show that given a coherent chain homotopy, $\widetilde{\eta}$ models $\eta$ on $C_*(M,C_*(F_1,F_1'))$. Finally in Subsection \ref{subsec:inv cch}, we line out how to invert certain coherent chain homotopies.\\

We define coherent chain homotopies as algebraic analogues of Riegel's coherent homotopies. These witness that a morphism $\widetilde{\eta}\colon C_*(M,C_*(F_1))\to C_*(M,C_{m+*}(F_2))$ models a map $\eta\colon C_*(E_1)\to C_{m+*}(E_2)$. Riegel proves analogous statements in the case where $\eta=\varphi_*$ for a map of fibrations $\varphi\colon E_1\to E_2$. We generalize to the setting where $\eta$ is a purely algebraic map which might not be degree-preserving.

\begin{Def}\label{def:coherent chain hpty}
    For $i=1,2$, let $\pi_i\colon E_i\to M$ be two fibrations with fibres $F_i$ and transitive lifting functions $\Phi_i$. Let $\eta_1:=\eta\colon C_*(E_1)\to C_{m+*}(E_2)$ be a chain map of degree $m\in \Z$.
    
    A \textit{coherent chain homotopy $\boldsymbol{\eta}\colon C_*(E_1)\to C_{m+*}(E_2)$ for $\eta$} is a collection of maps
    \begin{align*}
        \eta_k\colon C_*(F_1)\otimes C_*(\Omega)^{\otimes (k-2)}\otimes C_*(\mathcal{P}_{x_0\to M}M)\to C_{k-1+m+*}(E_2)
    \end{align*}
    for $k\geq 2$ that satisfy the $\mathcal{A}_\infty$-equations
    \begin{align}\label{eq:Ainfty coh chain hpty}
        \begin{split}
            \eta_{N+1}&\circ d+(-1)^{N+1+m}d \circ \eta_{N+1}\\
            &=(-1)^{N+1}\eta_N\circ (\Phi_1\otimes \id ^{\otimes N-1})-\Phi_2\circ (\eta_N\otimes \id)+\sum_{r=1}^{N-1}(-1)^{N+1+r}\eta_N\circ (\id^{\otimes r}\otimes \mu\otimes \id ^{\otimes {N-1-r}})
        \end{split}
    \end{align}
    for all $N\geq 1$ and restrict to maps
    \begin{align*}
        \eta_k\colon C_*(F_1)\otimes C_*(\Omega)^{\otimes k-1}\to C_{k-1+m+*}(F_2).
    \end{align*}
\end{Def}

The following proposition makes precise that a coherent chain homotopy is a witness that a map $\widetilde{\eta}\colon C_*(M,C_*(F_1))\to C_*(M,C_{m+*}(F_2))$ models the map $\eta\colon C_*(E_1)\to C_*(E_2)$. 

\begin{Prop}\label{prop:coh chain hpty witness}
    If there exists a coherent chain homotopy $\eta_k$ for $\eta\colon C_*(E_1)\to C_{m+*}(E_2)$, the restrictions 
    \begin{align*}
        \eta_k\colon C_*(F_1)\otimes C_*(\Omega M)^{\otimes (k-1)}\to C_{k-1+m+*}(F_2)
    \end{align*}
    define a morphism of $\mathcal{A}_\infty$-modules and the following diagram commutes
    \begin{center}
        \begin{tikzcd}
            H_*(M,C_*(F_1))\ar[r, "\widetilde{\eta}"] \ar[d, "\Psi_1"] & H_*(M,C_{m+*}(F_2) \ar[d, "\Psi_2"]\\
            H_*(E_1) \ar[r, "\eta"]& H_*(E_2).
        \end{tikzcd}
    \end{center}
\end{Prop}
\begin{proof}
    Because the equation for morphisms of $\mathcal{A}_\infty$-modules in Definition \ref{def:morph of Ainfty modules} is simply a restriction of the $\mathcal{A}_\infty$-equation in Definition \ref{def:coherent chain hpty}, it is immediate that the restrictions
    \begin{align*}
        \eta_k\colon C_*(F_1)\otimes C_*(\Omega M)^{k-1}\to C_{k-1+m+*}(F_2)
    \end{align*}
    define a morphism of $\mathcal{A}_\infty$-modules.\\
    
    We recall that the quasi-isomorphisms $\Psi_i\colon C_*(M,C_*(F_i))\to C_*(E_i)$ are given by compositions
    \begin{align*}
        C_*(M,C_*(F_i))\to C_*(M/T, C_*(F_i'))\overset{\Psi'_i}{\to} C_*(E_i')\overset{\theta_*}{\to} C_*(E_i),
    \end{align*}
    where $E_i'$ is the pullback of $E_i$ along the homotopy equivalence $\theta\colon M/T\to M$ as in Subsection \ref{subsec:Morse homo with DG}. We may identify $F_i\cong F_i'$ and then the first map becomes an isomorphism. We may also assume that $\theta$ is a homeomorphism. Then $\theta_*\colon C_*(E_i')\to C_*(E_i)$ is thus an isomorphism.
    
    We now consider the following diagram:
    \begin{center}
        \begin{tikzcd}
            C_*(M,C_*(F_1))\ar[r] \ar[d,"\widetilde{\eta}"] & C_*(M/T,C_*(F_1'))\ar[r, "{\Psi_1'}"]\ar[d,"\theta^*\widetilde{\eta}"] & C_*(E_1')  \ar[r, "\theta_*"]\ar[d,"\theta^*\eta"]& C_*(E_1)\ar[d,"\eta"]\\
            C_*(M,C_{m+*}(F_2))\ar[r] & C_*(M/T,C_{m+*}(F_2'))\ar[r, "{\Psi_2'}"] & C_{m+*}(E_2')  \ar[r, "\theta_*"]& C_{m+*}(E_2).
        \end{tikzcd}
    \end{center}
    The first and third square commute as the horizontal maps are given by isomorphisms and the vertical maps commute with these identifications.
    
    It thus remains to give a chain homotopy making the middle square commute. By the above identifications $C_*(F_1)\cong C_*(F_1')$ and $C_*(E_i)\cong C_*(E_i')$, we may interpret $\eta_k$ as maps
    \begin{align*}
        C_*(F_1')\otimes C_*(\Omega(M/T))^{\otimes (k-2)}\otimes C_*(\mathcal{P}_{x_0\to M/T}(M/T))\to C_{k-1+m+*}(E_2')
    \end{align*}
    and $\eta_1$ as a map $C_*(E_1')\to C_*(E_2')$.
    
    We recall some notation from Subsection \ref{subsec:Morse homo with DG}. We have for $x,y\in \Crit(f)$ the chains $m_x'\in C_{*}(\mathcal{P}_{x_0\to M/T} (M/T))$ and $m_{x,y}'\in C_*(\Omega(M/T))$. Moreover, $\boldsymbol{m}'$ denotes the degree $-1$ map
    \begin{align*}
        \Z\Crit(f)&\to C_*(\Omega (M/T))\otimes \Z\Crit(f),\\
        x&\mapsto \sum_y m_{x,y}'\otimes y
    \end{align*}
    which we extend to a map 
    \begin{align*}
        \widetilde{\boldsymbol{m}}\colon C_*(F_1')\otimes  T(C_*(\Omega (M/T)))\otimes \Z\Crit(f)\to  C_*(F_1')\otimes T(C_*(\Omega (M/T)))\otimes \Z\Crit(f)
    \end{align*}
    by $\id\otimes \id\otimes \boldsymbol{m}'$.
    
    We also denote 
    \begin{align*}
        \widetilde{\boldsymbol{l}}\colon C_*(F_1')\otimes T(C_*(\Omega(M/T)))\otimes \Z\Crit(f)&\to C_*(F_1')\otimes T(C_*(\Omega(M/T)))\otimes C_*(\mathcal{P}_{x_0\to M/T}(M/T)),\\
        \alpha\otimes \lambda_1\otimes \dots \otimes \lambda_k\otimes x&\mapsto \alpha\otimes \lambda_1\otimes \dots \otimes \lambda_k\otimes m_x'.
    \end{align*}
    Then we define
    \begin{align*}
        v:=\sum_{N\geq 1} (-1)^{N+1}\eta_{N+1}\circ \widetilde{\boldsymbol{l}}\circ \widetilde{\boldsymbol{m}}^{N-1}\colon C_*(M,C_*(F_1))\to C_{m+1+*}(E_2).
    \end{align*}
    This is well-defined because $\widetilde{\boldsymbol{m}}^{n+1}=0$. The proposition now follows from the following lemma:
    \begin{Lem}\label{lem:v is hpty}
        It holds
        \begin{align*}
            v d+(-1)^{m}d v=(\theta^*\eta) \circ \Psi_1'-\Psi_2'\circ (\theta^*\widetilde{\eta}).
        \end{align*}
    \end{Lem}
    \begin{proof}[Proof of Lemma \ref{lem:v is hpty}]
        We compute by \eqref{eq:Ainfty coh chain hpty}
        \begin{align*}
            (-1)^{m}d v=&\ \sum_{N\geq 1} (-1)^{N+1+m}d\eta_{N+1} \widetilde{\boldsymbol{l}}\widetilde{\boldsymbol{m}}^{N-1}\\
            =& -\sum_{N\geq 1}\eta_{N+1}d\widetilde{\boldsymbol{l}}\widetilde{\boldsymbol{m}}^{N-1}\\
            &+\sum_{N\geq 1}(-1)^{N+1}\eta_{N}(\Phi_1\otimes \id^{\otimes N-1})\widetilde{\boldsymbol{l}}\widetilde{\boldsymbol{m}}^{N-1}\\
            &-\sum_{N\geq 1} \Phi_2(\eta_N\otimes \id)\widetilde{\boldsymbol{l}}\widetilde{\boldsymbol{m}}^{N-1}\\
            &+\sum_{N\geq 1}\sum_{r=1}^{N-1}(-1)^{N+1+r}\eta_N(\id^{\otimes r}\otimes \mu\otimes \id^{\otimes N-1-r})\widetilde{\boldsymbol{l}}\widetilde{\boldsymbol{m}}^{N-1}.
        \end{align*}
        We turn our attention to the first term:
        \begin{align*}
            \eta_{N+1}d\widetilde{\boldsymbol{l}}\widetilde{\boldsymbol{m}}^{N-1}=\sum_{r=0}^{N}\eta_{N+1}(\id^{\otimes r}\otimes d\otimes \id^{\otimes N-r})\widetilde{\boldsymbol{l}}\widetilde{\boldsymbol{m}}^{N-1}.
        \end{align*}
        Because $\widetilde{\boldsymbol{l}}$ and $\widetilde{\boldsymbol{m}}$ are the identity on all components except the last, we have
        \begin{align*}
            \eta_{N+1}(d\otimes \id^{\otimes N})\widetilde{\boldsymbol{l}}\widetilde{\boldsymbol{m}}^{N-1}=(-1)^{N-1}\eta_{N+1}\widetilde{\boldsymbol{l}}\widetilde{\boldsymbol{m}}^{N-1}(d\otimes \id)
        \end{align*}
        and, for $0<r<N$,
        \begin{align*}
            \eta_{N+1}(\id^{\otimes r}\otimes d\otimes \id^{\otimes N-r})\widetilde{\boldsymbol{l}}\widetilde{\boldsymbol{m}}^{N-1}&=(-1)^{N-1-r}\eta_{N+1}\widetilde{\boldsymbol{l}}\widetilde{\boldsymbol{m}}^{N-1-r}(\id^{\otimes r}\otimes d\otimes \id)\widetilde{\boldsymbol{m}}^{r}\\
            &= (-1)^{N+r}\eta_{N+1}\widetilde{\boldsymbol{l}}\widetilde{\boldsymbol{m}}^{N-1-r}(\id^{\otimes r}\otimes\mu\otimes \id)\widetilde{\boldsymbol{m}}^{r+1}\\
            &=(-1)^{N+r}\eta_{N+1}(\id^{\otimes r}\otimes\mu\otimes \id^{\otimes N-r})\widetilde{\boldsymbol{l}}\widetilde{\boldsymbol{m}}^{N},
        \end{align*}
        where the second equality follows from the compatibility of the $m_{x,y}'$:
        \begin{align*}
            d (m'_{x,y})=\sum_{z}(-1)^{|x|-|z|}\mu(m_{x,z}'\otimes m_{z,y}')
        \end{align*}
        (see \cite[Lemma 4.4]{riegel2024chain} for a detailed treatment of the signs). Similarly, we have $d(m_x')=\sum_{y}\mu(m_{x,y}'\otimes m_y')$ and thus
        \begin{align*}
            \eta_{N+1}(\id^{\otimes N}\otimes d)\widetilde{\boldsymbol{l}}\widetilde{\boldsymbol{m}}^{N-1}=\eta_{N+1}(\id^{\otimes N}\otimes \mu)\widetilde{\boldsymbol{l}}\widetilde{\boldsymbol{m}}^{N}.
        \end{align*}
        We therefore have:
        \begin{align*}
                \sum_{N\geq 1}\eta_{N+1}d\widetilde{\boldsymbol{l}}\widetilde{\boldsymbol{m}}^{N-1}=&\ \sum_{N\geq 1}\eta_{N+1}\widetilde{\boldsymbol{l}}\widetilde{\boldsymbol{m}}^{N-1}(d\otimes \id)\\
                &+\sum_{N\geq 1}\sum_{r=1}^{N-1}(-1)^{N+r}\eta_{N+1}(\id^{\otimes r}\otimes\mu\otimes \id^{\otimes N-r})\widetilde{\boldsymbol{l}}\widetilde{\boldsymbol{m}}^{N}\\
                &+\sum_{N\geq 1}\eta_{N+1}(\id^{\otimes N}\otimes \mu)\widetilde{\boldsymbol{l}}\widetilde{\boldsymbol{m}}^{N}\\
                =&\ v(d\otimes \id)+\sum_{N\geq 1}\sum_{r=1}^{N}(-1)^{N+r}\eta_{N+1}(\id^{\otimes r}\otimes\mu\otimes \id^{\otimes N-r})\widetilde{\boldsymbol{l}}\widetilde{\boldsymbol{m}}^{N}.
        \end{align*}
        We therefore compute
        \begin{align*}
            (-1)^{m}d v
                =& -v(d\otimes \id)\\
                &-\sum_{N\geq 1}\sum_{r=1}^{N}(-1)^{N+r}\eta_{N+1}(\id^{\otimes r}\otimes\mu\otimes \id^{\otimes N-r})\widetilde{\boldsymbol{l}}\widetilde{\boldsymbol{m}}^{N}\\
                &+\sum_{N\geq 1}(-1)^{N+1}\eta_{N}(\Phi_1\otimes \id^{\otimes N-1})\widetilde{\boldsymbol{l}}\widetilde{\boldsymbol{m}}^{N-1}\\
                &-\sum_{N\geq 1} \Phi_2(\eta_N\otimes \id)\widetilde{\boldsymbol{l}}\widetilde{\boldsymbol{m}}^{N-1}\\
                &+\sum_{N\geq 1}\sum_{r=1}^{N-1}(-1)^{N+1+r}\eta_N(\id^{\otimes r}\otimes \mu\otimes \id^{\otimes N-1-r})\widetilde{\boldsymbol{l}}\widetilde{\boldsymbol{m}}^{N-1}.
        \end{align*}
        The second and last term cancel. Moreover by the fact that $\widetilde{\boldsymbol{l}}$ and $\widetilde{\boldsymbol{m}}$ are the identity on all components except the last, we find that:
        \begin{align*}
            \sum_{N\geq 1}(-1)^{N+1}\eta_{N}(\Phi_1\otimes \id^{\otimes N-1})\widetilde{\boldsymbol{l}}\widetilde{\boldsymbol{m}}^{N-1}&=\eta_1\Phi_1\widetilde{\boldsymbol{l}}+\sum_{N\geq 2}(-1)^{N+1}\eta_{N}\widetilde{\boldsymbol{l}}\widetilde{\boldsymbol{m}}^{N-2}(\Phi_1\otimes \id)\widetilde{\boldsymbol{m}}\\
            &=(\theta^*\eta_1)\Psi_1'-v(\Phi_1\otimes \id)\widetilde{\boldsymbol{m}}
        \end{align*}
        and
        \begin{align*}
            \sum_{N\geq 1} \Phi_2(\eta_N\otimes \id)\widetilde{\boldsymbol{l}}\widetilde{\boldsymbol{m}}^{N-1}&=\sum_{N\geq 1}\Phi_2\widetilde{\boldsymbol{l}}(\eta_N\otimes \id)\widetilde{\boldsymbol{m}}^{N-1}\\
            &=\sum_{N\geq 1}\Psi_2(\eta_N\otimes \id)\widetilde{\boldsymbol{m}}^{N-1}\\
            &=\Psi_2' (\theta^*\widetilde{\eta}).
        \end{align*}
        We thus have
        \begin{align*}
             (-1)^{m}d v
                =& -v(d_{F_1}\otimes \id)+(\theta^*\eta_1)\Psi_1'-v(\Phi_1\otimes \id)\widetilde{\boldsymbol{m}}-\Psi_2 (\theta^*\widetilde{\eta})\\
                =&\ -vd +(\theta^*\eta_1)\Psi_1'-\Psi_2' (\theta^*\widetilde{\eta}).
        \end{align*}
        This proves Lemma \ref{lem:v is hpty} and thus Proposition \ref{prop:coh chain hpty witness}
    \end{proof}
\end{proof}

\begin{Cor}\label{cor:coh chain hpty relative}
    For $i=1,2$, let $\pi_i\colon E_i'\subseteq E_i\to M$ be fibrations with fibres $F_i$ and $F_i'$. Let $\Phi_i$ be transitive lifting functions on $E_i$ that restrict to transitive lifting functions on $E_i'$. Let $\eta\colon C_*(E_1)\to C_{m+*}(E_2)$ be a map of degree $m\in \Z$ that restricts to a map $C_*(E_1')\to C_{m+*}(E_2')$.
    
    If there exists a coherent chain homotopy $\boldsymbol{\eta}\colon C_*(E_1)\to C_{m+*}(E_2)$ for $\eta$ that restricts to a coherent chain homotopy $\boldsymbol{\eta}\colon C_*(E_1')\to C_{m+*}(E_2')$, then there exists a $\boldsymbol{\eta}\colon C_*(E_1,E_1')\to C_{m+*}(E_2,E_2')$ a relative coherent chain homotopy and the following square commutes
    \begin{center}
        \begin{tikzcd}
            H_*(M,C_*(F_1,F_1'))\ar[d, "\Psi_1"] \ar[r, "\widetilde{\eta}"] & H_*(M,C_{m+*}(F_2,F_2'))\ar[d, "\Psi_2"]\\
            H_*(E_1,E_1')\ar[r, "\eta"] & H_{m+*}(E_2, E_2').
        \end{tikzcd}
    \end{center}
\end{Cor}
\begin{proof}
    This follows from Proposition \ref{prop:coh chain hpty witness} and the five lemma.
\end{proof}

\subsection{Inverting Coherent Chain Homotopies}\label{subsec:inv cch}
In this subsection, we construct inverses to coherent chain homotopies when $\eta_1$ is an isomorphism or a quasi-isomorphism. This follows an analogous strategy to Appendix \ref{app:inv qi of Ainfty modules}.\\

To achieve this, we replace the topological context of $\Omega M$ and $\mathcal{P}_{x_0\to M}M$ with an algebra $A$ and a left module:

\begin{Def}
    Let $(A_*,d_A,\mu)$ be a strict differential graded algebra and $(P_*,d_P,\mu)$ a strict left $A_*$-module and $A_*\subseteq P_*$ as $A_*$-modules. A \textit{path module over the pair $(A_*,P_*)$} is given by
    \begin{itemize}
        \item a complex $(\mathcal{E}_*,d_\mathcal{E})$;
        \item a subcomplex $\mathcal{F}_*\subseteq \mathcal{E}_*$;
        \item a collection of maps
        \begin{align*}
            m_k^\mathcal{E}\colon \mathcal{F}_*\otimes A_*^{\otimes k-2}\otimes \mathcal{P}_*\to \mathcal{E}_{k-2+*},
        \end{align*}
        for $k\geq 2$;
    \end{itemize}
    such that $m_2^\mathcal{E}$ is a chain map and for $N\geq 3$, the following $\mathcal{A}_\infty$-equation is satisfied:
    \begin{align*}
        (-1)^{N+1}m_N^\mathcal{E}d+d m_N^\mathcal{E}=&\sum_{r=0}^{N-3}(-1)^r m^\mathcal{E}_{N-1}(\id^{\otimes r+1}\otimes \mu\otimes \id^{\otimes N-r-3})\\
        &-\sum_{s=2}^{N-1}(-1)^{s(N-s)}m_{N-s+1}^\mathcal{E}(m_s^\mathcal{F}\otimes \id^{\otimes N-s})
    \end{align*}
    as maps $\mathcal{F}_*\otimes A_*^{\otimes N-2}\otimes P_*\to \mathcal{E}_{N-3+*}$ and $m_k^\mathcal{E}$ restricts to maps
    \begin{align*}
        \mathcal{F}_*\otimes A_*^{\otimes k-1}\to \mathcal{F}_{k-2+*}
    \end{align*}
    A path module over $(A_*,P_*)$ is \textit{strict} if $m_k^\mathcal{E}=0$ for $k\geq 3$.
\end{Def}

\begin{Rem}
    In particular, $m_k^\mathcal{E}$ makes $\mathcal{F}_*$ into a right $\mathcal{A}_\infty$-module over $A_*$.
\end{Rem}

\begin{Def}
    Let $(\mathcal{E}^1_*,\mathcal{F}^1_*,m_k^{\mathcal{E}^1})$ and $(\mathcal{E}^2_*,\mathcal{F}^2_*,m_k^{\mathcal{E}^2})$ be path modules over $(A_*,P_*)$. A \textit{morphism of path modules $\boldsymbol{\eta}\colon \mathcal{E}^1_*\to \mathcal{E}_{m+*}^2$ of degree $m\in \Z$} is given by a chain map $\eta_1\colon \mathcal{E}^1_*\to \mathcal{E}^2_{m+*}$ and
    \begin{align*}
        \eta_k\colon \mathcal{F}^1_*\otimes A_*^{\otimes (k-2)}\otimes P_*\to \mathcal{E}^2_{k-1+m+*},
    \end{align*}
    for $k\geq 2$, that satisfy the $\mathcal{A}_\infty$-equations for $N\geq 1$:
    \begin{align*}
        \eta_{N+1}\circ d +(-1)^{N+1+m}d\circ\eta_{N+1}=& \sum_{s=1}^N (-1)^{s(N-s)}\eta_{N-s+1}(m_{s+1}^{\mathcal{E}_1}\otimes \id^{N-s})\\
        &+(-1)^N\sum_{k=1}^{N}(-1)^{k(N-k)}m_{k+1}^{\mathcal{E}^2}(\eta_{N-k}\otimes \id^{\otimes k})\\
        &+\sum_{r=1}^{N-1}(-1)^{N+1+r}\eta_N(\id^{\otimes r}\otimes \mu\otimes \id^{\otimes N-1-r})
    \end{align*}
    and restrict to maps 
    \begin{align*}
        \mathcal{F}_*^1\otimes A_*^{k-1}\to \mathcal{F}_{m+k-1+*}^2.
    \end{align*}
    The \textit{identity morphism} is given by $\eta_1=\id$ and $\eta_k=0$ for $k\geq 2$.
\end{Def}

\begin{Rem}
    On $\mathcal{F}_*^1$, the map $\boldsymbol{\eta}$ defines a morphism of $\mathcal{A}_\infty$-modules.
\end{Rem}

\begin{Exmp}\label{exmp:path modules}
    Let $\pi\colon E'\subseteq E\to M$ be fibrations with fibre $F'\subseteq F$ and transitive lifting function $\Phi$. The data $\mathcal{E}_*=C_*(E,E')$, $\mathcal{F}_*:=C_*(F,F')$ and $m_2^\mathcal{E}:=\Phi_*$ defines a strict path module over $(C_*(\Omega),C_*(\mathcal{P}_{x_0\to M}))$.
    
    A coherent chain homotopy defines a morphism of path modules.
\end{Exmp}

\begin{Lem}
    For $i=1,2,3$, let $(\mathcal{E}^i_*,\mathcal{F}^i_*,m_k^{\mathcal{E}^i})$ be path modules over $(A_*,P_*)$. Let $\boldsymbol{\eta}\colon \mathcal{E}^1_*\to \mathcal{E}^2_{m_1+*}$ and $\boldsymbol{\eta}'\colon \mathcal{E}^2_*\to \mathcal{E}^3_{m_2+*}$ be morphisms of path modules of degree $m_1$ and $m_2$, respectively. The maps $(\eta'\circ \eta)_1:=\eta_1'\circ \eta_1$ and 
    \begin{align*}
        (\eta'\circ\eta)_k:= \sum_{r=1}^k(-1)^{(r-1)(k-r)}\eta'_{k-r+1}\circ (\eta_r\otimes \id^{\otimes k-r}),
    \end{align*}
    for $k\geq 2$, define a morphism of path modules $\boldsymbol{\eta}'\circ \boldsymbol{\eta}\colon \mathcal{E}^1_*\to \mathcal{E}^3_{m_1+m_2+*}$ of degree $m_1+m_2$.
\end{Lem}
\begin{proof}
    This is a direct computation analogous to the proof that the above formulas defines composition for morphisms of $\mathcal{A}_\infty$-modules.
\end{proof}

\begin{Rem}\label{rem:compose cch}
    By Example \ref{exmp:path modules}, the previous lemma shows that we can compose coherent chain homotopies.
\end{Rem}

On this level of algebraic abstraction, we can now invert certain morphisms of path modules.

\begin{Prop}\label{prop:invert path mod iso}
    For $i=1,2$, let $(\mathcal{E}^i_*,\mathcal{F}^i_*,m_k^{\mathcal{E}^i})$ be path modules over $(A_*,P_*)$. Let $\boldsymbol{\eta}\colon \mathcal{E}^1_*\to \mathcal{E}^2_{m+*}$ be a morphism of path modules of degree $m$ such that $\eta_1$ is an isomorphism and $\eta_1(\mathcal{F}_*^1)=\mathcal{F}^2$. The map $\boldsymbol{\zeta}\colon \mathcal{E}^2_*\to \mathcal{E}^1_*$ given by $\zeta_1:=\eta^{-1}$
    and
    \begin{align*}
        \zeta_{N+1}:=-\eta_1^{-1}\circ \left(\sum_{r=0}^{N-1}(-1)^{r(N-r)} \eta_{N-r+1}\circ \left(\zeta_{r+1}\otimes \id^{\otimes N-r}\right)\right),
    \end{align*}
    for $k\geq 2$ is an inverse to $\boldsymbol{\eta}$.
\end{Prop}
\begin{proof}
    By our assumptions, this map is well-defined and restricts to a map
    \begin{align*}
        \mathcal{F}_*^2\otimes A_*^{k-1}\to \mathcal{F}_{-m+k-1+*}^1.
    \end{align*}
    It is a straight-forward computation analogous to Proposition \ref{prop:inverting infty-iso} to check that $\zeta_N$ satisfy the relevant $\mathcal{A}_\infty$-equations and that $\boldsymbol{\eta}\circ\boldsymbol{\zeta}=\id$ and $\boldsymbol{\zeta}\circ\boldsymbol{\eta}=\id$. 
\end{proof}

The next step in finding an inverse to morphisms of path modules such that $\eta_1$ is a proposition akin to the homotopy transfer theorem.

\begin{Prop}\label{prop:htt for path mod}
    For $j=1,2$, let $\mathcal{F}_*^j\subseteq \mathcal{E}_*^j$ be chain complexes. Let $(i,p,h)$ exhibit $\mathcal{E}^1_*$ as a homotopy retract of $\mathcal{E}^2_*$ (see Definition \ref{def:hpty retract}) such that their restrictions to $\mathcal{F}^i_*$ exhibit $\mathcal{F}^1_*$ as a homotopy retract of $\mathcal{F}^2_*$. If $\mathcal{E}_1$ is a strict path module over $(A_*,P_*)$, then there exists a path module structure ${m}_k^{\mathcal{E}^2}$ on $\mathcal{E}^2_*$ such that $i$ and $p$ can be extended to morphisms of path modules.
\end{Prop}
\begin{proof}
    The path module structure and the higher maps are given analogous to Proposition \ref{prop:htt for Ainfty modules}. The proof is also an analogous computation.
\end{proof}

\begin{Prop}\label{prop:invert path mod qi}
    For $i=1,2$, let $(\mathcal{E}^i_*,\mathcal{F}^i_*,m_k^{\mathcal{E}^i})$ be strict path modules over $(A_*,P_*)$. Let $\boldsymbol{\eta}\colon \mathcal{E}^1_*\to \mathcal{E}^2_{*}$ be a morphism of path modules such that $\eta_1$ is a quasi-isomorphism. If $(i^{\mathcal{E}^i},p^{\mathcal{E}^i},h^{\mathcal{E}^i})$ exhibits $(\mathcal{E}_*^i,d_{\mathcal{E}^i})$ and $(\mathcal{F}_*^i,d_{\mathcal{E}^i})$ as a homotopy retract of $(H_*(\mathcal{E}^i),0)$ and $(H_*(\mathcal{F}^i),0)$, respectively, then there exists a morphism of path modules $\boldsymbol{\zeta}\colon \mathcal{E}^2_*\to \mathcal{E}^1_*$ such that $H_*(\zeta_1)=H_*(\eta_1)^{-1}$.
\end{Prop}
\begin{proof}
    By Proposition \ref{prop:htt for path mod}, we can define the composition $\boldsymbol{\varepsilon}$ as 
    \begin{center}
        \begin{tikzcd}
            H_*(\mathcal{E}^1)\ar[r, "\boldsymbol{i^1}"] & \mathcal{E}^1_* \ar[r, "\boldsymbol{\eta}"] & \mathcal{E}^2_{*} \ar[r," \boldsymbol{p}^2"] & H_{*}(\mathcal{E}^2).
        \end{tikzcd}
    \end{center}
    This is such that $\varepsilon_1=H_*(\eta_1)$ is an isomorphism. By Proposition \ref{prop:invert path mod iso}, we can invert $\boldsymbol{\varepsilon}$ and can thus define $\boldsymbol{\zeta}$ as the composition
    \begin{center}
        \begin{tikzcd}
            \mathcal{E}^2_{*} \ar[r," \boldsymbol{p}^2"] & H_{*}(\mathcal{E}^2) \ar[r, "\boldsymbol{\varepsilon}^{-1}"] &H_*(\mathcal{E}^1)\ar[r, "\boldsymbol{i^1}"] & \mathcal{E}^1_*.
        \end{tikzcd}
    \end{center}
\end{proof}

\begin{Cor}\label{cor:invert quasi iso cch}
    If there exists a coherent chain homotopy $\boldsymbol{\eta}\colon C_*(E_1,E_1')\to C_{*}(E_2,E_2')$ with real coefficients for a quasi-isomorphism $\eta$, then there exists coherent chain homotopy $\boldsymbol{\zeta}\colon C_*(E_2,E_2')\to C_{*}(E_1,E_1')$ for a map $\zeta$ where $H_*(\zeta)=H_*(\eta)^{-1}$.
\end{Cor}
\begin{proof}
    By \cite[Lemma 9.4.4]{loday2012algebraic}, we can exhibit $H_*(E_i,E_i')$ as a homotopy retract of $C_*(E_i,E_i')$ over field coefficients. Then Example \ref{exmp:path modules} and Proposition \ref{prop:invert path mod qi} conclude the proof.
\end{proof}

\section{Orthogonal Half-Constant Embeddings}\label{sec:ortho hc emb}
In this Section, we define a certain class of embeddings of $U_{x_0}\hookrightarrow U_\Delta$ that play the analogous role to $\SOEmb(U_{x_0},U_{x_0})$. In Section \ref{sec:ainfty struc}, we constructed higher homotopies to exhibit the cap product with $e_\Omega^*\tau_{x_0}$ as a morphism of $\mathcal{A}_\infty$-morphisms. The two key ingredients were a Thom class $\tau_{x_0}\in C^n(U_{x_0},V_{x_0}^c)$ which is invariant under $\SOEmb_*(U_{x_0},U_{x_0})$ and exhibiting $\SOEmb_*(U_{x_0},U_{x_0})$ as a deformation retract of $\Emb_*^+((U_{x_0},V_{x_0}^c),(U_{x_0},V_{x_0}^c))$.

In this section, we construct a deformation retract of certain embeddings $U_{x_0}\hookrightarrow U_\Delta$. Together with the invariant Thom class we construct in Section \ref{sec:inv Thom class}, this lets us define a coherent chain homotopy for capping with the Thom class.\\

We recall the notation
\begin{align*}
    \Theta_0\colon \Omega&\to \Emb_*^+((U_{x_0},V_{x_0}^c),(U_{x_0},V_{x_0}^c)),\\
    \lambda&\mapsto \Theta(\lambda)|_{U_{x_0}}.
\end{align*}
The role that this plays in Section \ref{sec:ainfty struc} is that $e_\Omega$ relates the $\Omega$-action on $(e_\Omega^{-1}(U_{x_0}),e_\Omega^{-1}(V_{x_0}^c))$ to the $\Emb_*^+((U_{x_0},V_{x_0}^c),(U_{x_0},V_{x_0}^c))$-action on $(U_{x_0},V_{x_0}^c)$ given by $\Theta_0$, i.e., the following diagram commutes for $\lambda\in \Omega$:
\begin{center}
    \begin{tikzcd}
        (e_\Omega^{-1}(U_{x_0}),e_\Omega^{-1}(V_{x_0}^c))\ar[r, "{\Theta(\lambda)_*}"] \ar[d, "e_\Omega"] & (e_\Omega^{-1}(U_{x_0}),e_\Omega^{-1}(V_{x_0}^c)) \ar[d, "e_\Omega"]\\
        (U_{x_0},V_{x_0}^c) \ar[r, "{\Theta_0(\lambda)}"] &(U_{x_0},V_{x_0}^c).
    \end{tikzcd}
\end{center}
We now define $\Theta_1$ such that the following diagram commutes
\begin{center}
    \begin{tikzcd}
        (e_\Omega^{-1}(U_{x_0}),e_\Omega^{-1}(V_{x_0}^c))\ar[r, "\Theta(\lambda)_*"] \ar[d, "e_\Omega"] & (e_\Lambda^{-1}(U_{x_0}),e_\Lambda^{-1}(V_{x_0}^c)) \ar[d, "e_\Lambda"]\\
        (U_{x_0},V_{x_0}^c) \ar[r, "\Theta_1(\lambda)"] &(U_{\Delta},V_{\Delta}^c)
    \end{tikzcd}
\end{center}
for all $\lambda\in \mathcal{P}_{x_0\to M}$.

We recall that $e_\Lambda$ sends $(t,\gamma)\in I\times \Lambda$ to $(e_0(\gamma),e_t(\gamma))$. By Lemma \ref{lem:ass trans lift fun} \eqref{item:ass trans lift fun (b)}, we have $\Theta(\lambda)(e_t(t,\gamma))=e_t(\Theta(\lambda)_*(t,\gamma))$. Moreover for $\lambda\in \mathcal{P}_{x_0\to M}$, we recall that $\Theta(\lambda)$ sends $x_0$ to $e_1(\lambda)$. We thus have
\begin{align*}
    \Theta_1(\lambda)\colon (U_{x_0},V_{x_0}^c)& \to (U_{\Delta},V_{\Delta}^c),\\
    x&\mapsto(e_1(\lambda), \Theta(\lambda)(x)).
\end{align*}
Therefore, the only maps $(U_{x_0},V_{x_0}^c)\to (U_{\Delta},V_{\Delta}^c)$ that appear as $\Theta_1(\lambda)$ are embeddings which are constant in the first variable. We call them half-constant embeddings.

\begin{Not}\label{Not:iota}
    For $x\in M$, we denote $U_x=\{y\in M\mid d_M(x,y)<\varepsilon\}$ and
    \begin{align*}
        \iota_x\colon U_x &\hookrightarrow U_\Delta,\\
        y&\mapsto(x,y).
    \end{align*}
\end{Not}

\begin{Def}
    An embedding $\varphi\colon U_{x_0}\hookrightarrow U_\Delta$ is \textit{half-constant} if there exists an $x\in M$ such that
    \begin{align*}
        \varphi(y)=(x,\psi(y))=\iota_x\circ \psi
    \end{align*}
    for some $\psi\colon U_{x_0}\hookrightarrow U_x$ continuous piecewise smooth, orientation-preserving and $\psi(x_0)=x$. The space of all half-constant embeddings is denoted by $\hcEmb(U_{x_0},U_\Delta)$ and
    \begin{align*}
        \hcEmb((U_{x_0},V_{x_0}^c),(U_\Delta,V_\Delta^c))
    \end{align*}
    denotes the subspace of all relative, half-constant embeddings that send $V_{x_0}^c$ to $V_\Delta^c$.
\end{Def}

The relative, half-constant embeddings play the role of $\Emb_*^+((U_{x_0},V_{x_0}^c),(U_{x_0},V_{x_0}^c))$ from Section \ref{sec:ainfty struc}. The role of $\SOEmb(U_{x_0},U_{x_0})$ is played by orthogonal half-constant embeddings. 

\begin{Def}
    The exponential map $\exp_x$ identifies $\{v\in T_xM\mid \|v\|<\varepsilon\}$ with $U_x$ and thus the Riemannian metric induces a scalar product on $U_x$.
    
    A half-constant embedding $\varphi=\iota_x\circ \psi$ is \textit{orthogonal} if $\psi$ commutes with the scalar product on $U_{x_0}$ and $U_x$. The space of all orthogonal half-constant embeddings is denoted by $\SOEmb(U_{x_0},U_\Delta)$.
\end{Def}

Now we can introduce a homotopy that exhibits $\SOEmb(U_{x_0},U_\Delta)$ as a deformation retract of $\hcEmb((U_{x_0},V_{x_0}^c),(U_\Delta,V_\Delta^c))$.

\begin{Lem}\label{lem:orth are homotopy retract}
    There exists a homotopy
    \begin{align*}
        H_1\colon I\times \hcEmb((U_{x_0},V_{x_0}^c),(U_\Delta,V_\Delta^c))\to \hcEmb((U_{x_0},V_{x_0}^c),(U_\Delta,V_\Delta^c))
    \end{align*}
    such that for all $\varphi,\varphi_1\in \hcEmb((U_{x_0},V_{x_0}^c),(U_\Delta,V_\Delta^c))$, $\varphi_2\in \Emb((U_{x_0},V_{x_0}^c),(U_{x_0},V_{x_0}^c))$ and $t\in I$
    \begin{enumerate}[(i)]
        \item\label{item:H1 Lambda i} $H_1(0,\varphi)=:r(\varphi)\in \SOEmb(U_{x_0},U_\Delta)$;
        \item $H_1(1,\varphi)=\varphi$;
        \item\label{item:H1 Lambda iii} $H_1(t, r(\varphi_1)\circ \varphi_2)= r(\varphi_1)\circ H_1(t,\varphi_2)$;
        \item\label{item:H1 Lambda iv} if $\varphi=\iota_x\circ \psi$ then $H_1(t,\varphi)(v)\in \iota_x(U_x)$ for all $v\in U_{x_0}$.
    \end{enumerate}
    Moreover, it restricts to the homotopy
    \begin{align*}
        H_1\colon I\times \Emb_*^+((U_{x_0},V_{x_0}^c),(U_{x_0},V_{x_0}^c))\to \Emb_*^+((U_{x_0},V_{x_0}^c),(U_{x_0},V_{x_0}^c))
    \end{align*}
    from Lemma \ref{lem:H1 for Omega}.
\end{Lem}
\begin{proof}
    We fix identifications $\alpha_x\colon U_{x_0}\cong U_x$ for all $x\in M$ that send $x_0$ to $x$ and commutes with the scalar product. These identifications are not assumed to satisfy any compatibility relation for different $x$.
    
    For $\varphi=\iota_x\circ \psi$, we define
    \begin{align}\label{eq:def of H1 Lambda}
        H_1(t,\varphi):=\iota_x\circ \alpha_x\circ H_1(t,\alpha_x^{-1}\circ \psi),
    \end{align}
    where $H_1(t,\alpha_x^{-1}\circ \psi)$ is defined by Lemma \ref{lem:H1 for Omega} because $\alpha_x^{-1}\circ \psi$ is in $\Emb_*((U_{x_0},V_{x_0}^c),(U_{x_0},V_{x_0}^c))$.
    
    The definition \eqref{eq:def of H1 Lambda} does not depend on our choice of $\alpha_x$: in fact two different identifications $\alpha_x,\alpha_x'$ differ by a unique $\varphi\in \SOEmb(U_{x_0},U_{x_0})$:
    \begin{align*}
        \alpha_x=\alpha_x'\circ \varphi.
    \end{align*}
    We thus have
    \begin{align*}
        \iota_x\circ \alpha_x\circ H_1(t,\alpha_x^{-1}\circ \psi)&=\iota_x\circ \alpha_x'\circ \varphi\circ H_1(t, \varphi^{-1}\circ(\alpha_x')^{-1}\circ \psi)\\
        &=\iota_x\circ \alpha_x'\circ \varphi\circ \varphi^{-1}\circ H_1(t,(\alpha_x')^{-1}\circ \psi)\\
        &=\iota_x\circ \alpha_x'\circ H_1(t,(\alpha_x')^{-1}\circ \psi).
    \end{align*}
    This in particular shows that our definition of $H_1$ is continuous. Moreover, Condition \eqref{item:H1 Lambda iv} is immediate and Conditions \eqref{item:H1 Lambda i}-\eqref{item:H1 Lambda iii} are inherited from the analogous conditions in Lemma \ref{lem:H1 for Omega}.
\end{proof}

\section{Invariant Thom Classes}\label{sec:inv Thom class}

In this section, we construct an explicit Thom class that is invariant under orthogonal half-constant embeddings. To construct an explicit Thom class, we review the construction of the Mathai-Quillen Thom class in Subsection \ref{subsec:Mathai Quillen Thom class}. This gives a Thom class in de Rham homology. In Subsection \ref{subsec:smooth}, we discuss how to smooth singular chains such that evaluating the Mathai-Quillen Thom class is well-defined and satisfies the desired invariants.\\ 

We denote the map
\begin{align*}
    \eta\colon (U_{x_0},V_{x_0}^c)\times \hcEmb((U_{x_0},V_{x_0}^c),(U_\Delta,V_\Delta^c))&\to (U_{\Delta},V_\Delta^c),\\
    (x,\varphi)&\mapsto \varphi(x).
\end{align*}

The goal of this section is to prove the following proposition:
\begin{Prop}\label{prop:invariant Thom class}
    There exists a Thom class (with real coefficients) $\tau_\Delta\in C^n(U_\Delta,V_\Delta^c)$ such that 
    \begin{align*}
        \tau_\Delta\circ \eta\circ \times =\tau_\Delta\circ  \mathrm{pr}_1\circ \times
    \end{align*}
    holds as maps $\bigoplus_{i+j=n} C_i(U_{x_0},V_{x_0})\otimes C_j(\SOEmb(U_{x_0},U_\Delta))\to \R$.
\end{Prop}

Restricting along $\iota_{x_0}$, Proposition \ref{prop:invariant Thom class} implies the following corollary:
\begin{Cor}\label{cor:invariant taux0}
    There exists a Thom class (with real coefficients) $\tau_{x_0}\in C^n(U_{x_0},V_{x_0}^c)$ such that 
    \begin{align*}
        \tau_{x_0}\circ \eta\circ \times =\tau_{x_0}\circ  \mathrm{pr}_1\circ \times
    \end{align*}
    holds as maps $C_*(U_{x_0},V_{x_0}^c)\otimes C_*(\SOEmb(U_{x_0},U_{x_0}))\to \R$.
\end{Cor}

\begin{Rem}
    We use de Rham homology in the construction of such a Thom class. This is the bottleneck that forces us to work with real coefficients. In Lemma \ref{lem:bary subd} and Corollary \ref{cor:invert quasi iso cch}, we need field coefficients. Otherwise, our proof strategy works for any Thom class satisfying Proposition \ref{prop:invariant Thom class} over any ground ring.

    However if such an invariant Thom class exists, the coefficient ring must contain $\Q$ as the following proposition shows. Without bigger changes to our proof strategy, one can thus only hope to improve the statement from real coefficient to rational coefficients.
\end{Rem}

\begin{figure}
    \centering
    \includegraphics[width=0.3\linewidth]{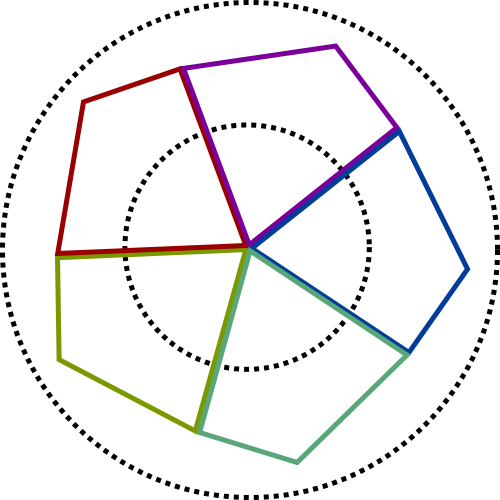}
    \caption{A cycle that cuts out a fifth of the $V_{x_0}$.}
    \label{fig:cycle forcing Q}
\end{figure}

\begin{Prop}
    Let $R$ be a ring and assume that $n>1$. If $\tau_{x_0}\in C^n(U_{x_0},V_{x_0}^c;R)$ is a Thom class that satisfies that $\tau_{x_0}\circ \eta\circ \times =\tau_{x_0}\circ  \mathrm{pr}_1\circ \times$ as maps $C_n(U_{x_0},V_{x_0}^c)\otimes C_0(\SOEmb(U_{x_0},U_{x_0}))\to R$, then $R$ contains $\Q$. 
\end{Prop}
\begin{proof}
    After intersecting $U_{x_0}$ with a plane containing the origin, we can reduce to the statement for $n=2$. 

    In dimension $2$, there exists a cycle $\sigma_0$ which cuts out $\frac{2\pi}{m}$ of the circle (see Figure \ref{fig:cycle forcing Q}). Acting on it by $\varphi:=\exp(\frac{2\pi ij}{m})\in \SO(2)$ for $j=0,\dots,m-1$ closes the circle. Let $\sigma_1$ be the chain obtained by gluing together the $m$ rotated copies of $\sigma_0$. Then there exists a $3$-chain $\widehat{\sigma}$ such that
    \begin{align*}
        d \widehat{\sigma}=\sigma_1-\sum_{i=0}^{m-1} \eta(\sigma_0\times \varphi^i).
    \end{align*}
    Moreover, $\sigma_1$ is a cycle that generates $H_2(U_{x_0},V_{x_0}^c;\Z)$.

    For such $\sigma_0,\sigma_1,\widehat{\sigma}$ and $\varphi$, we find
    \begin{align*}
        1=\tau_{x_0}(\sigma_1)=\tau_{x_0}(d \widehat{\sigma})+\sum_{i=0}^{m-1} \tau_{x_0}(\eta(\sigma_0\times \varphi^i))=0+\sum_{i=0}^{m-1} \tau_{x_0}(\mathrm{pr}_1(\sigma_0\times \varphi^i))=m\cdot\tau_{x_0}(\sigma_0).
    \end{align*}
    This shows that the element $m\in R^\times$ and thus $\Q\subseteq R$.
\end{proof}

\subsection{Mathai-Quillen Thom Class}\label{subsec:Mathai Quillen Thom class}
We recall the construction of the Mathai-Quillen Thom class. This is a de Rham cocycle $\omega_\Delta\in C^n_{\mathrm{dR}}(U_\Delta,V^c_\Delta)$. This class was originally introduced in \cite{mathai1986superconnections}. We follow the construction as in \cite[Section 1.6]{getzler1991thom}.\\

We denote $\pi\colon U_\Delta\to M$ the projection to the first coordinate on $U_\Delta\subseteq M\times M$. This is a vector bundle. We denote for $i,j\geq 0$ 
\begin{align*}
    \mathcal{A}^{i,j}:=\Gamma(U_\Delta, \Lambda^i T^*U_\Delta\otimes \Lambda^j \pi^* U_\Delta).
\end{align*}
The orientation on $M$ defines an orientation on $U_\Delta$. This thus defines a map
\begin{align*}
    B\colon \Lambda^n U_\Delta\to M\times \R
\end{align*}
called the \textit{Berezin integral}. The Berezin integral can be extended to a function
\begin{align*}
    B\colon \mathcal{A}^{i,j}&\to \Gamma(U_\Delta,\Lambda^i T^*U_\Delta)=C^i_{\mathrm{dR}}(U_\Delta),\\
    \chi\otimes \xi &\mapsto \chi \cdot B(\pi_*\xi).
\end{align*}
This map is zero if $j\neq n$.

The Riemannian metric defines a metric on the bundle $\pi\colon U_\Delta\to M$. We fix some connection $\nabla\colon \Gamma(M, \Lambda^iT^*M\otimes U_\Delta)\to \Gamma(M, \Lambda^{i+1}T^*M\otimes U_\Delta)$ that is compatible with the metric. The curvature $R:=\nabla^2$ can be interpreted as an element in
\begin{align*}
    \Gamma(M, \Lambda^2 T^* M\otimes \Lambda^2 U_\Delta).
\end{align*}
Therefore $\pi^*R$ can be interpreted as an element in $A^{2,2}$.

The canonical section
\begin{align*}
    \zeta\colon U_\Delta&\to \pi_*U_\Delta,\\
    x&\mapsto (x,x)
\end{align*}
defines an element in $\Gamma(U_\Delta, \Lambda^{0}T^*U_\Delta\otimes U_\Delta)$. We can thus interpret $(\pi^*\nabla)(\zeta)\in \Gamma(U_\Delta,  \Lambda^{1}T^*M\otimes U_\Delta)\subseteq \mathcal{A}^{1,1}$. Finally, we interpret $\|\cdot \|\colon U_\Delta\to \R$ as an element in $\mathcal{A}^{0,0}$.

We define $\omega_0:=\pi^* R+ (\pi^*\nabla)(\zeta)+\|\cdot \|\in \mathcal{A}^{2,2}\oplus \mathcal{A}^{1,1}\oplus \mathcal{A}^{0,0}$.

For $f\colon \R\to \R$, we define
\begin{align*}
    f(\omega_0)&\in \mathcal{A}^{*,*},\\
    f(\omega_0)(x)&:=\sum_{k=0}^n \frac{f^{(k)}(\|x\|)}{k!}((\pi^*\nabla)(\zeta)_{\pi(x)}+  \pi^*R_{\pi(x)})^k.
\end{align*}

We fix some $f\colon \R\to \R$ smooth such that $f|_{[\varepsilon,\infty)}\equiv 0$ and $\int_{\R^n}f^{(n)}(\tfrac{\|x\|}{2})dx=1$.

\begin{Def}
    The \textit{Mathai-Quillen Thom class} is defined as
    \begin{align*}
        \omega_\Delta:=(-1)^{\frac{n(n+1)}{2}}B(f(\omega_0))\in C^*_{\mathrm{dR}}(U_\Delta,V_\Delta^c).
    \end{align*}
\end{Def}

\begin{Prop}[{\cite[Proposition 1.3]{getzler1991thom}}]\label{prop:Mathai-Quillen Thom}
    The Mathai-Quillen Thom class $\omega_\Delta$ satisfies
    \begin{enumerate}[(i)]
        \item $\omega_\Delta\in C^n_{\mathrm{dR}}(U_\Delta,V_\Delta^c)$ is a cocycle;
        \item for any $x\in M$, we have  $\int_{B_\varepsilon(x)} \iota_x^* \omega_\Delta=1$.
    \end{enumerate} 
\end{Prop}

\begin{Rem}\label{rem:Mathai Quillen is nice}
    On a fibre $U_x$, the class $\iota_x^*\omega_\Delta$ can be described explicitly. We fix some coordinates on $B_\varepsilon(x)\cong B_\varepsilon(0)\subseteq \R^n$ such that the Riemannian metric is given by the standard scalar product on $\R^n$. Then $\iota_x^*\omega_\Delta$ is given by
    \begin{align}\label{eq:omega Delta on iota x}
        \iota_x^*\omega_\Delta(v)=f^{(n)}(\|v\|) d v_1\cdots d v_n
    \end{align}
    for $v\in \R^n$.
    
    Indeed the term $\pi^*R_{\pi(x)}$ is zero on $T^*U_x$. Moreover, we have $(\pi^*\nabla)(\zeta)=\sum_{i=1}^n v_i\otimes d v_i$. We therefore have on $U_x$ that 
    \begin{align*}
        \iota_x\omega_0(v)&= \sum_{k=0}^n \frac{f^{(k)}(\|v\|)}{k!}\left(\sum_{i=1}^nv_i\otimes d v_i\right)^k.
    \end{align*}
    The term $\left(\sum_{i=1}^mv_i\otimes d v_i\right)^k$ lives in $\mathcal{A}^{k,k}$. Thus the only term that survives after evaluating it under the Berezin integral $B$ is when $k=n$. We thus have 
    \begin{align*}
        \iota_x^*\omega_\Delta&=(-1)^{\frac{n(n+1)}{2}}\frac{f^{(n)}(\|v\|)}{n!}B\left(\left(\sum_{i=1}^mv_i\otimes d v_i\right)^n\right)\\
        &=(-1)^{\frac{n(n+1)}{2}}\frac{f^{(n)}(\|v\|)}{n!}B\left((-1)^{\frac{n(n+1)}{2}} n! v_1\wedge\dots \wedge w_n\otimes dv_1\cdots dv_n\right)\\
        &=f^{(n)}(\|v\|) d v_1\cdots d v_n.
    \end{align*}
    A consequence of this description is that the Mathai-Quillen Thom class is invariant orthogonal half-constant maps $\varphi=\iota_x\circ\psi$. Because $D_{x_0}\varphi$ preserves the Riemannian metric, it follows that 
    \begin{align}\label{eq:omega Delta is invariant}
        \iota_{x_0}^*\omega_\Delta=\psi^* \iota_x^*\omega_\Delta=\varphi^*\omega_\Delta.
    \end{align}
\end{Rem}

\subsection{Smoothing Chains}\label{subsec:smooth}
The problem with using the Mathai-Quillen Thom class is that we can only evaluate $\omega_\Delta$ on smooth chains in $C_n(U_\Delta,{V}_\Delta^c)$. To solve this issue we smooth chains in $C_n(U_\Delta,{V}_\Delta^c)$. To keep the invariance of the Mathai-Quillen Thom class, we smooth chains in a suitably equivariant way.\\

We fix the following terminology:
\begin{Def}
    A map $\sigma\colon I^k\to X$ is \textit{independent of $J\subseteq \{1,\dots,k\}$} if the composition
    \begin{align*}
        I^k\twoheadrightarrow I^{\{1,\dots,k\}\setminus J} \cong I^{\{1,\dots,k\}\setminus J}\times \{0\}^J\overset{\sigma}{\to} X
    \end{align*}
    equals $\sigma$.
\end{Def}

For $1\leq i\leq n$ and $v\in I$, we denote
\begin{align*}
    \lambda_i^v\colon I^{n-1}&\to I^n,\\
    (v_1,\dots,v_{n-1})&\mapsto (v_1,\dots, v_{i-1},v,v_i,\dots,v_n).
\end{align*}
We recall the notation 
\begin{align*}
    \mu \colon \SOEmb( U_{x_0},U_{x_0})\times \SOEmb(U_{x_0},U_\Delta)&\to \SOEmb(U_{x_0},U_\Delta),\\
    (\varphi_1,\varphi_2)&\mapsto \varphi_2\circ \varphi_1.
\end{align*}

We first smooth the chains on $\SOEmb(U_{x_0},U_\Delta)$. This makes sense because $\SOEmb(U_{x_0},U_\Delta)$ is itself a manifold with a canonical smooth structure.

\begin{Lem}\label{lem:smoothing embeddings}
    For all continuous maps $\sigma\colon I^k\to \SOEmb(U_{x_0},U_\Delta)$ there exists a homotopy 
    \begin{align*}
        H_\sigma\colon I\times I^k&\to \SOEmb(U_{x_0},U_\Delta)
    \end{align*}
    such that
    \begin{enumerate}[(i)]
        \item\label{item:smoothing emb i} $H_{\sigma}(0,v)=\sigma(v)$;
        \item\label{item:smoothing emb ii} $H_{\sigma}(1,v)=:s(\sigma)(v)$ is smooth in $v$;
        \item\label{item:smoothing emb iii} for $k=0$, it holds $H_{\sigma}(t)\equiv\sigma$;
        \item\label{item:smoothing emb iv} the homotopies commute with face maps, that is, for all $i=1,\dots, k$, it holds $H_\sigma\circ (\id\times \lambda_i^0)=H_{\sigma\circ \lambda_i^0}$ and $H_\sigma\circ (\id\times \lambda_i^1)=H_{\sigma\circ \lambda_i^1}$ as maps from $I\times I^{k-1}$;
        \item\label{item:smoothing emb independent} if $J\subseteq \{1,\dots,k\}$ is such that for $\sigma(v)$ is independent of $J$, then $H_\sigma(t,x)$ is independent of $J$ for all $t\in I$;
        \item\label{item:smoothing emb in Ux0} for $\sigma\colon I^k\to  \SOEmb(U_{x_0},U_{x_0})\subseteq \SOEmb(U_{x_0},U_{\Delta})$, it holds $H_{\sigma}(t,v)\in \SOEmb(U_{x_0},U_{x_0})$ for all $(t,v)\in I\times I^k$;
        \item\label{item:smoothing emb equiv} if $\sigma:=\mu\circ (\varphi\times \sigma_0)$ for some $\sigma_0\colon I^k\to \SOEmb(U_{x_0},U_\Delta)$ and $\varphi\in \SOEmb(U_{x_0},U_{x_0})$, we have $H_\sigma=\mu\circ (\varphi\times H_{\sigma_0})$.
    \end{enumerate}
\end{Lem}
\begin{proof}
    The proof follows the idea of smoothing singular chains on a smooth manifold in a way that commutes with faces. This is a standard result from smoothing theory (see for example \cite[Lemma 18.8]{lee2003smooth}).\\
    
    We prove this statement by induction over $k$. For $k=0$, we define $H_\sigma(t,v):=\sigma(v)$.
    
    Assume we have constructed $H_\sigma$ as above for all $\sigma\colon I^{k'}\to X$ for $k'<k$. We fix some $\sigma\colon I^k\to X$. If $\sigma$ is smooth $H_\sigma(t,v):=\sigma(v)$ does the trick.\\

    Let $J\subseteq \{1,\dots,k\}$ be the maximal subset such that $\sigma$ is independent of $J$. Up to renaming the coordinates we may assume that $J=\{1,\dots,j\}$ for some $j\leq k$.
    
    If $J\neq \varnothing$, we set $\sigma':=\sigma|_{\{0\}^J\times I^{\{j+1,\dots,k\}}}\colon I^{k-j}\to U$. By induction $H_{\sigma'}$ is defined and 
    \begin{align*}
        H_\sigma(t,v_1,\dots, v_k)&=H_{\sigma'}(t, v_{j+1},\dots, v_k)
    \end{align*}
    does the trick.\\
    
    From now on, we can thus assume that $\sigma$ is not smooth and is not independent of any $\varnothing \neq J\subseteq \{1,\dots,k\}$.
    
    We note that $\SOEmb(U_{x_0},U_{x_0})$ acts freely (via $\mu$) on all such $\sigma$. We can thus fix a system of representatives $\mathcal{X}$ among all such $\sigma$.
    
    If $\sigma\in \mathcal{X}$ and $\sigma(I^k)\subseteq \SOEmb(U_{x_0},U_{x_0})$, we do a construction of $H_\sigma$ that commutes with face maps as in \cite[Lemma 18.8]{lee2003smooth} for the manifold $\SOEmb(U_{x_0},U_{x_0})$. This ensures \eqref{item:smoothing emb in Ux0}.
    
    If $\sigma\in \mathcal{X}$ and $\sigma(I^k)\not\subseteq \SOEmb(U_{x_0},U_{x_0})$, we do a construction of $H_\sigma$ that commutes with face maps as in \cite[Lemma 18.8]{lee2003smooth} for the manifold $\SOEmb(U_{x_0},U_{\Delta})$.
    
    If $\sigma\notin \mathcal{X}$, there exists a unique $\varphi\in \SOEmb(U_{x_0},U_{x_0})$ and $\sigma'\in \mathcal{X}$ such that $\sigma=\mu(\varphi\times \sigma')$. We then define $H_{\sigma}:= \mu(\varphi\times H_{\sigma'})$. This construction thus ensures \eqref{item:smoothing emb equiv}.
\end{proof}

We now turn our attention to smoothing the chains in $C_*(U_\Delta,V_\Delta^c)$. For this, we need the following definition. 

\begin{Def}\label{def:cutting}
    A chain $\sigma\colon I^k\to U_\Delta$ has a \textit{cutting} at $0\leq i< k$ if there exist chains $\sigma_0\colon I^i\to U_{x_0}$ and $\varphi_1\colon I^{k-i}\to \SOEmb(U_{x_0},U_\Delta)$ with $\sigma=\eta\circ (\sigma_0\times \varphi_1)$.
\end{Def}

The following lemma shows that we can exhibit all cuttings of a chain at once:

\begin{Lem}\label{lem:cutting}
    Let $\sigma\colon I^k\to U_\Delta$ be a chain with cuttings at $0\leq i_1<\dots<i_l\leq k$. There exist $\sigma_0\colon I^{i_1}\to U_{x_0}$, $\varphi_j\colon I^{i_{j+1}-i_j}\to\SOEmb(U_{x_0},U_{x_0})$ for $1\leq j<l$ and $\varphi_l\colon I^{k-i_l}\to \SOEmb(U_{x_0},U_\Delta)$ with
    \begin{align*}
        \sigma=\eta\circ (\sigma_0\times \mu\circ (\varphi_1\times \mu(\sigma_2\times \dots\times\mu(\varphi_{l-1}\times \varphi_l)\dots))).
    \end{align*}
\end{Lem}
\begin{proof}
    We prove this by induction on $l$. The statement is the definition of cutting when $l=1$. We can thus assume that we have constructed $\sigma_0, \varphi_1,\dots,\varphi_{l-1}$ as in the statement for $l-1$. After composing $\varphi_1,\dots,\varphi_{l-1}$, the induction step reduces to the statement for $l=2$.
    
    We thus assume that $\sigma=\eta(\sigma_0'\times \varphi_1')=\eta(\sigma_0''\times \varphi_1'')$.
    
    We can thus compute
    \begin{align*}
        \eta\circ (\sigma_0'(v_1,\dots,v_{i_{1}})\times \varphi_{1}'(v_{i_{1}+1},\dots,v_{i_2},0,\dots,0))=\eta\circ (\sigma_0''(v_1,\dots,v_{i_2})\times \varphi_1''(0,\dots,0)).
    \end{align*}
    We note that the maps in $\SOEmb(U_{x_0},U_\Delta)$ have the property that their images agree if and only if their images intersect. This calculation thus shows $\im(\varphi_1''(0,\dots,0))=\im(\varphi_{1}'(v_{i_{1}+1},\dots,v_{i_2},0,\dots,0)$ for all $v_{i_{1}+1},\dots,v_{i_2}$. Moreover, we have by injectivity of maps in $\SOEmb(U_{x_0},U_\Delta)$:
    \begin{align*}
        \sigma_0''(v_1,\dots,v_{i_2})&=\eta(\sigma_0'(v_1,\dots,v_{i_{1}})\times (\varphi_1''(0,\dots,0)^{-1}\circ\varphi'_{1}(v_{i_{1}+1},\dots,v_{i_2},0,\dots,0))).
    \end{align*}
    We therefore find that $\sigma_0:=\sigma_0'$, $\varphi_{1}(v_{i_{1}+1},\dots, v_{i_2}):=\varphi'_{1}(v_{i_{1}+1},\dots,v_{i_2},0,\dots,0)$ and $\varphi_2:=\varphi_1''$ does the trick.
\end{proof}

We are now in a position to prove and state how we smooth chains in $C_*(U_\Delta,V_\Delta^c)$:

\begin{Lem}\label{lem:smoothing}
    For all continuous maps $\sigma\colon I^k\to U_\Delta$ there exists a homotopy 
    \begin{align*}
        H_\sigma\colon I\times I^k&\to U_\Delta
    \end{align*}
    such that
    \begin{enumerate}[(i)]
        \item\label{item:smoothing i} $H_\sigma(0,v)=\sigma(v)$;
        \item\label{item:smoothing ii} $H_\sigma(1,v)=:s(\sigma)(v)$ is smooth in $v$;
        \item\label{item:smoothing iii} if $\sigma$ is smooth $H_\sigma(t,v)=\sigma(v)$ for all $(t,v)\in I\times I^k$;
        \item\label{item:smoothing face commute} the homotopies commute with face maps, that is, for all $i=1,\dots, k$, it holds $H_\sigma\circ (\id\times \lambda_i^0)=H_{\sigma\circ \lambda_i^0}$ and $H_\sigma\circ (\id\times \lambda_i^1)=H_{\sigma\circ \lambda_i^1}$ as maps from $I\times I^{k-1}$;
        \item\label{item:smoothing independent} if $\sigma$ is independent of $J\subseteq \{1,\dots,k\}$, then $H_\sigma(t,x)$ is independent of $J$ for all $t\in I$;
        \item\label{item:smoothing Vc} if $\sigma(I^k)\subseteq V^c_\Delta$, then $H_\sigma(t,v)\in V^c_\Delta$ for all $(t,v)\in I\times I^k$;
        \item\label{item:smoothing equivariant} if $\sigma=\eta(\sigma_0\times \varphi)$ for some $\sigma_0\colon I^k\to U_{x_0}$ and $\varphi\in \SOEmb(U_{x_0},U_\Delta)$, then $H_\sigma=\eta(H_{\sigma_0}\times \varphi)$;
        \item\label{item:smoothing cutting} if $\sigma$ has a cutting at $i$, then so does $H_\sigma(t,v)$ for all $t\in I$.        
    \end{enumerate}
\end{Lem}
\begin{proof}
    The proof follows similar lines as the proof of Lemma \ref{lem:smoothing embeddings}. We do induction on $k$ and can thus assume it is proven for all $k'<k$. Moreover we may assume that $\sigma\colon I^k\to U_\Delta$ is not smooth and not independent for any $\varnothing\neq J\subseteq \{1,\dots, k\}$.\\

    We identify $\alpha \colon U_{x_0}\cong B_\varepsilon(0)$ via orthonormal coordinates.\\
    We make a case distinction, according to the criteria:
    \begin{enumerate}[(1)]
        \item whether or not $\sigma$ has a cutting;
        \item whether or not $\sigma(I^k)\subseteq U_{x_0}$;
        \item whether or not $\sigma(I^k)\subseteq V_{\Delta}^c$;
        \item\label{item:case distinc smoothing} if there exists $\sigma_0\colon I^k\to U_{x_0}$ and $\varphi\in \SOEmb(U_{x_0},U_\Delta)$ with $\sigma=\eta(\sigma_0\times \varphi)$.
    \end{enumerate}
    
    We first consider $\sigma$ with no cuttings and $\sigma(I^k)\subseteq U_{x_0}$. In particular, there always exist $\sigma_0$ and $\varphi$ as in \eqref{item:case distinc smoothing}.
    
    We define $T_\sigma\subseteq \R^n$ the linear subspace spanned by the image of $\alpha\circ \sigma$ and set $W_\sigma:=\alpha^{-1}(T_\sigma)$. We note that $\SOEmb(U_{x_0},U_{x_0})$ acts on all such $\sigma$. We fix $\mathcal{X}$ a representation system under this action. For such a $\sigma$, its stabilizer is all $\varphi\in \SOEmb(U_{x_0},U_{x_0})$ that are the identity on $W_\sigma$.\\

    We now assume that additionally $\sigma(I^k)\subseteq V_{\Delta}^c$. If $\sigma\in \mathcal{X}$, we construct $H_\sigma$ as in \cite[Lemma 18.8]{lee2003smooth} for the manifold $W_\sigma\cap V_{x_0}^c$. We note that if $\varphi$ is in the stabilizer of $\sigma$, then it also leaves $H_\sigma$ invariant.
    
    If $\sigma\notin \mathcal{X}$, there exists a $\sigma_0\in \mathcal{X}$ and $\varphi\in \SOEmb(U_{x_0},U_{x_0})$ such that $\sigma=\eta(\sigma_0\times \varphi)$. We define $H_\sigma:=\mu(H_{\sigma_0}\times \varphi)$. While $\sigma_0$ is unique, $\varphi$ is not. In fact two different $\varphi$ may exist but their difference is in the stabilizer of $\sigma_0$ and thus leaves $H_{\sigma_0}$ invariant. Therefore our definition of $H_\sigma$ is independent of our choice of $\varphi$ and \eqref{item:smoothing equivariant} is satisfied. This concludes with the first case.\\

    We now assume that $\sigma$ has no cuttings and $\sigma(I^k)\subseteq U_{x_0}$, but not $\sigma(I^k)\subseteq V_{\Delta}^c$. If $\sigma\in \mathcal{X}$, we do a construction for $H_\sigma$ as in \cite[Lemma 18.8]{lee2003smooth} for the manifold $W_\sigma$. If $\sigma\in \mathcal{X}$, we find $\sigma_0\in \mathcal{X}$ and $\varphi\in \SOEmb(U_{x_0},U_{x_0})$ such that $\eta(\sigma_0\times \varphi)$. We define $H_\sigma:=\mu(H_{\sigma_0}\times \varphi)$. As before this is independent of our choices and ensures \eqref{item:smoothing equivariant}.\\

    We now assume that $\sigma$ has no cuttings, $\sigma(I^k)\not\subseteq U_{x_0}$, $\sigma(I^k)\subseteq V_{\Delta}^c$ and there exists no $\sigma_0$ and $\varphi$ as in \eqref{item:case distinc smoothing}. We then construct $H_\sigma$ as in \cite[Lemma 18.8]{lee2003smooth} for the manifold $V_\Delta^c$.\\

    We now assume that $\sigma$ has no cuttings, $\sigma(I^k)\not\subseteq U_{x_0}$, $\sigma(I^k)\not\subseteq V_{\Delta}^c$ and there exists no $\sigma_0$ and $\varphi$ as in \eqref{item:case distinc smoothing}. We then construct $H_\sigma$ as in \cite[Lemma 18.8]{lee2003smooth} for the manifold $U_\Delta$.\\

    We now assume that $\sigma$ has no cuttings and that we can write $\sigma=\eta(\sigma_0\times \varphi)$ for some $\sigma_0\colon I^k\to  U_{x_0}$ and $\varphi\in \SOEmb(U_{x_0},U_\Delta)$. We note that if $\sigma(I^k)\subseteq V_\Delta^c$, then also $\sigma_0(I^k)\subseteq V_{x_0}^c$. We define $H_\sigma:=\mu(H_{\sigma_0}\times \varphi)$. By \eqref{item:smoothing equivariant} for chains with image in $U_{x_0}$, this definition is independent of our choice of $\sigma_0$ and $\varphi$.\\

    Finally we assume that $\sigma$ has cuttings. In particular we assume that it has precisely cuttings at $0\leq i_1<\dots<i_l<k$. We then choose $\sigma_0$ and $\varphi_j$ as in Lemma \ref{lem:cutting} for $1\leq j\leq l$.
    
    We then define 
    \begin{align*}
        H_\sigma:=\eta (H_{\sigma_0}\times \mu (H_{\varphi_1}\times \mu(H_{\varphi_2}\times \dots\times\mu(H_{\varphi_{l-1}}\times H_{\varphi_l})\dots))),
    \end{align*}
    where $H_{\varphi_j}$ is as in Lemma \ref{lem:smoothing embeddings}.
    
    This ensures \eqref{item:smoothing cutting} and the Conditions \eqref{item:smoothing i}-\eqref{item:smoothing equivariant} are inherited from $\sigma_0$ and $\varphi_i$ by Lemma \ref{lem:smoothing embeddings}. Moreover equivariance for $\sigma_0$ and $\varphi_i$ ensures that this definition is independent of choice.
\end{proof}

We can now define the following Thom class:

\begin{Def}
    The class $\tau_\Delta$ is defined as $s^*\omega_\Delta$ and $\tau_{x_0}$ is defined as $\iota_{x_0}^*\tau_\Delta$.
\end{Def}

This construction gives a Thom class as described in Proposition \ref{prop:invariant Thom class}:

\begin{proof}[Proof of Proposition \ref{prop:invariant Thom class}]
    We recall that Proposition \ref{prop:invariant Thom class} states that $\tau_{x_0}$ is a Thom class and 
    \begin{align}\label{eq:pi1 and invariant}
        \tau_\Delta\circ \eta\circ \times =\tau_{x_0}\circ  \mathrm{pr}_1\circ \times
    \end{align}
    holds as maps $C_*(U_{x_0},V_{x_0})\otimes C_*(\SOEmb(U_{x_0},U_\Delta))\to \R$.
    
    We first proof that $\tau_{x_0}$ is a Thom class. We recall that $\omega_\Delta$ is a Thom class of smooth chains (see e.g.\  \cite[Proposition 1.48]{berline2003heat}).
    
    Next, we note that the map $s_*\colon C_*(U_\Delta)\to C_*(U_\Delta)$ is a well-defined map as it sends degenerate chains to degenerate chains by Lemma \ref{lem:smoothing} \eqref{item:smoothing independent}. Because smoothing commutes with face maps as in Lemma \ref{lem:smoothing} \eqref{item:smoothing face commute}, $s_*$ is a chain map and it lands in the image of all smooth chains. Moreover, the map
    \begin{align*}
        h\colon C_*(U_\Delta)&\to C_{1+*}(U_\Delta),\\
        \sigma&\mapsto  H_\sigma(I\otimes -)
    \end{align*}
    defines a chain homotopy between $s_*$ and $\id$. On the other hand for $\sigma\in C^{\mathrm{sm}}_*(U_\Delta^c)$, we have $s_*(\sigma)=\sigma$ by Lemma \ref{lem:smoothing} \eqref{item:smoothing iii}, where $C_*^{\mathrm{sm}}$ denotes the smooth chains. We therefore have $s\colon C_*(U_\Delta)\simeq C_*^{\mathrm{sm}}(U_\Delta)$.
    
    By Lemma \ref{lem:smoothing} \eqref{item:smoothing Vc}, the map $s$ sends $C_*(V^c_\Delta)$ to $C_*^{\mathrm{sm}}(V_\Delta^c)$ and it is also a quasi-isomorphism. Therefore, by the five lemma $s\colon C_*(U_\Delta,V_\Delta^c)\to C_*^{\mathrm{sm}}(U_\Delta,V_\Delta^c)$ is also a quasi-isomorphism. This proves that $\tau_\Delta=s^*\omega_\Delta$ is a Thom class.
    
    It remains to check \eqref{eq:pi1 and invariant}. Because $\mathrm{pr}_1\circ \times$ is degenerate if $j>0$, we note that \eqref{eq:pi1 and invariant} entails two statements:
    \begin{enumerate}[(a)]
        \item\label{item:proof of inv Thom a} on $C_i(U_{x_0},V_{x_0})\otimes C_j(\SOEmb(U_{x_0},U_\Delta))$ with $j>0$ and $i+j=n$, we have $\tau_\Delta\circ \eta\circ \times =0$;
        \item\label{item:proof of inv Thom b} for $\sigma\otimes \varphi\in C_n(U_{x_0},V_{x_0})\otimes C_0(\SOEmb(U_{x_0},U_\Delta))$, we have $\tau_\Delta(\eta(\sigma\times \varphi))=\tau_{x_0}(\sigma)$.
    \end{enumerate}
    For \eqref{item:proof of inv Thom a} we fix $\sigma_0\otimes \varphi_1\in C_i(U_{x_0},V_{x_0})\otimes C_j(\SOEmb(U_{x_0},U_\Delta))$. By Definition \ref{def:cutting}, $\eta(\sigma_0\times \varphi_1)$ has a cutting at $i<n$. By Lemma \ref{lem:smoothing} \eqref{item:smoothing cutting}, $s(\eta(\sigma_0\times \varphi_1))$ has a cutting at $i<n$. That means that there exist $\sigma_0'\otimes \varphi_1'\in C_i^{\mathrm{sm}}(U_{x_0},V_{x_0})\otimes C_j^{\mathrm{sm}}(\SOEmb(U_{x_0},U_\Delta))$ with $s(\eta(\sigma_0\times \varphi_1))=\eta(\sigma_0'\times \varphi_1')$. We therefore compute:
    \begin{align*}
        \tau_\Delta\circ \eta(\sigma_0\times \varphi_1)=\omega_\Delta\circ s\circ \eta(\sigma_0\times \varphi_1)=\omega_\Delta\circ \eta(\sigma_0'\times \varphi_1')=\omega_\Delta\circ  \mathrm{pr}_1(\sigma_0'\times \varphi_1'),
    \end{align*}
    where we used that $\omega_\Delta$ is invariant in the sense of \eqref{eq:omega Delta is invariant}. But $\varphi_1'$ is not a zero chain, therefore $ \mathrm{pr}_1(\sigma_0'\times \sigma_1')$ is degenerate an the above equation evaluates to zero. This shows \eqref{item:proof of inv Thom a}.
    
    For \eqref{item:proof of inv Thom b} we note that, by Lemma \ref{lem:smoothing} \eqref{item:smoothing equivariant}, we have $s(\eta (\sigma\times \varphi))=\eta(s(\sigma)\times \varphi)$. Using again the invariance of $\omega_\Delta$ as in \eqref{eq:omega Delta is invariant}, we compute:
    \begin{align*}
        \tau_\Delta(\eta(\sigma\times \varphi))=\omega_\Delta\circ s \circ \eta(\sigma\times \varphi)=\omega_\Delta\circ \eta(s(\sigma)\times \varphi)=\varphi^*\omega_\Delta(s(\sigma))\overset{\eqref{eq:omega Delta is invariant}}{=} \iota_{x_0}^*\omega_\Delta(s(\sigma))=\tau_{x_0}(\sigma)
    \end{align*}
    which is what we wanted to proof.
\end{proof}

\section{Main Theorem}\label{sec:mainthm}
In this section, we state and prove our main theorem. We do this by constructing a coherent chain homotopy for the Goresky-Hingston coproduct in \ref{subsec:cch for GH}. We use our main theorem to compute the coproduct for manifolds $M=S^n/G$ in Subsection \ref{subsec:Sn/G}. We conclude the section by applying our main theorem to the question on the invariance of the Goresky-Hingston coproduct under maps of manifolds in Subsection \ref{subsec:inv under maps}.

\subsection{Coherent Chain Homotopies for the Goresky-Hingston Coproduct}\label{subsec:cch for GH}
In this subsection, we construct a coherent chain homotopy $\boldsymbol{\nu}\colon C_*({\Lambda}M,M)\to C_{1-n+*}(\mathcal{F}\Lambda,\mathcal{H})$ for the map $\nu\colon C_*({\Lambda}M,M)\to C_{1-n+*}(\mathcal{F}\Lambda,\mathcal{H})$ such that $\mathrm{c}\circ \nu\colon C_*({\Lambda}M,M)\to C_{1-n+*}((\Lambda M,M)^2)$ computes the Goresky-Hingston coproduct as in Subsection \ref{subsec:GH coprod}.

The process of constructing this coherent chain homotopy is similar to constructing the morphism of $\mathcal{A}_\infty$-modules $\boldsymbol{\nu}\colon C_*(\Omega,x_0)\to C_{1-n+*}((\Omega,x_0)^2)$ in \eqref{eq:Ainfty nu}. In fact the morphism of $\mathcal{A}_\infty$-modules is simply a restriction of the coherent chain homotopy.

\begin{Def}
    A coherent chain homotopy $\boldsymbol{\eta}\colon C_*(E_1)\to C_{m+*}(E_2)$ is \textit{trivial} if $\eta_k=0$ for all $k\geq 2$.
\end{Def}

We use the notation $\mathcal{P}:=\mathcal{P}_{x_0\to M}M$, $\mathcal{H}\Omega:=\partial I\times \Omega\cup I\times \{x_0\}$ and $\mathcal{H}\Lambda:=\partial I\times \Lambda\cup I\times M$. We consider the following diagram:
\begin{center}
    \begin{tikzcd}
        C_*(\Omega,x_0)  \otimes C_*(\mathcal{P}) \ar[d, "(I\times) \otimes \id"]\ar[r,"\Phi_*"] \ar[rd,phantom ,"\eqref{sq:cchI}"]& C_*(\Lambda, M)\ar[d, "I\times "]\\
        C_{1+*}(I\times \Omega, \mathcal{H}\Omega)\otimes C_*(\mathcal{P}) \ar[d] \ar[r, "\Phi_*"] \ar[rd,phantom ,"\eqref{sq:cchII}"] & C_{1+*}(I\times \Lambda,\mathcal{H}\Lambda)\ar[d]\\
        C_{1+*}(I\times \Omega, \mathcal{H}\Omega\cup e^{-1}_\Omega(V_{x_0}^c)) \otimes C_*(\mathcal{P}) \ar[d, "\simeq"]  \ar[rd,phantom ,"\eqref{sq:cchIII}"]\ar[r, "\Phi_*"] & C_{1+*}(I\times \Lambda, \mathcal{H}\Lambda \cup e^{-1}_\Lambda(V_{\Delta}^c)) \ar[d, "\simeq"]\\
        C_{1+*}(e^{-1}_\Omega(U_{x_0}),\mathcal{H}\Omega \cup e^{-1}_\Omega(V_{x_0}^c))\otimes C_*(\mathcal{P}) \ar[d, "e_\Omega^*\tau_{x_0}\cap \otimes \id "] \ar[r, "\Phi_*"]  \ar[rd,phantom ,"\eqref{sq:cchIV}"] & C_{1+*}(e^{-1}_\Lambda(U_{\Delta}),\mathcal{H}\Lambda \cup e^{-1}_\Lambda(V_{\Delta}^c)) \ar[d, "e_\Lambda^*\tau_\Lambda\cap "] \\
        C_{1-n+*}(e^{-1}_\Omega(U_{x_0}),\mathcal{H}\Omega)\otimes C_*(\mathcal{P}) \ar[d, "R_*\otimes \id"] \ar[r, "\Phi_*"]  \ar[rd, phantom, "\eqref{sq:cchV}"]&C_{1-n+*}(e^{-1}_\Lambda(U_{x_0}),\mathcal{H}\Lambda)\ar[d, "R_*"]\\
        C_{1-n+*}(\mathcal{F}\Omega,\mathcal{H}\Omega) \otimes C_*(\mathcal{P})\ar[d, "\mathrm{c}\otimes \id"]\ar[r, "\Phi_*"]  \ar[rd, phantom, "\eqref{sq:cchVI}"] & C_{1-n+*}(\mathcal{F}\Lambda,\mathcal{H}\Lambda) \ar[d, "\mathrm{c}"]\\
        C_{1-n+*}((\Omega,x_0)^2)\otimes C_*(\mathcal{P}) \ar[r, "\Phi_*"] & C_{1-n+*}((\Lambda,M)^2).
    \end{tikzcd}
\end{center}

\begin{enumerate}[(1)]
    \item\label{sq:cchI} This square commutes strictly thus the trivial coherent chain homotopy does the trick.
    \item\label{sq:cchII} This square commutes strictly thus the trivial coherent chain homotopy does the trick.
    \item\label{sq:cchIII} A coherent chain homotopy exists for this square because, by Corollary \ref{cor:invert quasi iso cch}, we can invert the trivial coherent chain homotopy for the inclusion
    \begin{align*}
        C_{*}(e^{-1}_\Lambda(U_{\Delta}),\mathcal{H} \cup e^{-1}_\Lambda(V_{\Delta}^c))\overset{\simeq}{\hookrightarrow}  C_{*}(I\times \Lambda, \mathcal{H} \cup e^{-1}_\Lambda(V_{\Delta}^c)).
    \end{align*}
    \item\label{sq:cchIV} We construct a coherent chain homotopy for this square in \ref{subsubsec:cap}.
    \item\label{sq:cchV} This square comes from a morphism of fibrations. Riegel shows that in this case a coherent homotopy exists (\cite[Lemma 5.9]{riegel2024chain}). Moreover by \cite[Lemma 5.12]{riegel2024chain}, a coherent homotopy defines a coherent chain homotopy.
    \item\label{sq:cchVI} This square commutes strictly thus the trivial coherent chain homotopy does the trick.
\end{enumerate}

\subsubsection{Capping with the Thom Class}\label{subsubsec:cap}
We define 
\begin{align*}
    q_1\colon C_*(e^{-1}_\Lambda(U_\Delta),\mathcal{H}\cup e^{-1}_\Lambda(V_\Delta))&\to C_{-n+*}(e_\Lambda^{-1}(U_\Delta,\mathcal{H})),\\
    \alpha&\mapsto e_\Lambda^*\tau_\Lambda\cap \alpha.
\end{align*}
We recall some notation: 
\begin{align*}
    \eta\colon U_{x_0}\times \Emb(U_{x_0},U_\Delta)&\to U_\Delta,\\
    (x,\varphi)&\mapsto\varphi(x)
\end{align*}
and
\begin{align*}
    \Theta_1\colon \mathcal{P}&\to \hcEmb((U_{x_0},V_{x_0}^c),(U_\Delta,V_{\Delta}^c))\subseteq \Emb(U_{x_0},U_\Delta),\\
    \lambda&\mapsto \Theta(\lambda)|_{U_{x_0}}.
\end{align*}
A similar computation to Lemma \ref{lem:rewritten q1 Phi} \eqref{item:rewritten a} shows for $\alpha\otimes \lambda\in C_*(e^{-1}_\Omega(U_{x_0}),\mathcal{H}\cup e^{-1}_\Omega(V_{x_0}))\otimes C_*(\mathcal{P})$
\begin{align*}
    q_1\circ \Phi(\alpha\otimes \lambda)=((\tau_\Lambda\circ \eta\circ (e_\Omega\times \Theta_1))\otimes \Phi)\circ \Delta_*(\alpha\times \lambda)
\end{align*}
and a calculation analogous to Lemma \ref{lem:rewritten q1 Phi} \eqref{item:rewritten b} shows
\begin{align*}
    \Phi\circ (q_1\otimes \id)(\alpha\otimes \lambda)=((\tau_\Lambda\circ \eta\circ (e_\Lambda\times (r\circ \Theta_1)))\otimes \Phi)\circ \Delta_*(\alpha\times \lambda).
\end{align*}

We recall that we constructed a homotopy $H_1$ between $r$ and $\id$ in Lemma \ref{lem:orth are homotopy retract}. We define inductively for $k\geq 2$
\begin{align*}
    H_k\colon I^k\times \hcEmb((U_{x_0},V_{x_0}^c),(U_{x_0},V_{x_0}^c))^{\times k-1}&\times \hcEmb((U_{x_0},V_{x_0}^c),(U_{\Delta},V_{\Delta}^c))\\
    &\to \hcEmb((U_{x_0},V_{x_0}^c),(U_{\Delta},V_{\Delta}^c)),\\
    (t_1,\dots,t_k,\varphi_1,\dots,\varphi_k)&\mapsto H_1(t_1,H_{k-1}(t_2,\dots,t_k,\varphi_2,\dots,\varphi_k)\circ \varphi_1)
\end{align*}
and for $k\geq 1$
\begin{align*}
    h_k\colon C_*(\hcEmb(U_{x_0},U_{x_0})^{k-1}\times \hcEmb(U_{x_0},U_\Delta))&\to C_{k+*}(\hcEmb(U_{x_0},U_\Delta),\\
    \varphi_1\times \dots\times \varphi_k& \mapsto H_k(I^k\times \varphi_1\times\dots\times \varphi_k).
\end{align*}

We then define for $k\geq 1$
\begin{align*}
    q_{k+1}\colon C_*(e^{-1}_\Omega(U_{x_0}),\mathcal{H}\cup e^{-1}_\Omega(V_{x_0}))&\otimes C_*(\Omega)^{\otimes k-1}\otimes C_*(\mathcal{P})\\
    &\to C_{k-n+*}(e^{-1}_\Lambda(U_\Delta),\mathcal{H}),\\
    \alpha\otimes \lambda_1\otimes \dots\otimes \lambda_k
    &\mapsto ((\tau_\Lambda\circ \eta\circ (e_\Omega\times (h_k\circ (\Theta_0^{\times k-1}\times \Theta_1))))\otimes \Phi)\circ \Delta(\alpha\times \dots\times \lambda_k).
\end{align*}
This restricts to the definition of $q_{k+1}$ in Section \ref{sec:ainfty struc}. An analogous computation to Theorem \ref{thm:q is an Ainfty morphism} shows the following:
\begin{Thm}
    The tuple $q_k$ defines a coherent chain homotopy $\boldsymbol{q}$ for the map
    \begin{align*}
        e_\Lambda^*\tau_\Lambda\cap \colon C_*(e^{-1}_\Lambda(U_\Delta),\mathcal{H}\cup e^{-1}_\Lambda(V_\Delta))&\to C_{-n+*}(e_\Lambda^{-1}(U_\Delta),\mathcal{H}).
    \end{align*}
    It restricts to the maps 
    \begin{align*}
        q_k\colon C_*(e^{-1}_\Omega(U_{x_0}),\mathcal{H}\cup e^{-1}_\Omega(V_{x_0}))\otimes C_*(\Omega)^{\otimes k}&\to C_{k-1-n+*}(e_\Omega^{-1}(U_{x_0}),\mathcal{H})
    \end{align*}
    as defined in \eqref{eq:def of qk}.
\end{Thm}

We can now compose the coherent chain homotopies of the above squares and thus find:

\begin{Cor}\label{cor:cch nu}
    The composition of the coherent chain homotopies of the above squares \eqref{sq:cchI}-\eqref{sq:cchVI} define a coherent chain homotopy 
    \begin{align*}
        \boldsymbol{\nu}\colon C_*(\Omega, x_0)\to  C_{1-n+*}((\Omega ,x_0)^2)
    \end{align*}
    for a map $\nu\colon C_*({\Lambda}M,M)\to C_{1-n+*}(\mathcal{F}\Lambda,\mathcal{H})$ such that $\mathrm{c}\circ \nu\colon C_*({\Lambda}M,M)\to C_{1-n+*}((\Lambda M,M)^2)$ computes the Goresky-Hingston coproduct. 
\end{Cor}
\begin{proof}
    This follows from the fact that we can compose coherent chain homotopies (see Remark \ref{rem:compose cch}).
\end{proof}

\begin{Thm}\label{thm:model of GH}
    Let $(M,x_0)$ be a smooth, pointed manifold of dimension $n$. The following diagram
    commutes
    \begin{center}
        \begin{tikzcd}
            H_*(M,C_*(\Omega M,x_0))\ar[d, "\widetilde{\nu}"] \ar[r, "\Psi"] & H_*(\Lambda M,M) \ar[dd, "\nu_{GH}", bend left=70] \ar[d] \\
            H_*(M,\Delta_* C_{1-n+*}((\Omega M,x_0)^2)\ar[d, "\Delta_*"]\ar[r, "\Psi"] & H_{1-n+*}(\mathcal{F}\Lambda, \mathcal{H})\ar[d, "\mathrm{c}_*"]\\
            H_*(M\times M, C_{1-n+*}((\Omega M,x_0)^2)\ar[r,"\Psi"]   &H_{1-n+*}((\Lambda M,M)^2).
        \end{tikzcd}
    \end{center}
    Moreover, the morphism of $\mathcal{A}_\infty$-modules $\boldsymbol{\nu}$ defining $\widetilde{\nu}$ is such that 
    \begin{align*}
        \nu_1\colon C_*(\Omega M,x_0)\to C_{1-n+*}((\Omega M,x_0)^2)
    \end{align*}
    computes the Goresky-Hingston coproduct on $\Omega M$.
\end{Thm}
\begin{proof}
    The statement follows from Corollary \ref{cor:coh chain hpty relative} and the fact that the coherent chain homotopy defined in Corollary \ref{cor:cch nu} restricts to the morphism of $\mathcal{A}_\infty$-module defined in \eqref{eq:Ainfty nu}.
\end{proof}

As explained in \cite[Section 4.2]{barraud2023morse}, there exists a canonical spectral sequence $E_{p,q}^r$ converging to $H_*(M,C_*(F))\cong H_*(E)$ which for $r\geq 2$ corresponds to the Leray-Serre spectral sequence \cite[Theorem 7.2.1]{barraud2023morse}. An analogous statement follows for the relative case and we thus have the following:
\begin{Cor}\label{cor:coprod spec seq}
    The Goresky-Hingston coproduct is computed on the second page of the Leray-Serre spectral sequence $E_{p,q}^2=H_p(M,H_q(\Omega M,x_0))$ as the following composition
    \begin{align*}
        H_p(M,H_q(\Omega M,x_0))\overset{\nu_1}{\to}H_p(M,H_{q+1-n}((\Omega M,x_0)^2))\overset{\Delta_*}{\to} H_p(M\times M,H_{q+1-n}((\Omega M,x_0)^2)).
    \end{align*}
    Dually, the Goresky-Hingston cohomology product is computed by first taking the cup product with local coefficients in $H_{q+1-n}((\Omega M,x_0)^2)$ and then computing the cohomology product on the based loop space.
\end{Cor}
\begin{proof}
    This comes down to the fact that the coproduct preserves 
    \begin{align}\label{eq:filtr GH domain}
        F_p(C_*(M,C_*(\Omega M, x_0)))=\bigoplus_{\substack{i+j=*\\ j\leq p}} C_{i}(\Omega, x_0)\otimes \Z \Crit_j(f).
    \end{align}
    By \cite[Remark 9.2.1]{barraud2023morse}, $\Delta_*$ preserves this filtration. It thus remains to show that the morphism of $\mathcal{A}_\infty$-modules
    \begin{align*}
        \widetilde{\nu}\colon C_*(M,C_*(\Omega M, x_0))\to C_*(M,C_{1-n+*}((\Omega M, x_0)^2))
    \end{align*}
    sends \eqref{eq:filtr GH domain} to
    \begin{align*}
        F_p(C_*(M,C_{1-n+*}((\Omega M, x_0)^2)))=\bigoplus_{\substack{i+j=*+1-n\\ j\leq p}} C_{i}((\Omega, x_0)^2)\otimes \Z \Crit_j(f).
    \end{align*}
    Indeed, we note that for $x\otimes \alpha\in F_p(C_*(M,C_*(\Omega M, x_0)))$
    \begin{align*}
        \widetilde{\nu}(\alpha\otimes x)=\sum_{k=0}^n\sum_{x=x_0>\dots>x_k}\pm \nu_{k+1}(\alpha\otimes m_{x_0,x_1}\otimes \dots\otimes m_{x_{k-1},x_k})\otimes x_k
    \end{align*}
    is in $F_p(C_*(M,C_{1-n+*}(\Omega M, x_0)))$ because $|x|\leq p$ and thus $|x_k|\leq |x|\leq p$.
\end{proof}

\begin{Rem}\label{rem:compute via spec seq}
    This corollary gives a description of the Goresky-Hingston coproduct on the associated graded of a filtration on $H_*(\Lambda M,M)$ up to extension:
    
    Because we are working with real coefficients, $H_*(\Lambda M,M)$ is determined by its dimension. In particular, we get a non-canonical isomorphism:
    \begin{align*}
        H_*(\Lambda M,M)\cong \bigoplus_{p+q=*} E^\infty_{p,q}
    \end{align*}
    where $E^\infty_{p,q}$ is the $E^\infty$-page of the Leray-Serre spectral with $E^2$-page $H_p(M,H_q(\Omega M,x_0))$.
    
    In particular, we get a filtration 
    \begin{align}\label{eq:Einfty-filtration}
        F^p_{*}=\bigoplus_{\substack{i+j=*\\ j\leq p}}  E^\infty_{i,j}
    \end{align}
    of $H_*(\Lambda M,M)$.
    
    In the proof of Corollary \ref{cor:coprod spec seq}, we showed that the coproduct preserves the filtration \eqref{eq:filtr GH domain} and thus also preserves the filtration \eqref{eq:Einfty-filtration}. In particular, the coproduct on $E^\infty_{p,q}$ lands in
    \begin{align*}
        \bigoplus_{\substack{p_1+p_2+q_1+q_2=p+q+1-n\\ p_1+p_2\leq p}} E_{p_1,q_1}^\infty\otimes E^\infty_{p_2,q_2}.
    \end{align*}
    The description of the coproduct on the spectral sequence as in Corollary \ref{cor:coprod spec seq} only describes the image of the coproduct in 
    \begin{align*}
        \bigoplus_{\substack{p_1+p_2+q_1+q_2=p+q+1-n\\ p_1+p_2= p}} E_{p_1,q_1}^\infty\otimes E^\infty_{p_2,q_2}.
    \end{align*}
    In the next subsection, we discuss an example where for degree reasons this image describes the whole coproduct. In other cases such as Example \ref{exmp:g>1 surface}, this description is missing all information.
\end{Rem}

\begin{Exmp}\label{exmp:g>1 surface}
    We consider the case for $M=\Sigma_g$ of a surface of genus $g\geq 2$.
    
    Because $M$ has contractible universal covering, we have $\pi_*(\Omega M)\cong \pi_{1+*}(M)\cong 0$ for all $*\geq 1$. Therefore $\Omega M$ has contractible path components and $H_*(\Omega M)$ is concentrated in degree $0$. In particular, the Goresky-Hingston coproduct on $H_*(\Omega M,x_0)$ is of degree $n-1=1$ and thus zero.
    
    Therefore Corollary \ref{cor:coprod spec seq} shows that the coproduct is zero on the Leray-Serre spectral sequence. This is in contrast to the fact that the coproduct on $H_*(\Lambda M,M)$ is not zero (see \cite{hartenstein2025string}).
\end{Exmp}

\subsection{The Sphere Quotiened by Finite Groups}\label{subsec:Sn/G}
In this subsection, we study the example of $M=S^n/G$ where $G$ is a finite group. We show that the coproduct with real coefficients is completely determined by the description on the $E^2$-page of the Leray-Serre spectral sequence as in Corollary \ref{cor:coprod spec seq}.\\

Let $G$ be a non-trivial, finite group that acts orientation-preserving on $S^n$. In particular $n$ is odd. We additionally assume that $n>1$. We set $M:=S^n/G$.

We note that the connected components of $\Lambda M$ correspond to conjugacy classes of $\pi_1(M)=G$. We denote for $[g]\in \mathrm{ccl}(G)$:
\begin{align*}
    \Lambda_{[g]}M:=\{\gamma\in \Lambda M\mid [\gamma]=h \text{ for some }h\in [g]\subseteq \pi_1(X)\}.
\end{align*}
We denote by $E^r_{pq}([g])$ the $E^\infty$-page of the Leray-Serre spectral sequence of the fibration 
\begin{center}
    \begin{tikzcd}
        \bigsqcup_{h\in [g]} \ar[r] & \ar[d]\Lambda_{[g]}M\\
        & M.
    \end{tikzcd}
\end{center}
In \cite{clivio2025GHonbased}, we showed that 
\begin{align*}
    E_{p,q}^\infty([g])\cong \begin{cases}
        \Z & \text{if } p=0,n \text{ and } q=k(n-1);\\
        H_p(G_g)&\text{if }1\leq p\leq n-1 \text{ and } q=k(n-1);\\
        0&\text{else,}
    \end{cases}
\end{align*}
where $G_g\subseteq G$ is the finite subgroup of elements in $G$ that commute with $g$. Thus $H_p(G_g)$ is zero after tensoring with $\R$ and we find 
\begin{align*}
    E_{p,q}^\infty\otimes \R\cong  \begin{cases}
        \R & \text{if } p=0,n\text{ and } q=k(n-1);\\
        0&\text{else.}
    \end{cases}
\end{align*}
This shows that 
\begin{align*}
    H_*(\Lambda M;\R)\cong \begin{cases}
        \bigoplus\limits_{[g]\in\mathrm{ccl}(G)}\R &\text{if } *=k(n-1), k(n-1)+n \text{ for some }k\geq 0;\\
        0&\text{ else.}
    \end{cases}
\end{align*}

\begin{Thm}\label{thm:GH on Sn/G}
    Let $G$ be a non-trivial, finite group acting orientation-preserving on $S^n$ for $n>1$. The homology of the free loop space of $M=S^n/G$ with real coefficients is given by
    \begin{align*}
        H_*(\Lambda M;\R)\cong&\  \bigoplus_{[g]\in \mathrm{ccl}(G)}\bigoplus_{k\geq 0} \R x_{[g],k}\oplus \R y_{[g],k} 
    \end{align*}
    with $|x_{[g],k}|=k(n-1)$ and $|y_{[g],k}|=k(n-1)+n$. The relative homology $H_*(\Lambda M,M;\R)$ has the same description except missing the $[g]=[1]$ and $k=0$ summand.
    
    The Goresky-Hingston coproduct lifts to a map 
    \begin{align*}
        \widehat\lor\colon H_*(\Lambda M;\R)&\to H_*(\Lambda M;\R)\otimes H_*(\Lambda M;\R)
    \end{align*}
    computed by
    \begin{align*}
        \widehat\lor(x_{[g],k})=\sum_{i+j=k-1} \sum_{h\in G} x_{[gh^{-1}],i}\otimes x_{[h],j}
    \end{align*}
    and
    \begin{align*}
        \widehat\lor(y_{[g],k})=\sum_{i+j=k-1} \sum_{h\in G} (x_{[gh^{-1}],i}\otimes y_{[h],j}+y_{[gh^{-1}],i}\otimes x_{[h],j}).
    \end{align*}
\end{Thm}
\begin{proof}
    The description as vector spaces follows from the discussion above. For the coproduct, we use Corollary \ref{cor:coprod spec seq} together with the computation of the coproduct on $\Omega M$ as in \cite{clivio2025GHonbased}. We showed that the Pontryagin ring structure on $H_*(\Omega M)$ is isomorphic $\Z[G][x]$ for a central element $x$ of degree $n-1$. Moreover, we showed that the coproduct on $H_*(\Omega M)$ lifts to a map
    \begin{align}\label{eq:GH on Omega Sn/G}
        \begin{split}
            \lor_\Omega\colon \Z[G][t]&\to \Z[G][t]\otimes \Z[G][t],\\
            gx^k&\mapsto \sum_{i+j=k-1}\sum_{h\in H} gh^{-1}x^i\otimes hx^j.
        \end{split}
    \end{align}
    This map is $\pi_1(M)$-equivariant and thus gives a map on the $E^2$-page:
    \begin{align*}
        \nu^{E^2}\colon H_*(M,H_{1-n+*}(\Omega M))&\to H_*(M,\Delta^*(H_*(\Omega M)\otimes H_*(\Omega M)))
    \end{align*}
    The $E^2$-page is isomorphic the $E^\infty$-page and thus we have a map $\nu^{E^\infty}$ on the $E^\infty$-page. By Corollary \ref{cor:coprod spec seq}, the map $\Delta_*\circ \nu^{E^\infty}$ is a lift of the Goresky-Hingston coproduct on the $E^\infty$-page. We now argue why $\Delta_*\circ \nu^{E^\infty}$ describes $\widehat\lor$ uniquely.
    
    We have $H_{n(k-1)}(\Lambda_{[g]}M)\cong E^\infty_{0,n(k-1)}$. By Remark \ref{rem:compute via spec seq} the image of $E^\infty_{0,n(k-1)}$ under $\widehat{\lor}$ lands in 
    \begin{align*}
        \bigoplus_{\substack{p_1+p_2+q_1+q_2=k(n-1)+1-n\\ p_1+p_2\leq 0}} E_{p_1,q_1}^\infty\otimes E^\infty_{p_2,q_2}.
    \end{align*}
    Because $E_{p,q}^\infty\cong 0$ for $p<0$, we have that the image of $E^\infty_{0,n(k-1)}$ under $\widehat{\lor}$ lands in 
    \begin{align*}
        \bigoplus_{\substack{p_1+p_2+q_1+q_2=k(n-1)+1-n\\ p_1+p_2= 0}} E_{p_1,q_1}^\infty\otimes E^\infty_{p_2,q_2}
    \end{align*}
    and as argued in Remark \ref{rem:compute via spec seq} is determined by the map it induces on the $E^\infty$-page i.e.\ $\Delta_*\circ \nu^{E^\infty}$.
    
    Similarly, $H_{k(n-1)+n}(\Lambda M)\cong E^\infty_{n,k(n-1)}$ lands in 
    \begin{align}\label{eq:image of y[g],k}
        \bigoplus_{\substack{p_1+p_2+q_1+q_2=k(n-1)+n+1-n\\ p_1+p_2\leq n}} E_{p_1,q_1}^\infty\otimes E^\infty_{p_2,q_2}.
    \end{align}
    But $H_*(\Lambda M)$ is concentrated in degree $*\equiv 0,1$ (mod $n-1$). Therefore \eqref{eq:image of y[g],k} is concentrated in degrees where $p_1+q_1\equiv 1$ (mod $n-1$) and $p_2+q_2\equiv 0$ (mod $n-1$) or vice-versa. But $E^\infty_{p,q}$ is concentrated in degrees where $q\equiv 0$ (mod $n-1$). This shows that $p_1+p_2\equiv 1$ (mod $n-1$) which in turn shows $p_1+p_2=n$.
    
    This finally shows that $\widehat\lor$ sends $H_{k(n-1)+n}(\Lambda M)$ to 
    \begin{align*}
        \bigoplus_{\substack{p_1+p_2+q_1+q_2=k(n-1)+n+1-n\\ p_1+p_2= n}} E_{p_1,q_1}^\infty\otimes E^\infty_{p_2,q_2}.
    \end{align*}
    and as argued in Remark \ref{rem:compute via spec seq} is determined by $\Delta_*\circ \nu^{E^\infty}$.
    
    We now compute what $\widehat{\lor}=\Delta_*\circ \nu^{E^\infty}$ does on $E^\infty_{p,q}\cong E^2_{p,q}$. On the $E^2$-page we have that $E_{0,k(n-1)}([g])$ is represented by $ gx^k\otimes x_0$ and $E_{n,k(n-1)}([g])$ is represented by $gt^k\otimes [M]$. We call the corresponding elements in $H_*(\Lambda M)$ $x_{[g],k}$ and $y_{[g],k}$, respectively.
    
    Using \eqref{eq:GH on Omega Sn/G}, $\Delta_*(x_0)=x_0\otimes x_0$ and $\Delta_*([M])=x_0\otimes [M]+[M]\otimes x_0$, we compute
    \begin{align*}
        \widehat{\lor}(x_{[g],k})=\Delta_*(\lor_\Omega(gt^k)\otimes x_0)&=\Delta_*\left(\left(\sum_{i+j=k-1}\sum_{h\in H} gh^{-1}x^i\otimes hx^j\right)\otimes x_0\right)\\
        &=\sum_{i+j=k-1}\sum_{h\in H} ( gh^{-1}x^i\otimes x_0)\otimes ( hx^j\otimes x_{0})\\
        &=\sum_{i+j=k-1} \sum_{h\in G} x_{[gh^{-1}],i}\otimes x_{[h],j}
    \end{align*}
    and 
    \begin{align*}
        \widehat{\lor}(y_{[g],k})&=\Delta_*( \lor_\Omega(gx^k)\otimes [M])\\
        &=\Delta_*\left(\left(\sum_{i+j=k-1}\sum_{h\in H} gh^{-1}x^i\otimes hx^j\right)\otimes [M]\right)\\
        &=\sum_{i+j=k-1}\sum_{h\in H} \left(( gh^{-1}x^i\otimes x_0)\otimes ( ht^j\otimes [M])+( gh^{-1}x^i\otimes [M])\otimes ( ht^j\otimes x_0)\right)\\
        &=\sum_{i+j=k-1} \sum_{h\in G} (x_{[gh^{-1}],i}\otimes y_{[h],j}+y_{[gh^{-1}],i}\otimes x_{[h],j}).
    \end{align*}
\end{proof}

\subsection{Invariance under Maps}\label{subsec:inv under maps}
In this subsection, we use our main theorem to study the invariance or non-invariance of the Goresky-Hingston coproduct loop space under maps of manifolds.\\

For this we import the following result that shows that the Goresky-Hingston coproduct commutes with degree $1$ maps. This was originally proved in \cite{hingston2010loop} and we state and prove here the result for the coproduct as we defined it:
\begin{Prop}\label{prop:degree 1 commutes with nu1}
    Let $f\colon (M,x_0)\to (N,y_0)$ be a degree $1$ map and denote 
    \begin{align*}
        \boldsymbol{\nu}^N\colon  C_{*}(\Omega N, y_0)\to C_{1-n+*}((\Omega N,y_0)^2)
    \end{align*}
    for a morphism of $\mathcal{A}_\infty$-modules as in \eqref{eq:Ainfty nu}. There exists a morphism of $\mathcal{A}_\infty$-modules as constructed in \eqref{eq:Ainfty nu} 
    \begin{align*}
        \boldsymbol{\nu}^M\colon  C_{*}(\Omega M, x_0)\to C_{1-n+*}((\Omega N,x_0)^2)
    \end{align*}
    a map $g\colon (M,x_0)\to (N,y_0)$ homotopic to $f$ which satisfies
    \begin{align*}
        \nu_1^N\circ \Omega g=(\Omega g\times \Omega g)\circ \nu_1^M.
    \end{align*}
\end{Prop}
\begin{proof}
    By \cite{hopf1930topologie} for $n\neq 2$ and \cite{kneser1928glattung} for $n=2$, we can find $g$ homotopic to $f$ such that $g^{-1}(U_{y_0})=U_{x_0}$. Moreover, we can assume that $g|_{U_{x_0}}$ commutes with the Riemannian metric. We then have $\tau_{x_0}:=g^*\tau_{y_0}$ is a Thom class that satisfies Corollary \ref{cor:invariant taux0}. We can thus define $\boldsymbol{\nu}^M$ using this Thom class.
    
    In particular, we find that $\Omega g$ commutes with $\nu_1^M$ and $\nu_1^N$, because the cap product commutes with $\Omega g$ by naturality of the cap product and
    \begin{align*}
        e_\Omega^*\tau_{x_0}=e_\Omega^*g^*\tau_{y_0}=(\Omega g)^*e_\Omega^*\tau_{y_0}.
    \end{align*}
\end{proof}

\begin{Cor}\label{cor:commute with deg1}
    Let $f\colon (M,x_0)\to (N,y_0)$ be a degree $1$ map and denote $\nu^M$ and $\nu^N$ the Goresky-Hingston coproduct on $H_*(\Lambda M,M)$ and $H_*(\Lambda N, N)$, respectively. On the $E^2$-page of the Serre spectral sequence it holds
    \begin{align*}
        \nu^N\circ \Lambda f=(\Lambda f\times \Lambda f)\circ \nu^M
    \end{align*}
    as maps of spectral sequences.
\end{Cor}
\begin{proof}
    The diagonal $\Delta_*$ commutes with all maps $f\colon (M,x_0)\to (N,y_0)$. Moreover, on the $E^1$-page of the canonical spectral sequence of $C_*(M,C_*(\Omega,x_0))$, and thus on all higher pages, $\widetilde{\nu}$ is computed by 
    \begin{align*}
        H_q(\Omega,x_0)\otimes \Z\Crit_p(f)&\to H_{q+1-n}((\Omega,x_0)^2)\otimes \Z \Crit_p(f),\\
        \alpha\otimes x&\mapsto \nu_1(\alpha)\otimes x,
    \end{align*}
    where $\nu_1$ is the coproduct on the based loop space. By Proposition \ref{prop:degree 1 commutes with nu1}, $f$ commutes with $\nu_1$ in homology.
\end{proof}

\begin{Rem}\label{rem:commute with maps}
    This result is notable because it is in contrast to the fact that the Goresky-Hingston coproduct does not commute in general with maps that are not simple homotopies (see \cite{naef2021string, naef2024simple, wahl2019invariance, kenigsberg2024obstructions}).
    
    In fact for a map of manifolds $f\colon M\to N$, we can ask whether or not it commutes with the following variants of the Goresky-Hingston coproduct.
    \begin{center}
        \begin{tikzcd}
            \makecell{\Lambda f\text{ commutes with}\\ \text{the coproduct on }H_*(\Lambda M,M;\Z) \\ \text{with integer coefficients;}}\ar[d, Rightarrow] \\
            \makecell{\Lambda f\text{ commutes with}\\ \text{the coproduct on }H_*(\Lambda M,M;\R) \\ \text{with real coefficients;}}\ar[d, Rightarrow]\\
            \makecell{\Lambda f\text{ commutes with}\\ \text{the coproduct on the }E^\infty\text{-page} \\ \text{of the Leray-Serre spectral sequence}\\ \text{with real coefficients.}}
        \end{tikzcd}
    \end{center}

    Naef proved in \cite{naef2024simple} that the top statement is true if $f$ is a simple homotopy and showed in \cite{naef2024simple} that there exist examples of homotopy equivalences such that top statement is false.

    Corollary \ref{cor:commute with deg1} shows that the bottom statement is true if $f$ is a degree $1$ map.

    To our knowledge, there exist no known examples of degree $1$ maps where the statement in the middle is false. In particular, Theorem \ref{thm:GH on Sn/G} shows that, for maps between lens spaces, the middle statement holds. Lens spaces were used by Naef in \cite{naef2021string} to give a homotopy equivalence where the top left statement does not hold.
\end{Rem}

\appendix

\section{Inverting Quasi-isomorphisms of $\mathcal{A}_\infty$-Modules}\label{app:inv qi of Ainfty modules}

In this section, we prove that, for an $\infty$-quasi-isomorphism of $\mathcal{A}_\infty$-modules $\boldsymbol{f}\colon \mathcal{M}_*\to \mathcal{N}_*$, we can find a morphism $\boldsymbol{g}\colon \mathcal{N}_*\to \mathcal{M}_*$ where $g_1$ is an inverse to $f_1$ on homology. This is a well-known result for $\mathcal{A}_\infty$-algebras (see \cite{lefevre2003infini} and \cite[Theorem 10.4.4]{loday2012algebraic} for algebras over any Koszul operad). However, the statement for $\mathcal{A}_\infty$-modules, while known to experts, does not have a proof in the literature to our knowledge.

\begin{Def}\label{def:infty iso}
    Let $(A_*,d)$ be a (strict) differential graded algebra. A morphism of $\mathcal{A}_\infty$-modules $\boldsymbol{f}\colon \mathcal{M}_*\to \mathcal{N}_*$ given by $f_k\colon \mathcal{M}_*\otimes A_*^{\otimes k-1}\to \mathcal{N}_{k-1+*}$ is an \textit{$\infty$-isomorphism} if $f_1$ is an isomorphism and an \textit{$\infty$-quasi-isomorphism} if $f_1$ is a quasi-isomorphism.
\end{Def}

The goal of this section is to find for an $\infty$-quasi-isomorphism $\boldsymbol{f}\colon \mathcal{M}_*\to \mathcal{N}_*$ a map $\boldsymbol{g}\colon \mathcal{N}_*\to \mathcal{M}_*$ such that $H_*(g_1)$ is an inverse to $H_*(f_1)$. Our strategy is a classical strategy of $\mathcal{A}_\infty$-morphisms as for example in \cite[Theorem 10.4.4]{loday2012algebraic}:
\begin{enumerate}[(1)]
    \item invert $\infty$-isomorphisms;
    \item prove a version of the Homotopy Transfer Theorem for $\mathcal{A}_\infty$-modules.
\end{enumerate}

\subsection{Inverting $\infty$-Isomorphisms}

\begin{Prop}\label{prop:inverting infty-iso}
    Let $(A_*,d)$ be a (strict) differential graded algebra. Let $\boldsymbol{f}\colon \mathcal{M}_*\to \mathcal{N}_*$ be an $\infty$-isomorphism of $A$-modules given by
    \begin{align*}
        f_k\colon \mathcal{M}_*\otimes A_*^{\otimes k-1}\to \mathcal{N}_{k-1+*}.
    \end{align*}
    The map $\boldsymbol{g}\colon \mathcal{N}_*\to \mathcal{M}_*$ given by $g_1:=f_1^{-1}$ and for $k\geq 2$:
    \begin{align*}
        g_{N+1}:=-f_1^{-1}\circ \left(\sum_{r=0}^{N-1} (-1)^{r(N-r)}f_{N-r+1}\circ \left(g_{r+1}\otimes \id^{\otimes N-r}\right)\right)
    \end{align*}
    defines an inverse to $\boldsymbol{f}$. In other words, it holds $\boldsymbol{g}\circ \boldsymbol{f}=\id_{\mathcal{M}_*}$ and $\boldsymbol{f}\circ \boldsymbol{g}=\id_{\mathcal{N}_*}$.
\end{Prop}
\begin{proof}
    We have to check that the above formulas define a morphism of $\mathcal{A}_\infty$-modules and that $\boldsymbol{g}\circ \boldsymbol{f}=\id_{\mathcal{M}_*}$ and $\boldsymbol{f}\circ \boldsymbol{g}=\id_{\mathcal{N}_*}$ holds. We denote $m^\mathcal{M}_k$ and $m^{\mathcal{N}}_k$ for the $\mathcal{A}_\infty$-module structure on $\mathcal{M}_*$ and $\mathcal{N}_*$, respectively.
    
    We first check that $\boldsymbol{g}\circ \boldsymbol{f}=\id_{\mathcal{M}_*}$ and $\boldsymbol{f}\circ \boldsymbol{g}=\id_{\mathcal{N}_*}$. We recall that composition on morphisms of $\mathcal{A}_\infty$-modules is defined via
    \begin{align*}
        (\boldsymbol{\eta}\circ \boldsymbol{\zeta})_{N+1}=\sum_{k=0}^{N} (-1)^{k(N-k)} \eta_{k+1}(\zeta_{N-k+1}\otimes \id^{\otimes k}).
    \end{align*}
    We have $g_1\circ f_1=\id$ and $f_1\circ g_1=\id$ and for $N\geq 1$
    \begin{align*}
        (\boldsymbol{f}\circ \boldsymbol{g})_{N+1}=f_1\circ g_{N+1}+ \sum_{r=0}^{N-1} (-1)^{r(N-r)}f_{N-r+1}\circ \left(g_{r+1}\otimes \id^{\otimes N-r}\right)=0
    \end{align*}
    by construction. To show that $(\boldsymbol{g}\circ \boldsymbol{f})_{N+1}=0$, we define for $k_1,\dots,k_l\geq 1$:
    \begin{align*}
        \varphi_{k_1,\dots,k_l}:=f_1^{-1}  f_{k_1+1} (f_1^{-1}f_{k_2+1}\otimes \id^{\otimes k_1}) \dots (f_1^{-1}f_{k_l+1}\otimes \id^{\otimes k_1+\dots +k_{l-1}})(f_1^{-1}\otimes \id^{\otimes k_1+\dots+k_l})
    \end{align*}
    as in Figure \ref{fig:phik1...kl}.\\
    \begin{figure}[H]
        \centering
        \includegraphics[width=0.2\linewidth]{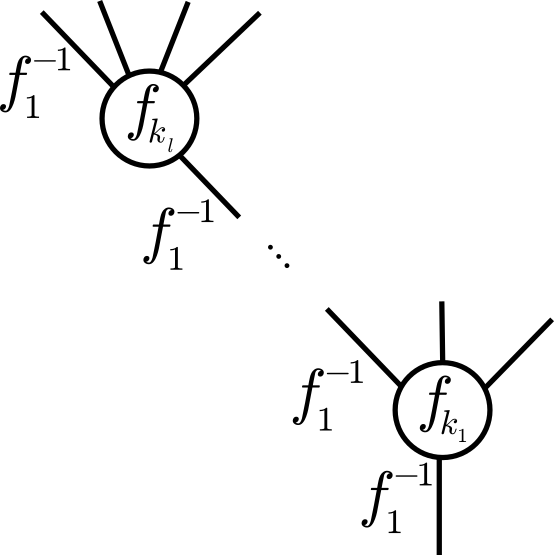}
        \caption{The map $\varphi_{k_1,\dots,k_l}$.}
        \label{fig:phik1...kl}
    \end{figure}
    We then find, by induction on $N$, the formula
    \begin{align*}
        g_{N+1}&=\sum_{k_1+\dots+k_l=N} (-1)^{l+\sum_{i<j}k_ik_j} \varphi_{k_1,\dots,k_l}\\
        &=\sum_{k_l=1}^N(-1)^{k_l(N-k_l)}\sum_{k_1+\dots+k_{l-1}=N-k_l} (-1)^{l+\sum_{i<j<l}k_ik_j} \varphi_{k_1,\dots,k_l}\\
        &=\left(\sum_{k_l=1}^{N}(-1)^{k_l(N-k_l)}g_{N-k_l+1}(f_{k_l+1}\otimes \id^{\otimes N-k_l}) \right)(f_1^{-1}\otimes \id^{\otimes N}).
    \end{align*}
    From this, it follows that $(\boldsymbol{g}\circ \boldsymbol{f})_{N+1}=0$.
    
    We now check the $\mathcal{A}_\infty$-equation which reads for $N\geq 0$ as:
    \begin{align}\label{eq:Ainfty morph eq}
        \begin{split}
            dg_{N+1}=&-\sum_{r = 1}^N (-1)^{(N-r)r} m^{\mathcal{M}}_{r+1}(g_{N-r+1} \otimes \id^{\otimes r})\\
            &+\sum_{r=1}^N (-1)^{(N-r)(r+1)} g_{N-r+1}(m_{r+1}^{\mathcal{N}} \otimes  \id^{\otimes N-r}) \\
            &+ (-1)^N g_{N+1}d\\
            &+ \sum_{r=1}^{N-1} (-1)^r g_{N}(\id^{\otimes r} \otimes \mu \otimes \id^{\otimes N-1-r}).
        \end{split}
    \end{align}
    We prove it by induction on $N$. For $N=0$, the statement is that $g_1=f_1^{-1}$ is a chain map. This is true because $f_1$ is a chain map.
    
    Because $\boldsymbol{f}$ is a morphism of $\mathcal{A}_\infty$-modules, we find:
    \begin{align*}
        (-f_1^{-1})df_{N-k+1}=&-m_{N-k+1}^\mathcal{M}\\
        &+\sum_{r=1}^{N-k} (-1)^{r(N-k-r)} f_1^{-1} m_{r+1}^\mathcal{N} (f_{N-k-r+1}\otimes \id^{\otimes r})\\
        &-\sum_{r=1}^{N-k-1} (-1)^{r(N-k-r+1)}f_1^{-1}f_{r+1}(m_{N-k-r+1}^\mathcal{M}\otimes \id^{\otimes r})\\
        &-(-1)^{N-k} f_1^{-1}f_{N-k+1} d\\
        &-\sum_{r=1}^{N-k-1} (-1)^r f_1^{-1}f_{N-k}(\id ^{\otimes r}\otimes \mu \otimes \id^{\otimes N-k-r-1})
    \end{align*}
    and thus
    \begin{align}\label{eq:dgN+1}
        \begin{split}
            dg_{N+1}=& -\sum_{k=0}^{N-1} (-1)^{k(N-k)}m_{N-k+1}^\mathcal{M}(g_{k+1}\otimes \id^{\otimes N-k})\\
            &+\sum_{k=0}^{N-1} (-1)^{k(N-k)}\sum_{r=1}^{N-k} (-1)^{r(N-k-r)} f_1^{-1} m_{r+1}^\mathcal{N} (f_{N-k-r+1}\otimes \id^{\otimes r})(g_{k+1}\otimes \id^{\otimes N-k})\\
            &-\sum_{k=0}^{N-1} (-1)^{k(N-k)}\sum_{r=1}^{N-k-1} (-1)^{r(N-k-r+1)}f_1^{-1}f_{r+1}(m_{N-k-r+1}^\mathcal{M}\otimes \id^{\otimes r})(g_{k+1}\otimes \id^{\otimes N-k})\\
            &-\sum_{k=0}^{N-1} (-1)^{(k+1)(N-k)} f_1^{-1}f_{N-k+1} (dg_{k+1}\otimes \id^{\otimes N-k})\\
            &+(-1)^N\sum_{k=0}^{N-1} (-1)^{k(N-k)} (-f_1^{-1})f_{N-k+1} (g_{k+1}\otimes \id^{\otimes N-k})(\id^{\otimes k+1}\otimes d)\\
            &+\sum_{r=k+1}^{N-1} (-1)^rg_N(\id ^{\otimes r}\otimes \mu \otimes \id^{\otimes N-r-1})
        \end{split}
    \end{align}
    The second term can be rewritten as
    \begin{align*}
        \sum_{r=1}^{N}(-1)^{r(N-r)}f_1^{-1} m_{r+1}^\mathcal{N} (\left(\boldsymbol{f}\circ \boldsymbol{g}\right)_{N-r+1}\otimes \id^{\otimes r})=g_1m_{N+1}^\mathcal{N}
    \end{align*}
    because $\boldsymbol{f}\circ \boldsymbol{g}=\id$.
    
    The third and fourth term in \eqref{eq:dgN+1} evaluate to
    \begin{align*}
        % -\sum_{k=0}^{N-1}& (-1)^{k(N-k)}\sum_{r=1}^{N-k-1} (-1)^{r(N-k-r+1)}f_1^{-1}f_{r+1}(m_{N-k-r+1}^\mathcal{M}\otimes \id^{\otimes r})(g_{k+1}\otimes \id^{\otimes N-k})\\
        % -\sum_{k=0}^{N-1}& (-1)^{(k+1)(N-k)} f_1^{-1}f_{N-k+1} (dg_{k+1}\otimes \id^{\otimes N-k})\\
        % =&-\sum_{r=0}^{N-1}(-1)^{rN}\sum_{k=0}^{N-r-1}(-1)^{(r-k)(N-k)}f_1^{-1}f_{r+1}(m_{N-k-r+1}^\mathcal{M}\otimes \id^{\otimes r})(g_{k+1}\otimes \id^{\otimes N-k})\\
        % &-\sum_{k=0}^{N-1} (-1)^{(k+1)(N-k)} f_1^{-1}f_{N-k+1} (dg_{k+1}\otimes \id^{\otimes N-k})\\
        % =&
        -\sum_{k=0}^{N-1}(-1)^{(k+1)(N-k)}f_1^{-1}f_{N-k+1} \left(\left(\sum_{r=1}^{k}(-1)^{(k-r)r }m_{r+1}^\mathcal{M}(g_{k-r+1}\otimes \id^{\otimes k}) + d g_{k+1} \right)\otimes \id^{\otimes N-k}\right).
    \end{align*}
    We can apply the induction hypothesis and find:
    \begin{align*}
        -\sum_{k=0}^{N-1}&(-1)^{(k+1)(N-k)}f_1^{-1}f_{N-k+1} \left(\left(\sum_{r=1}^{k}(-1)^{(k-r)r }m_{r+1}^\mathcal{M}(g_{k-r+1}\otimes \id^{\otimes k}) + d g_{k+1} \right)\otimes \id^{\otimes N-k}\right)\\
        % &=-\sum_{k=0}^{N-1}(-1)^{(k+1)(N-k)}f_1^{-1}f_{N-k+1} \left(\left(\sum_{r=1}^k(-1)^{(k-r)(r+1)} g_{k-r+1}(m_{r+1}^{\mathcal{N}}\otimes \id^{\otimes k-r}+(-1)^k g_{k+1}d +\sum_{r=1}^{k-1} g_r(\id^{\otimes r}\otimes \mu\otimes \id^{\otimes k-1-r}) \right)\otimes \id^{\otimes N-k}\right)\\
        % =
        =&\ \sum_{r=1}^{N-1}(-1)^{N(N-r)}(-f_1^{-1})((\boldsymbol{f}\circ \boldsymbol{g})_{N-r+1}-f_1g_{N-r+1})\circ (m_{r+1}^{\mathcal{N}}\circ \id^{\otimes N-r})\\
        &+(-1)^N\sum_{k=0}^{N-1} (-1)^{k(N-k)} (-f_1^{-1})f_{N-k+1} (g_{k+1}\otimes \id^{\otimes N-k})(d\otimes \id^{\otimes N-k})\\
        &+\sum_{r=1}^{N-k-1} (-1)^rg_N(\id ^{\otimes r}\otimes \mu \otimes \id^{\otimes N-r-1}).
    \end{align*}
    Plugging this into \eqref{eq:dgN+1}, we find:
    \begin{align*}
        dg_{N+1}=& -\sum_{k=0}^{N-1} (-1)^{k(N-k)}m_{N-k+1}^\mathcal{M}(g_{k+1}\otimes \id^{\otimes N-k})\\
        &+\sum_{r=1}^{N}(-1)^{N(N-r)}g_{N-r+1}\circ (m_{r+1}^{\mathcal{N}}\circ \id^{\otimes N-r})\\
        &+(-1)^N\sum_{k=0}^{N-1} (-1)^{k(N-k)} (-f_1^{-1})f_{N-k+1} (g_{k+1}\otimes \id^{\otimes N-k})(d\otimes \id^{\otimes N-k}+\id^{\otimes k+1}\otimes d)\\
        &+\sum_{r=1}^{N-k-1} (-1)^rg_N(\id ^{\otimes r}\otimes \mu \otimes \id^{\otimes N-r-1})+\sum_{r=k+1}^{N-1} (-1)^rg_N(\id ^{\otimes r}\otimes \mu \otimes \id^{\otimes N-r-1})
    \end{align*}
    which simplifies to \eqref{eq:Ainfty morph eq}.
\end{proof}

\subsection{Homotopy Transfer Theorem for $\mathcal{A}_\infty$-Modules}
In this subsection, we prove that we can transport a (right) module structure along a homotopy retract to an $\mathcal{A}_\infty$-module structure.

\begin{Def}\label{def:hpty retract}
    The data $(i,p,h)$ as in
    \begin{center}
        \begin{tikzcd}
            \mathcal{M}_* \ar[r, "p", shift left=1] \ar[loop left, "h"] & \mathcal{N}_* \ar[l, "i", shift left=1]
        \end{tikzcd}
    \end{center}
    exhibits $(\mathcal{N}_*,d_\mathcal{N})$ as a \textit{homotopy retract} of $(\mathcal{M}_*,d_\mathcal{M})$ if $i$ and $p$ are chain maps, $i$ is a quasi-isomorphism and $d_\mathcal{M}h+hd_\mathcal{M}=\id_\mathcal{M}-ip$.
\end{Def}

\begin{Prop}[Homotopy Transfer Theorem for $\mathcal{A}_\infty$-modules]\label{prop:htt for Ainfty modules}
    Let $(i,p,h)$ exhibit $(\mathcal{N}_*,d_\mathcal{N})$ as a homotopy retract of $(\mathcal{M}_*,d_\mathcal{M})$. Let $(A_*,d_A,\mu_A)$ be a (strict) differential graded algebra such that $\mu_\mathcal{M}\colon \mathcal{M}_*\otimes A_*\to \mathcal{M}_*$ exhibits $\mathcal{M}_*$ as a right $A_*$-module.
    
    The maps $m_1:=d_\mathcal{N}$, $m_2:=p\circ  \mu_\mathcal{M}(i\otimes \id)$ and
    \begin{align*}
        m_k:=(-1)^{\frac{k(k-1)}{2}}p \mu_\mathcal{M}(\overset{k-2}{\overbrace{h \mu_\mathcal{M}(\dots (h\mu_\mathcal{M}}}(i\otimes \id)\otimes \id)\otimes\dots)\otimes\id ),
    \end{align*}
    for $k\geq 2$, define a $\mathcal{A}_\infty$-module structure on $\mathcal{N}_*$.
    
    Moreover, the maps $\boldsymbol{p}\colon \mathcal{M}_*\to \mathcal{N}_*$ and $\boldsymbol{i}\colon \mathcal{N}_*\to \mathcal{M}_*$ given by $p_1:=p$, $i_1:=i$
    \begin{align*}
        p_k&:=(-1)^{\frac{(k-1)(k-2)}{2}}p \mu_\mathcal{M}(\overset{k-1}{\overbrace{h \mu_\mathcal{M}(\dots (h\mu_\mathcal{M}(h}}\otimes \id)\otimes \id)\otimes\dots)\otimes\id ),\\
        i_k&:=(-1)^{\frac{k(k-1)}{2}}\overset{k-1}{\overbrace{h \mu_\mathcal{M}(h \mu_\mathcal{M}(\dots (h\mu_\mathcal{M}}}(i\otimes \id)\otimes \id)\otimes\dots)\otimes\id ),
    \end{align*}
    for $k\geq 2$, define morphisms of $\mathcal{A}_\infty$-modules.
\end{Prop}
\begin{figure}[H]
        \centering
        \begin{minipage}{0.3\textwidth}
            \centering
            \captionsetup{width=0.85\textwidth}
            \includegraphics[width=0.85\textwidth]{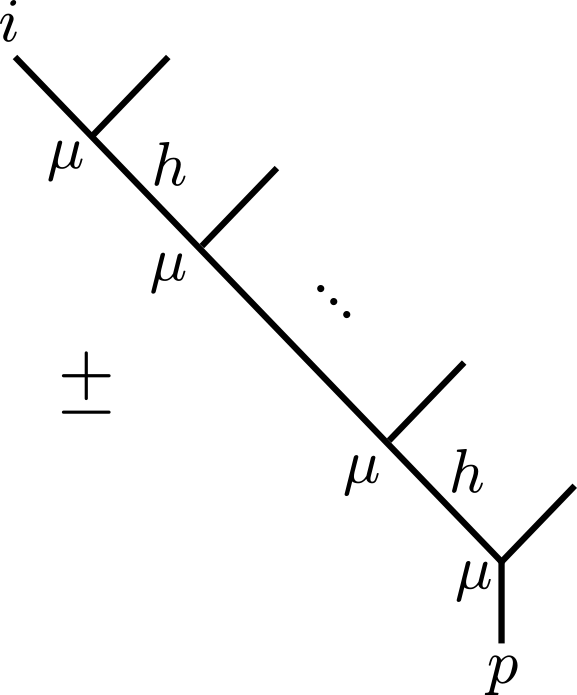} 
            \caption{The map $m_k$.}\label{fig:mk}
        \end{minipage}\hfill
        \begin{minipage}{0.3\textwidth}
            \centering
            \captionsetup{width=0.85\textwidth}
            \includegraphics[width=0.85\textwidth]{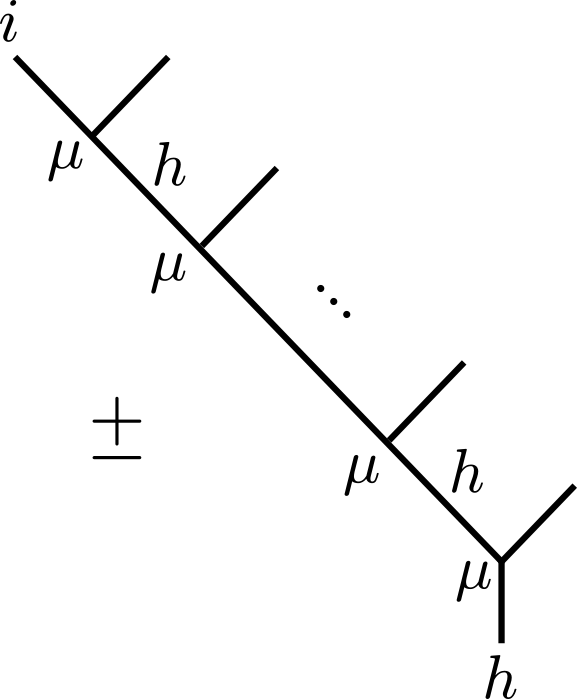} 
            \caption{The map $i_k$.}\label{fig:ik}
        \end{minipage}\hfill
        \begin{minipage}{0.3\textwidth}
            \centering
            \captionsetup{width=0.85\textwidth}
            \includegraphics[width=0.85\textwidth]{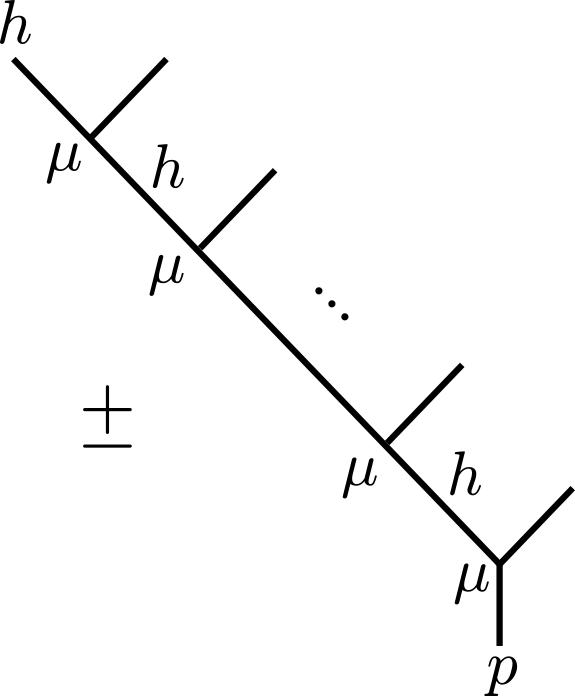} 
            \caption{The map $p_k$.}\label{fig:pk}
        \end{minipage}
    \end{figure}
\begin{proof}
    We need to check the $\mathcal{A}_\infty$-equation for $m_k$, which reads as:
    \begin{align*}
        d m_N+(-1)^{N+1}m_Nd&+\sum_{t=1}^{N-2} (-1)^{(N-t)t} m_{t+1}(m_{N-t}\otimes \id^{\otimes t})\\
        &+\sum_{r=1}^{N-2}(-1)^rm_{N-1}(\id^{\otimes r}\otimes \mu\otimes \id^{\otimes N-r-2})=0.
    \end{align*}
    Using that $dh=\id_\mathcal{M}-ip-hd$, we compute as in Figure \ref{fig:dh in tree}
    \begin{align*}
        p \mu_\mathcal{M}(&h \mu_\mathcal{M}(\dots dh\mu_\mathcal{M}(\dots (h\mu_\mathcal{M}(i\otimes \id)\otimes \id)\dots \otimes \id)\dots \otimes \id )\otimes\id )\\
        =&\ p \mu_\mathcal{M}(h \mu_\mathcal{M}(\dots \mu_\mathcal{M}(\dots (h\mu_\mathcal{M}(i\otimes \id)\otimes \id)\dots \otimes \id)\dots \otimes \id )\otimes\id )\\
        &-p \mu_\mathcal{M}(h \mu_\mathcal{M}(\dots ip\mu_\mathcal{M}(\dots (h\mu_\mathcal{M}(i\otimes \id)\otimes \id)\dots \otimes \id)\dots \otimes \id )\otimes\id )\\
        &-p \mu_\mathcal{M}(h \mu_\mathcal{M}(\dots hd\mu_\mathcal{M}(\dots (h\mu_\mathcal{M}(i\otimes \id)\otimes \id)\dots \otimes \id)\dots \otimes \id )\otimes\id ).       
    \end{align*}
    \begin{figure}[H]
        \centering
        \includegraphics[width=0.7\linewidth]{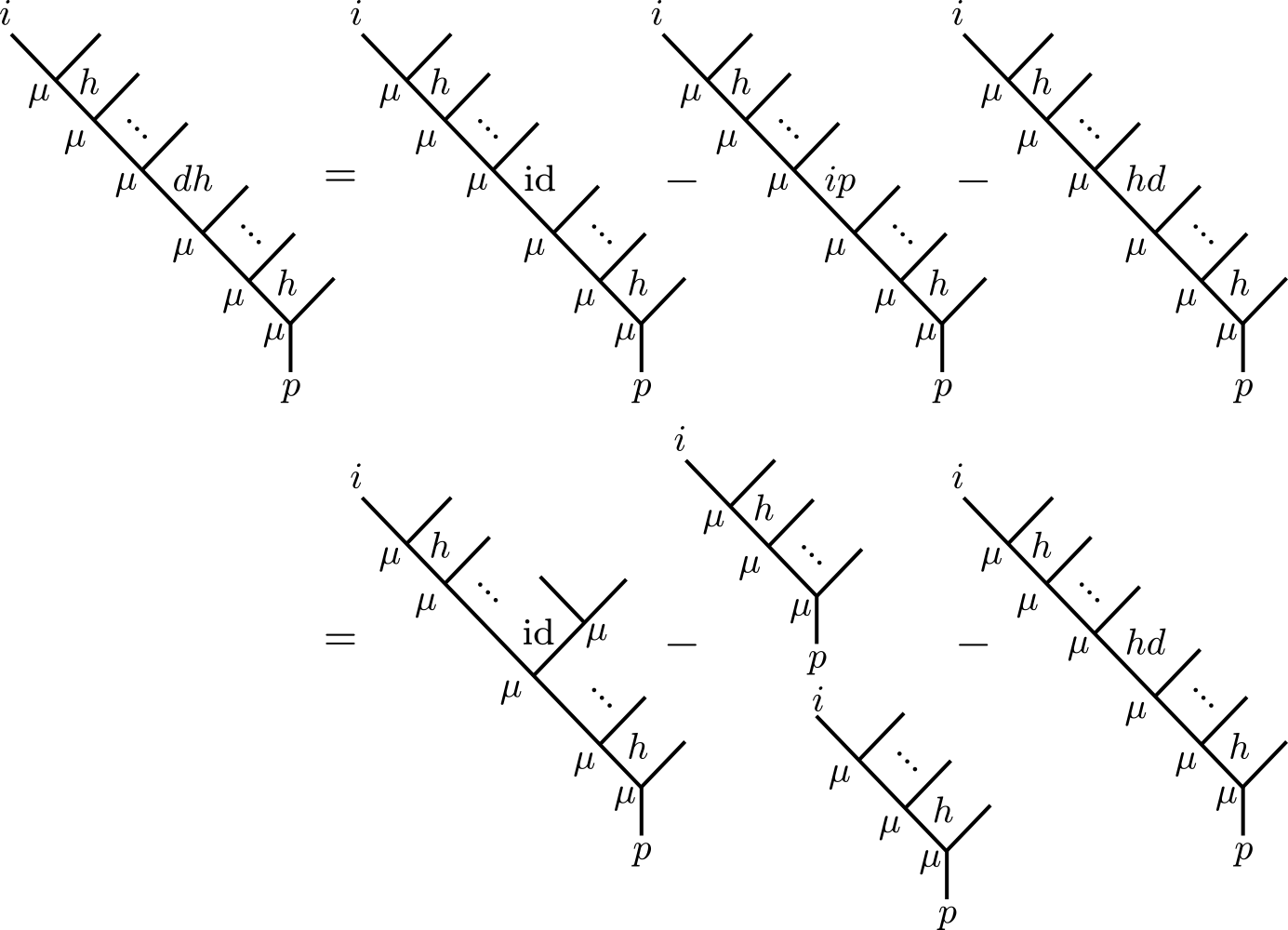}
        \caption{Putting a $dh$ on an edge in the tree gives the same tree with $\id_\mathcal{M}-ip-hd$ at the edge.}
        \label{fig:dh in tree}
    \end{figure}
    It therefore follows that
    \begin{align*}
        (-1)^{\frac{N(N-1)}{2}}d m_N=&\  d p \mu_\mathcal{M}(h \mu_\mathcal{M}(\dots (h\mu_\mathcal{M}(i\otimes \id)\otimes \id)\dots \otimes \id)\otimes\id )\\
        =&\ p \mu_\mathcal{M}(d {h \mu_\mathcal{M}(\dots (h\mu_\mathcal{M}}(i\otimes \id)\otimes \id)\dots \otimes \id)\otimes\id )\\
        &+(-1)^{N-2}p \mu_\mathcal{M}( {h \mu_\mathcal{M}(\dots (h\mu_\mathcal{M}}(i\otimes \id)\otimes \id)\dots \otimes \id)\otimes d)\\
        =&\ p \mu_\mathcal{M}( \mu_\mathcal{M}(\dots (h\mu_\mathcal{M}(i\otimes \id)\otimes \id)\dots \otimes \id)\otimes\id )\\
        &-p \mu_\mathcal{M}( ip\mu_\mathcal{M}(\dots (h\mu_\mathcal{M}(i\otimes \id)\otimes \id)\dots \otimes \id)\otimes\id )\\
        &-p \mu_\mathcal{M}( hd\mu_\mathcal{M}(\dots (h\mu_\mathcal{M}(i\otimes \id)\otimes \id)\dots \otimes \id)\otimes\id )\\
        &-(-1)^{\frac{N(N-1)}{2}}(-1)^{N+1}m_N(\id^{\otimes N}\otimes  d)\\
        =&\ (-1)^{\frac{(N-1)(N-2)}{2}}m_{N-1}(\id^{\otimes N-2}\otimes \mu_\mathcal{M})\\
        &-(-1)^{\frac{(N-1)(N-2)}{2}}m_{2}(m_{N-1}\otimes \id)\\
        &-p \mu_\mathcal{M}( hd\mu_\mathcal{M}(\dots (h\mu_\mathcal{M}(i\otimes \id)\otimes \id)\dots \otimes \id)\otimes\id )\\
        &+(-1)^{N-2}p \mu_\mathcal{M}( {h \mu_\mathcal{M}(\dots (h\mu_\mathcal{M}}(i\otimes \id)\otimes \id)\dots \otimes \id)\otimes d)\\
        =&\dots\\
        =&-(-1)^{\frac{N(N-1)}{2}}\sum_{r=1}^{N-2}(-1)^{r}m_{N-1}(\id^{\otimes r}\otimes \mu_\mathcal{M}\otimes \id^{\otimes N-r-2})\\
        &-(-1)^{\frac{N(N-1)}{2}}\sum_{t=1}^{N-2} (-1)^{(N-t)t} m_{t+1}(m_{N-t}\otimes \id^{\otimes t})\\
        &-(-1)^{\frac{N(N-1)}{2}}(-1)^{N+1}m_Nd.
    \end{align*}
    This shows that $m_k$ defines a $\mathcal{A}_\infty$-structure on $\mathcal{N}_*$.
    
    Analogous computations show that the above formulas for $p_k$ and $i_k$ define morphism of $\mathcal{A}_\infty$-modules.
\end{proof}

\begin{Prop}\label{prop:inverting infty-quasi}
    Let $(A_*,d_A)$ be a (strict) differential graded algebra. Let $(\mathcal{M}_*,d_\mathcal{M})$ and $(\mathcal{N}_*,d_\mathcal{N})$ be strict right $A_*$-modules with an $\infty$-quasi-isomorphism $\boldsymbol{f}\colon \mathcal{M}_*\to \mathcal{N}_*$.
    
    If $(i^\mathcal{M},p^\mathcal{M},h^\mathcal{M})$ and $(i^\mathcal{N},p^\mathcal{N},h^\mathcal{N})$ exhibit $(\mathcal{M}_*,d_\mathcal{M})$ and $(\mathcal{N}_*,d_\mathcal{N})$ as a homotopy retract of $(H_*(\mathcal{M}),0)$ and $(H_*(\mathcal{N}),0)$, respectively, then there exists an $\infty$-quasi-isomorphism $\boldsymbol{g}\colon \mathcal{N}_*\to \mathcal{M}_*$ such that $H_*(g_1)=H_*(f_1)^{-1}$.
\end{Prop}
\begin{proof}
    By Proposition \ref{prop:htt for Ainfty modules}, we can define the following composition of morphisms of $\mathcal{A}_\infty$-modules as $\boldsymbol{j}$:
    \begin{center}
        \begin{tikzcd}
            H_*(\mathcal{M}) \ar[r, "\boldsymbol{i}^\mathcal{M}"] & \mathcal{M}_* \ar[r, "\boldsymbol{f}"] &  \mathcal{N}_* \ar[r, "\boldsymbol{p}^\mathcal{N}"] & H_*(\mathcal{N}).
        \end{tikzcd}
    \end{center}
    This morphism is such that $j_1=H_*(f_1)$ is an isomorphism. Therefore, we can invert it by Proposition \ref{prop:inverting infty-iso}. We now define $\boldsymbol{g}$ as the composition
    \begin{center}
        \begin{tikzcd}
            \mathcal{N}_*  \ar[r, "\boldsymbol{p}^\mathcal{N}"] & H_*(\mathcal{N}) \ar[r, "\boldsymbol{j}^{-1}"] &H_*(\mathcal{M}) \ar[r, "\boldsymbol{i}^\mathcal{M}"] & \mathcal{M}_*.
        \end{tikzcd}
    \end{center}
\end{proof}

\begin{Cor}\label{cor:inverting infty-quasi over field}
    Let $(A_*,d_A)$ be a (strict) differential graded algebra over a field. Let $(\mathcal{M}_*,d_\mathcal{M})$ and $(\mathcal{N}_*,d_\mathcal{N})$ be strict right $A_*$-modules with an $\infty$-quasi-isomorphism $\boldsymbol{f}\colon \mathcal{M}_*\to \mathcal{N}_*$.
    
    There exists an $\infty$-quasi-isomorphism $\boldsymbol{g}\colon \mathcal{N}_*\to \mathcal{M}_*$ such that $H_*(g_1)=H_*(f_1)^{-1}$.
\end{Cor}
\begin{proof}
    By \cite[Lemma 9.4.4]{loday2012algebraic} any complex over a field admits its homology as a homotopy retract. Together with Proposition \ref{prop:inverting infty-quasi}, this shows the corollary.
\end{proof}

\section{Homotopies of Relative Embeddings}\label{app:hpt of rel emb}
We construct a homotopy of the space of smooth, pointed, orientations-preserving, relative embeddings $(\R^n,B_1(0)^c)\hookrightarrow (\R^n,B_1(0)^c)$ that exhibits the group $\SO(n)$ as a deformation retract of the space of all such embeddings.\\

We denote $\Emb_*^+((\R^n,B_1(0)^c),(\R^n,B_1(0)^c))$ as the space of all smooth embeddings $\varphi\colon \R^n\hookrightarrow \R^n$ such that
\begin{enumerate}[(i)]
    \item $\varphi(0)=0$;
    \item $\varphi(B_1(0)^c)\subseteq \varphi(B_1(0)^c)$;
    \item $\varphi$ is orientation-preserving.
\end{enumerate}
The group $\SO(n)$ is a subspace of $\Emb_*^+((\R^n,B_1(0)^c),(\R^n,B_1(0)^c))$. We define a homotopy that retracts $\Emb_*^+((\R^n,B_1(0)^c),(\R^n,B_1(0)^c))$ to $\SO(n)$.

The naive strategy is to first homotope $\Emb_*^+((\R^n,B_1(0)^c),(\R^n,B_1(0)^c))$ to $\GL^+(n)$ using a radial homotopy and then to homotope $\GL^+(n)$ to $\SO(n)$. This however fails because $A\in \GL^+(n)$ does not send $B_1(0)^c$ to itself in general.

Our strategy is to flip the two steps:
\begin{enumerate}[(1)]
    \item\label{item:hpt of rel emb 1} We homotope $\Emb_*^+((\R^n,B_1(0)^c),(\R^n,B_1(0)^c))$ to $\Emb^\SO((\R^n,B_1(0)^c),(\R^n,B_1(0)^c))$, the subspace of all relative embeddings which have derivative in $\SO(n)$ at $0$.
    \item\label{item:hpt of rel emb 2} We homotope $\Emb^\SO((\R^n,B_1(0)^c),(\R^n,B_1(0)^c))$ to $\SO(n)$ using a radial homotopy.
\end{enumerate}
For step \eqref{item:hpt of rel emb 1}, we recall that any matrix $A\in \GL^+(n)$ can be uniquely written as $QP$ for $Q\in \SO(n)$ and $P\in \SPD(n)$, where $\SPD(n)$ denotes the space of symmetric, positive definite matrices. For $\varphi \in \Emb_*^+((\R^n,B_1(0)^c),(\R^n,B_1(0)^c))$, we denote $Q(\varphi)\in \SO(n)$ and $P(\varphi)\in \SPD(n)$ for the matrices such that $D_0\varphi=Q(\varphi)P(\varphi)$. Step \eqref{item:hpt of rel emb 1} thus is about homotoping $P(\varphi)$ to the identity.

\subsection{Homotopies of Positive Definite Matrices}
It is well-known that the space $\SPD(n)$ is contractible. We construct a particularly nice contraction.

\begin{Lem}\label{lem:H1 in app}
    For all $n>0$, there exists a map
    \begin{align*}
        H^{(n)}\colon I\times \mathrm{SPD}(n)&\to \mathrm{SPD}(n)
    \end{align*}
    such that, for all $A\in \mathrm{SPD}(n)$, the map
    \begin{align*}
        f_A\colon \R^n &\to \R^n,\\
        x&\mapsto \begin{cases} 
            H^{(n)}(\|x\|,A)\cdot x & \text{if } \|x\|\leq 1;\\
            x &\text{if }\|x\|\geq 1;
        \end{cases}
    \end{align*}
    defines a diffeomorphism that sends $0$ to $0$ with $D_0f_A=A$. In particular, $H^{(n)}$ is a homotopy between the identity and the map that is constantly the identity matrix.
\end{Lem}
\begin{proof}
    We first prove the statement for $n=1$. Then $\mathrm{SPD}(1)=\R^{>0}$. We define a symmetric, smooth function $f_\lambda$ directly.

    By defining $f_{\lambda}:=f_{\lambda^{-1}}^{-1}$ for $\lambda> 1$, we may assume that $\lambda\leq 1$. We then define
    \begin{align*}
        f_\lambda(x)=x+(1-\lambda) x\phi(x)
    \end{align*}
    for a bump function $\phi\colon \R\to \R$ such that 
    \begin{enumerate}[(i)]
        \item $\phi(0)=1$;
        \item $\phi'(0)=0$;
        \item $\phi\equiv 0$ outside $(-1,1)$;
        \item $\phi(x)+x\phi'(x)\geq 0$.
    \end{enumerate}
    For example, we can take $\phi(x):=e^{1-\frac{1}{1-x^2}}$. We then have for $x\in [-1,1]$
    \begin{align*}
        f'_\lambda(x)=1+(1-\lambda)(\phi(x)+x\phi'(x)).
    \end{align*}
    This is positive and thus $f_\lambda$ is a diffeomorphism. Moreover $f'_\lambda(0)=1+(1-\lambda)(1+0)=\lambda$ as desired.

    Then the homotopy 
    \begin{align*}
        H^{(1)}\colon I\times \R^{>0}&\to \R^{>0},\\
        (t,\lambda)&\mapsto \begin{cases}f_\lambda(t)/t &\text{if } t>0;\\
            \lambda &\text{if }t=0;
        \end{cases}
    \end{align*}
    is a homotopy as in the statement.
    
    We now consider the case $n>1$. Then $A\in \mathrm{SPD}(n)$ is diagonalizable with positive eigenvalues $\lambda_1,\dots,\lambda_n$. We can thus write
    \begin{align*}
        A=Q^{-1}\mathrm{Diag}(\lambda_1,\dots,\lambda_n)Q
    \end{align*}
    for some orthogonal matrix $Q$. We then define
    \begin{align*}
        H^{(n)}(t,A)=Q^{-1}\mathrm{Diag}(H^{(1)}(t,\lambda_1),\dots,H^{(1)}(t,\lambda_n))Q.
    \end{align*}
    This is independent of our choice of $Q$ because on each eigenspace of an eigenvalue $\lambda_i$, we act at time $t$ via $H^{(1)}(t,\lambda_i)$. 
    
    We thus have
    \begin{align*}
        f_A(x)=Q^{-1}\cdot \begin{pmatrix}
            f_{\lambda_1}(y_1)\\\vdots\\ f_{\lambda_n}(y_n)
        \end{pmatrix}
    \end{align*}
    where $y=Q\cdot x$.

    It then follows that $f_A$ is a diffeomorphism that sends $0$ to $0$ and $D_0f_A$ from the $1$-dimensional case.
\end{proof}

We can thus construct the homotopy 
\begin{align*}
    G\colon I\times \Emb_*^+((\R^n,B_1(0)^c),(\R^n,B_1(0)^c))&\to \Emb_*^+((\R^n,B_1(0)^c),(\R^n,B_1(0)^c)),\\
    (t,\varphi)&\mapsto \varphi\circ f_{H^{(n)}(t,P(\varphi)^{-1})}.
\end{align*}
We note that:
\begin{align*}
    D_0(G(0,\varphi))=D_0\left(\varphi\circ f_{H^{(n)}(t,P(\varphi)^{-1})}\right)=D_0(\varphi)\cdot P(\varphi)^{-1}=Q(\varphi)\in \SO(n).
\end{align*}
This finishes step \eqref{item:hpt of rel emb 1} of our strategy.

\subsection{Relative Radial Homotopy}
For step \eqref{item:hpt of rel emb 2}, we define a variation of a radial homotopy on $\Emb^\SO((\R^n,B_1(0)^c),(\R^n,B_1(0)^c))$.

\begin{Lem}\label{lem:H2 in app}
    There exists a homotopy
    \begin{align*}
        F\colon I\times \Emb^\SO((\R^n,B_1(0)^c),(\R^n,B_1(0)^c))\to \Emb^\SO((\R^n,B_1(0)^c),(\R^n,B_1(0)^c))
    \end{align*}
    such that
    \begin{enumerate}[(i)]
        \item $F(0,\varphi)=D_0(\varphi)$;
        \item $F(1,\varphi)=\varphi$.
    \end{enumerate}
\end{Lem}
\begin{proof}
    For $\varphi\in \Emb_*^{\SO}((\R^n,B_1(0)^c),((\R^n,B_1(0)^c))$ and $t\in I$, we define
    \begin{align*}
        r_t:=r_t(\varphi):=\sup_{\|v\|\leq t} \|\varphi^{-1}(v)\|=\max_{\|v\|=t}\|\varphi^{-1}(v)\|
    \end{align*}
    as in Figure \ref{fig:rt}.
    \begin{figure}[H]
        \centering
        \includegraphics[width=0.4\linewidth]{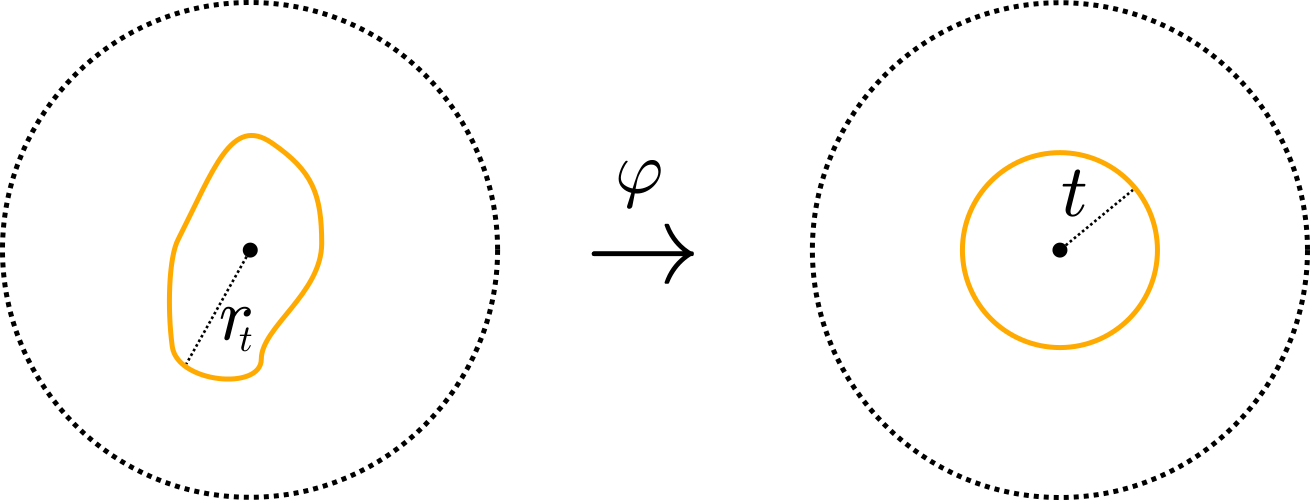}
        \caption{The radius $r_t$ is the smallest number such that $x$ with $\|x\|\geq r_t$ gets sent to outside of $B_t(0)$ under $\varphi$.}
        \label{fig:rt}
    \end{figure}
    We note that 
    \begin{align}\label{eq:r_t/t}
        \lim_{t\to 0}\frac{r_t(\varphi)}{t}=\max_{\|v\|=1}\|D_0(\varphi^{-1})\cdot v\|=1
    \end{align}
    because $D_0(\varphi)\in \SO(n)$ by assumption and thus $D_0(\varphi^{-1})\in \SO(n)$.
    
    We now define
    \begin{align*}
        F(t,\varphi)(x):= \begin{cases}
            \frac{\varphi((r_t(1-t^2)+t^2)x)}{t}&\text{if }t>0;\\
            D_0(\varphi)\cdot x &\text{if }t=0.
        \end{cases}
    \end{align*}
    By \eqref{eq:r_t/t}, this is continuous in $t$. Moreover, we have $F(1,\varphi)=\varphi$.
    
    It remains to check that $F(t,\varphi)$ sends $B_1(0)^c$ to itself. For $t=0$, this follows from the fact that $D_0(\varphi)\in \SO(n)$. We thus assume $t>0$.
    
    We note for $y$ with $\|y\|\geq r_t$, we have that $\varphi(y)\notin B_t(0)$ because 
    \begin{align*}
        r_t=\sup_{\|v\|<t}\|\varphi^{-1}(v)\|
    \end{align*}
    which is a supremum that is not realised by any $v$ with $\|v\|<t$. Therefore, $\varphi/t$ sends $B_{r_t}(0)^c$ to $B_1^c(0)$. Finally $r_t\leq 1$ shows that for $x\notin B^c_1(0)$, we have $(r_t(1-t^2)+t^2)x\in B_{r_t}^c(0)$ and thus $F(t,\varphi)$ sends $B_1(0)^c$ to itself.
\end{proof}

Combining the homotopies from Lemma \ref{lem:H1 in app} and Lemma \ref{lem:H2 in app}, we get the following homotopy:

\begin{Prop}\label{prop:hpt of rel emb}
    There exists a homotopy
    \begin{align*}
        H\colon I\times \Emb_*^+((\R^n,B_1(0)^c),(\R^n,B_1(0)^c))\to \Emb_*^+((\R^n,B_1(0)^c),(\R^n,B_1(0)^c))
    \end{align*}
    such that for $\varphi\in \Emb_*^+((\R^n,B_1(0)^c),(\R^n,B_1(0)^c))$, $Q\in \SO(n)$ and $t\in I$, it holds
    \begin{enumerate}[(i)]
        \item\label{item:H in app 1} $H(0,\varphi)=Q(\varphi)\in \SO(n)$;
        \item\label{item:H in app 2} $H(1,\varphi)=\varphi$;
        \item\label{item:H in app 3} $H(t,Q\cdot \varphi)=Q\cdot H(t,\varphi)$
    \end{enumerate}
\end{Prop}
\begin{proof}
    We define 
    \begin{align*}
        H(t,\varphi):=\begin{cases}
            F(2t,\varphi) &\text{if } t\leq \frac{1}{2};\\
            G(2t-1,F(1,\varphi)) &\text{if }t \geq \frac{1}{2}.
        \end{cases}
    \end{align*}
    Conditions \eqref{item:H in app 1} and \eqref{item:H in app 2} follow from Lemma \ref{lem:H1 in app} and Lemma \ref{lem:H2 in app}. It remains to check Condition \eqref{item:H in app 3}. We compute for $Q\in \SO(n)$:
    \begin{align*}
        G(t,Q\cdot \varphi)(x)&=Q\cdot \varphi\left(f_{H^{(n)}(t,P(Q\cdot \varphi)^{-1})}(x)\right)=Q\cdot \varphi\left(f_{H^{(n)}(t,P( \varphi)^{-1})}(x)\right)=Q\cdot G(t,\varphi)
    \end{align*}
    and
    \begin{align*}
        F(t,Q\cdot \varphi)=\left.\begin{cases}
            \frac{Q\cdot \varphi(r_tx)}{t}&\text{if }t>0;\\
            D_0(Q\cdot \varphi)\cdot x=Q\cdot D_0(\varphi)\cdot x &\text{if }t=0;
        \end{cases}\right\}=Q\cdot F(t,\varphi).
    \end{align*}
    This shows that $F$ and $G$ commute with multiplication by $Q$ and thus so does their composition $H$.
\end{proof}

\normalem
\bibliographystyle{alpha}
\bibliography{ref.bib}

\end{document}